\documentclass[ejs]{imsart}

\usepackage{amsthm,amsmath}
\RequirePackage[numbers]{natbib}
\RequirePackage[linkcolor=blue,citecolor=blue,urlcolor=black]{hyperref}

\startlocaldefs

\usepackage[utf8]{inputenc}
\usepackage[english]{babel}
\usepackage{amsthm,amsmath,mathtools,amsfonts,amssymb}

\usepackage{blindtext}

\usepackage{caption}
\usepackage{subcaption}

\usepackage{graphicx,float}
\usepackage{enumerate}

\usepackage{bbm}

\numberwithin{equation}{section}

\usepackage{dirtytalk}

\usepackage{url}

\usepackage[nameinlink,capitalize,noabbrev]{cleveref}

\usepackage{nicefrac}

\newcounter{algsubstate}

\def\ie{{\em i.e.,~}}
\def\st{{\em s.t.~}}
\def\eg{{\em e.g.,~}}

\newcommand{\wrt}{w.r.t.}

\DeclareMathOperator*{\argmin}{argmin}

\DeclareMathOperator*{\id}{I} 
\newcommand{\zero}{\ensuremath{\mathbf{0}}}

\DeclareMathOperator*{\vect}{vec} 
\DeclareMathOperator*{\rank}{rank} %

\DeclareMathOperator*{\kl}{KL}
\DeclareMathOperator*{\tkl}{KL^{\otimes n}}
\DeclareMathOperator*{\jtkl}{JKL_\rho^{\otimes n}}

\DeclareMathOperator*{\hel}{d^2}
\DeclareMathOperator*{\thel}{d^{2 \otimes n}}
\DeclareMathOperator*{\thell}{d^{\otimes n}}
\DeclareMathOperator*{\entropy}{\mathcal{H}_{d_{\norm{\sup}_\infty}}}
\DeclareMathOperator*{\dsup}{d_{\left\|\sup\right\|_{\infty}}}

\DeclareMathOperator*{\card}{card}

\DeclareMathOperator*{\disc}{disc} 
\DeclareMathOperator*{\Beta}{Beta} 

\DeclareMathOperator*{\sbm}{\cS^{\cB}_{(K,J)}} 

\let\inf\relax 
\DeclareMathOperator*\inf{\vphantom{p}inf}

\DeclareMathOperator*{\adj}{Adj}

\DeclarePairedDelimiter\norm{\lVert}{\rVert}

\newcommand*{\gkdb}{\cG_{K,d,\bfB}} 
\newcommand*{\gkd}{\cG_{K,d}} 
\newcommand*{\gdbk}{\cG_{d,\bfB_k}} 

\newcommand{\UpsilondK}{\Upsilonb_{K,d}} 

\newcommand*{\UpsilondKstar}{\Upsilonb^*_{K,d}} 

\newcommand*{\Upsilondk}{\Upsilonb_{k,d}} 

\newcommand{\tr}{tr}
\newcommand{\nn}{\nonumber} 
\newcommand{\R}{\mathbb{R}} 

\newcommand{\Indi}{\mathbb{I}} 
\newcommand{\Ns}{\mathbb{N}^\star} 
\newcommand{\pen}{\text{pen}} 
\newcommand{\Z}{\mathbb{Z}} 

\newcommand{\fC}{\mathfrak{C}}

\newcommand{\E}[2]{\mathbb{E}_{#1} \left[#2\right]}

\newcommand{\Ep}[1]{\mathbb{E}\left(#1\right)}

\newcommand{\varp}[1]{\text{var}\left(#1\right)}

\newcommand{\covp}[1]{\text{cov}\left(#1\right)}

\newcommand{\cA}{\mathcal{A}}
\newcommand{\cb}{\mathcal{b}}
\newcommand{\cB}{\mathcal{B}}

\newcommand{\cC}{\mathcal{C}}

\newcommand{\cD}{\mathcal{D}}

\newcommand{\cE}{\mathcal{E}}
\newcommand{\cF}{\mathcal{F}}

\newcommand{\cG}{\mathcal{G}}

\newcommand{\cH}{\mathcal{H}}

\newcommand{\cK}{\mathcal{K}}

\newcommand{\cL}{\mathcal{L}}

\newcommand{\cM}{\mathcal{M}}

\newcommand{\cN}{\mathcal{N}}
\newcommand{\cO}{\mathcal{O}}

\newcommand{\cP}{\mathcal{P}}

\newcommand{\cR}{\mathcal{R}}

\newcommand{\cS}{\mathcal{S}}

\newcommand{\cT}{\mathcal{T}}

\newcommand{\cU}{\mathcal{U}}

\newcommand{\cV}{\mathcal{V}}

\newcommand{\cW}{\mathcal{W}}

\newcommand{\cX}{\mathcal{X}}

\newcommand{\cY}{\mathcal{Y}}

\newcommand{\bfA}{\mathbf{A}}
\newcommand{\bfb}{\mathbf{b}}
\newcommand{\bfB}{\mathbf{B}}
\newcommand{\bfc}{\mathbf{c}}
\newcommand{\bfC}{\mathbf{C}}

\newcommand{\bfE}{\mathbf{E}}

\newcommand{\bfm}{\mathbf{m}}

\newcommand{\bfP}{\mathbf{P}}

\newcommand{\bfr}{\mathbf{r}}

\newcommand{\bft}{\mathbf{t}}

\newcommand{\bfv}{\mathbf{v}}
\newcommand{\bfV}{\mathbf{V}}
\newcommand{\bfw}{\mathbf{w}}
\newcommand{\bfW}{\mathbf{W}}

\def\alphab      {{\sbmm{\alpha}}\XS}

\def\Gammab    {{\sbmm{\Gamma}}\XS}

\def\thetab      {{\sbmm{\theta}}\XS}      \def\Thetab    {{\sbmm{\Theta}}\XS}

    \def\Lambdab   {{\sbmm{\Lambda}}\XS}
\def\mub         {{\sbmm{\mu}}\XS}
\def\nub         {{\sbmm{\nu}}\XS}
\def\xib         {{\sbmm{\xi}}\XS}                 
\def\pib         {{\sbmm{\pi}}\XS}                 \def\Pib        {{\sbmm{\Pi}}\XS}

      \def\Sigmab    {{\sbmm{\Sigma}}\XS}

\def\psib        {{\sbmm{\psi}}\XS}        \def\Psib       {{\sbmm{\Psi}}\XS}
\def\omegab      {{\sbmm{\omega}}\XS}      \def\Omegab    {{\sbmm{\Omega}}\XS}

\def\upsilonb    {{\sbmm{\upsilon}}\XS}   \def\Upsilonb  {{\sbmm{\Upsilon}}\XS}

\def\upsilonb       {{\sbmm{\upsilon}}\XS}

\newcommand{\Xv}{\ensuremath{\mathbf{X}}} 
\newcommand{\xv}{\ensuremath{\mathbf{x}}} 
\newcommand{\Yv}{\ensuremath{\mathbf{Y}}}
\newcommand{\yv}{\ensuremath{\mathbf{y}}}

\def\XS{\xspace}
\DeclareMathAlphabet{\mathb}{OML}{cmm}{b}{it}
\def\sbm#1{\ensuremath{\mathb{#1}}}
\def\sbmm#1{\ensuremath{\boldsymbol{#1}}}  

\def\Bb{{\sbm{B}}\XS}  
  \def\cb{{\sbm{c}}\XS}

\def\mb{{\sbm{m}}\XS}

\def\Ub{{\sbm{U}}\XS}  
\def\Vb{{\sbm{V}}\XS}  \def\vb{{\sbm{v}}\XS}
  
\def\Xb{{\sbm{X}}\XS}  \def\xb{{\sbm{x}}\XS}
  
\def\Zb{{\sbm{Z}}\XS}  \def\zb{{\sbm{z}}\XS}

\theoremstyle{plain}
\newtheorem{theorem}{Theorem}[section]
\newtheorem{lem}[theorem]{Lemma}

\newtheorem{prop}[theorem]{Proposition}

\newtheorem{assumption}{Assumption}[section]
\theoremstyle{definition}
\newtheorem{definition}[theorem]{Definition}

\newtheorem{example}[theorem]{Example}

\newtheorem{remark}[theorem]{Remark}

\crefname{lem}{Lemma}{Lemmas}
\crefname{theorem}{Theorem}{Theorems}
\crefname{prop}{Proposition}{Propositions}
\theoremstyle{plain}

\endlocaldefs

\begin{document}

\begin{frontmatter}
	
\title{A non-asymptotic approach for model selection via penalization in high-dimensional mixture of experts models\thanksref{t1}}

%
\runtitle{A non-asymptotic model selection in mixture of experts models}

\begin{aug}
	\author{\fnms{TrungTin} \snm{Nguyen}
		\corref{}\thanksref{t2}\ead[label=e1]{trung-tin.nguyen@inria.fr}\ead[label=e2,url]{http://trung-tinnguyen.github.io/}}
	\address{Normandie Univ, UNICAEN, CNRS, LMNO, 14000 Caen, France.\\
		Univ. Grenoble Alpes, Inria, CNRS, Grenoble INP, LJK,\\ Inria Grenoble Rhone-Alpes, 655 av. de l'Europe, 38335 Montbonnot, France.\\
		\printead{e1}\\ \printead{e2}}
	
	\author{\fnms{Hien Duy} \snm{Nguyen}
		\corref{}\ead[label=e3]{h.nguyen7@uq.edu.au}\ead[label=e4,url]{https://hiendn.github.io/}}
	\address{Department of Mathematical and Physical Sciences,\\ La Trobe University, Bundoora Melbourne 3066, Victoria Australia.\\ 
		School of Mathematics and Physics,\\
		University of Queensland, St. Lucia, Queensland, Australia.\\	
		\printead{e3}\\ \printead{e4}}
	
	\author{\fnms{Faicel} \snm{Chamroukhi}
		\corref{}\ead[label=e5]{faicel.chamroukhi@unicaen.fr}\ead[label=e6,url]{https://chamroukhi.com/}}
	\address{Normandie Univ, UNICAEN, CNRS, LMNO, 14000 Caen, France.\\ \printead{e5}\\ \printead{e6}}
	
	\author{\fnms{Florence} \snm{Forbes}
		\corref{}\ead[label=e7]{florence.forbes@inria.fr}\ead[label=e8,url]{http://mistis.inrialpes.fr/people/forbes/}}
	\address{Univ. Grenoble Alpes, Inria, CNRS, Grenoble INP, LJK,\\ Inria Grenoble Rhone-Alpes, 655 av. de l'Europe, 38335 Montbonnot, France.\\ \printead{e7}\\ \printead{e8}}
	
	%
	\thankstext{t1}{This work is partially supported by the French Ministry of Higher Education and Research (MESRI), French National Research Agency (ANR) grant SMILES ANR-18-CE40-0014, Australian Research Council grant number DP180101192, and the Inria LANDER project.}
	\thankstext{t2}{Corresponding author.}
	
	\runauthor{T. Nguyen et al.}
	
\end{aug}
	
	\begin{abstract}
		Mixture of experts (MoE) are a popular class of statistical and machine learning models that have gained attention over the years due to their flexibility and efficiency. 
	In this work, we consider Gaussian-gated localized MoE (GLoME) and block-diagonal covariance localized MoE (BLoME) regression models to present nonlinear relationships in heterogeneous data with potential hidden graph-structured interactions between high-dimensional predictors.
		These models pose difficult statistical estimation and model selection questions, both from a computational and theoretical perspective.
			This paper is devoted to the study of the problem of model selection among a collection of GLoME or BLoME models characterized by the number of mixture components, the complexity of Gaussian mean experts, and the hidden block-diagonal structures of the covariance matrices, in a penalized maximum likelihood estimation framework.
			In particular, we establish non-asymptotic risk bounds that take the form of weak oracle inequalities, provided that lower bounds for the penalties hold. The good empirical behavior of our models is then demonstrated on synthetic and real datasets.
		
		%
	\end{abstract}
	
	
	
	
	\tableofcontents
	
\end{frontmatter}


\section{Introduction}






%
\subsection{Mixture of experts models}
Mixture of experts (MoE) models, originally introduced as neural network architectures in \cite{jacobs1991adaptive,jordan1994hierarchical}, are flexible models that generalize the classical finite mixture models as well as finite mixtures of regression models \citep[Sec. 5.13]{mclachlan2000finite}. Their flexibility comes from allowing the mixture weights (or the gating functions) to depend on the explanatory variables, along with the  component densities (or experts). In the context of regression, MoE models with Gaussian experts and softmax or normalized Gaussian gating functions are the most popular choices and are powerful tools for modeling more complex nonlinear relationships between outputs (responses) and inputs (predictors) that arise from different subpopulations. 
The popularity of these conditional mixture density models arise largely due to their universal approximation properties, which have been studied in \citep{jiang1999hierarchical, norets2010approximation, nguyen2016universal, ho2019convergence,nguyen2019approximation,nguyen2020approximationMoE}, and which improve upon approximation capabilities of unconditional finite mixture models, as studied in \citep{genovese2000rates,rakhlin2005risk, nguyen2013convergence, ho2016convergence, ho2016strong, nguyen2020approximation,nguyen2021approximationLebesgue}.
%
%
Detailed reviews on the practical and theoretical aspects of MoE models can be found in \cite{yuksel2012twenty,masoudnia2014mixture,nguyen2018practical,nguyen_model_2021}.


	In this paper, we wish to investigate MoE models with Gaussian experts and normalized Gaussian gating functions for clustering and regression, first introduced by \cite{xu1995alternative}, which extended the original MoE models of \cite{jacobs1991adaptive}.
	From hereon in, we refer to these models as \emph{Gaussian-gated localized MoE} (GLoME) models and \emph{block-diagonal covariance for localized MoE} (BLoME) models, and we provide their precise definitions in \cref{sec_GLoME_Model}. 
	Furthermore, the original MoE models with softmax gating functions will be referred to as \emph{linear-combination-of-bounded-functions softmax-gated MoE} (LinBoSGaME) regression models. In particular, we simply refer to affine instances of LinBoSGaME models as \emph{softmax-gated MoE} (SGaME) regression models.
	It is worth pointing out that the BLoME models generalize GLoME models by utilizing a parsimonious covariance structure, via block-diagonal structures for covariance matrices in the Gaussian experts. It is also interesting to point out that supervised \emph{Gaussian locally-linear mapping} (GLLiM) and \emph{block-diagonal covariance for Gaussian locally-linear mapping} (BLLiM) models in \cite{deleforge2015high} and \cite{devijver2017nonlinear} are affine instances of GLoME and BLoME models, respectively, where the latter pair consider general linear combination of bounded functions instead of affine for mean functions of Gaussian experts, in the prior pair.

%

%
	One of the main disadvantages of LinBoSGaME models is the difficulty of applying an expectation–maximization (EM) algorithm \citep{dempster1977maximum,mclachlan1997algorithm}, which requires an internal iterative numerical optimization procedure (\eg iteratively-reweighted least squares, Newton-Raphson algorithm) to update the softmax parameters. To overcome this problem, we instead use the Gaussian gating network that enables us to link BLoME with hierarchical Gaussian mixture model. Then, the maximization with respect to the parameters of the gating network can be solved analytically with the EM algorithm framework, which decreases the computational complexity of the estimation routine. Furthermore, we can then also make use of well-established theoretical results for finite mixture models.

We note here that both GLoME and BLoME models have been thoroughly studied in the statistics and machine learning literatures in many different guises, including localized MoE \citep{ramamurti1996structural,ramamurti1998use,moerland1999classification,bouchard2003localised}, normalized Gaussian networks \citep{sato2000line}, MoE modeling of priors in Bayesian nonparametric regression \citep{norets2014posterior,norets2017adaptive}, cluster-weighted modeling \citep{ingrassia2012local}, 
deep mixture of linear inverse regressions \citep{Lathuiliere2017}, hierarchical Gaussian locally linear mapping structured mixture (HGLLiM) model \citep{tu2019prediction}, multiple-output Gaussian gated mixture of linear experts \citep{nguyen2019approximation}, and approximate Bayesian computation with surrogate posteriors using GLLiM \citep{forbes2021approximate}.

\subsection{Goals}
Our main goal is to learn potentially nonlinear relationships between a multivariate output and a high-dimensional input issued from a heterogeneous population. This involves performing regression, clustering and model selection, simultaneously. While estimation can be performed using standard EM algorithms, it crucially depends and requires data-driven hyperparameter choices, including the number of mixture components (or clusters), the degree of complexity of each Gaussian expert's mean function, and the hidden block-diagonal structures of the covariance matrices.
%
%
%
Recall that hyperparameter choices of data-driven learning algorithms belong to the class of model selection problems, which have attracted a lot of attention in statistics and machine learning over the past 50 years \citep{akaike1974new,mallows1973some,burnham2002model,massart2007concentration,arlot2019minimal}. Specifically, given a set of models, how do we select the one with the lowest possible risk from the data? Note that penalization is one of the main strategies proposed for model selection. It suggests choosing the estimator that minimizes the sum of its empirical risk and some penalty terms corresponding to the fit of the model to the data while avoiding overfitting.


\subsection{Related literature}
In general, model selection for MoE models is often performed using the Akaike information criterion (AIC) \citep{akaike1974new}, the Bayesian information criterion (BIC) \citep{schwarz1978estimating} or the BIC-like
approximation of integrated classification likelihood (ICL-BIC) \cite{biernacki_assessing_2000}. An important limitation of these criteria, however, is that they are only valid asymptotically. This implies that there is no finite-sample guarantee when using  AIC,  BIC or ICL-BIC, to choose between different levels of complexity. Their use in small samples is, therefore, ad hoc.
To overcome such difficulties, \cite{birge2007minimal} proposed a novel approach, called the slope heuristics, supported by a non-asymptotic oracle inequality. This method leads to an optimal data-driven choice of multiplicative constants for the penalties. Recent reviews and practical issues concerning the slope heuristic can be found in \cite{baudry2012slope,arlot2019minimal} and the references given therein.
%

It should be emphasized that a general model selection result, originally established by \citep[Theorem 7.11]{massart2007concentration}, ensures a penalized criterion leads to a good model selection and that the penalty is known only up to multiplicative constants and is proportional to the dimensions of models. In particular, these multiplicative constants can be calibrated by the slope heuristic approach in a finite-sample setting.
Then, in the spirit of the methods based on concentration inequalities developed in \cite{massart2007concentration,Massart:2011aa,cohen2011conditional,cohen2013partition}, a number of finite-sample oracle type inequalities have been established for the least absolute shrinkage and selection operator (LASSO) \citep{tibshirani1996regression} and general penalized maximum likelihood estimators (PMLE). These results include the works for high-dimensional Gaussian graphical models \citep{devijver2018block}, Gaussian mixture model selection \citep{maugis2011non,maugis2011data,maugis2013adaptive}, finite mixture regression models \citep{Meynet:2013aa, devijver2015l1, devijver2015finite,devijver2017model,devijver2017joint}, and LinBoSGaME models, outside the high-dimensional setting \citep{montuelle2014mixture}. 


%
\subsection{Main contributions}
%

%
To the best of our knowledge, no attempt has been made in the literature to develop a finite-sample oracle inequality for the framework of MoE regression models for high-dimensional data. In this paper, our first original contribution is to provide finite-sample oracle inequalities, \cref{thm_Oracle_Inequality_GLoME,thm_Oracle_Inequality_BLoME} for high-dimensional GLoME and BLoME models, respectively.  
More specifically, we establish non-asymptotic risk bounds that take the form of weak oracle inequalities, provided that the lower bounds on the penalties are true, in high-dimensional GLoME and BLoME models, based on an inverse regression strategy and block-diagonal structures of the Gaussian covariance matrices experts.
Unlike traditional criteria such as AIC, BIC, or ICL-BIC, which are based on asymptotic theory or a Bayesian approach, our contributed non-asymptotic risk bounds allow the number $n$ of observations to be fixed while the dimensionality and cardinality of the models, characterized by the number of covariates and the size of the response, are allowed to grow with respect to $n$ and can be much larger than $n$.

In particular, our oracle inequalities show that the Jensen--Kullback--Leibler loss performance of our PMLEs is comparable to that of oracle models if we take sufficiently large constants in front of the penalties, whose shapes are only known up to multiplicative constants and proportional to the dimensions of the models.
These theoretical justifications for the shapes of the penalties motivate us to make use of the slope heuristic criterion to select several hyperparameters, including the number of mixture components,  the degree of polynomial mean functions, and the potential hidden block-diagonal structures of the covariance matrices of the multivariate predictors.

Moreover, in \cref{thm_Oracle_Inequality_GLoME}, we extends a corollary of \citep[Theorem 1]{montuelle2014mixture}, which can be verified via Lemma 1 from \cite{nguyen2020approximationMoE}, which makes explicit the relationship between softmax and Gaussian gating classes.


Another significant contribution is the numerical experiments for simulated and real data sets in \cref{sec_numerical_Experiment} that support our theoretical results, and the statistical study of non-asymptotic model selection in GLoME models.
To the best of our knowledge, instead of using classical asymptotic approaches for model selection such as AIC, BIC and ICL-BIC, we are the first to illustrate that the slope heuristic works well for MoEs with Gaussian experts and normalized Gaussian gating networks such as GLLiM, BLLiM, GLoME, and BLoME models.  
Note that our main objective here is to investigate how well the empirical tensorized Kullback--Leibler divergence between the true model ($s_0^*$) and the selected model $\widehat{s}_{\widehat{\bfm}}^*$ follows the finite-sample oracle inequality of \cref{thm_Oracle_Inequality_GLoME}, as well as the convergence rate $\cO(n^{-1})$ of the error upper bound.
Therefore, we focus on  $1$-dimensional data sets. Beyond the statistical estimation and model selection purposes considered here, the dimensionality reduction capability of GLLiM and BLLiM models in high-dimensional regression data, can be found in \citep[Section 6]{deleforge2015high} and \citep[Sections 3 and 4]{devijver2017nonlinear}, respectively.

Specifically, besides the important theoretical issues regarding the tightness of upper bounds, how to integrate a priori information, and a minimax analysis of our proposed PMLEs, our finite-sample oracle inequalities and the corresponding illustrative numerical experiments help to partially answer the following two important questions raised in the field of MoE regression models: (1) What is the number of mixture components $K$ to choose, given the sample size $n$; and (2) Whether it is better to use a few complex experts or a combination of many simple experts, given the total number of parameters.
Note that such problems are also considered in the work of \citep[Proposition 1]{mendes2012convergence}, where the authors provided some qualitative guidance and only suggested a practical method for choosing $K$ and $d$, involving a complexity penalty or cross-validation. Furthermore, their model is only non-regularized maximum-likelihood estimation, and thus is not suitable for the high-dimensional setting.

Lastly, we emphasize that although the finite-sample oracle inequalities compare performances of our estimators with the best model in the collection, they also allow us to approximate a rich class of conditional densities if we take enough clusters and/or a high enough degree $d$ of polynomials in the mean of the functions of Gaussian experts, in the context of mixture of Gaussian experts \citep{jiang1999hierarchical,mendes2012convergence,nguyen2016universal,ho2019convergence,nguyen2020approximationMoE}. This leads to the upper bounds of the risks being small.

%
\subsection{Main  challenges}

For the Gaussian gating parameters, the technique for handling the logistic weights in the SGaME models of \cite{montuelle2014mixture} is not directly applicable to the GLoME or BLoME framework, due to the quadratic form of the canonical link. Therefore, we have to propose a \emph{reparameterization trick}\footnote{Note that we only use this nomenclature to perform a change of variables of the Gaussian  gating parameters space of GLoME models via the logistic weights of SGaME models. This reparameterization trick does not stand for the well-known one of variational autoencoders (VAEs) in deep learning literature (see \cite{kingma2013auto}, for more details).}  to bound the metric entropy of the Gaussian gating parameters space; see \eqref{eq_reparameterizationTrick_intro} and \cref{lem_Bracketing_Entropy_Gates_intro_mine} for more details.

To work with conditional density estimation in the BLLiM and BLoME models, it is natural to make use of a general conditional density model selection theorem from \cite{cohen2011conditional,cohen2013partition}, see also \cref{thm_weakOracleInequality_PennecEL_intro} for its adaption to our context.
However, it is worth mentioning that since the collection of models constructed by the BLLiM model \cite{devijver2017nonlinear} or by some appropriate procedures for BLoME models in practice is usually random, we have to utilize a model selection theorem for MLE among a random subcollection (cf. \citep[Theorem 5.1]{devijver2015finite} and \citep[Theorem 7.3]{devijver2018block}). 

%
In particular, our collections of BLoME models must satisfy certain regularity assumptions,  see \cref{proofLemma_nguyen2021nonBLoME_intro} for details, which are not easy to verify  due to the complexity of BLoME models and technical reasons.  For BLoME models, the main difficulty in proving our oracle inequality lies in bounding the bracketing entropy of the Gaussian gating functions and Gaussian experts with block-diagonal covariance matrices. 
To solve the first problem, we need to use our previous trick of reparametrizing the Gaussian gating parameter space. For the second one, based on some ideas of Gaussian mixture models from \cite{genovese2000rates,maugis2011non},  we provide a novel extension for standard MoE models with Gaussian gating networks. Note that our important contributions also extend the recent result on block-diagonal covariance matrices in \cite{devijver2018block}, which is only developed for Gaussian graphical models.

In addition, there are two possible ways to decompose the bracketing entropy of collection models based on different distances. As mentioning in Appendix B.2.1 from \cite{montuelle2014mixture}, their decomposition boils down to assuming that the predictor domain is bounded.  And this limiting assumption can then be relaxed by using a smaller distance, namely a tensorized extension of the Hellinger distance $\thell$ in \eqref{eq_tensorized_Hellinger}, but bounding the corresponding bracketing entropy becomes much more difficult. 
Nevertheless, it is important to point out that we are able to weaken this limiting assumption by using the smaller distance: $\thell$, for the bracketing entropy of our collection of BLoME models, see \cref{lem_bracketingEntropyDecomposition2_intro} for more details.  This reinforces our original contributions concerning the control of bracketing entropy of BLoME models.

%

\subsection{Organization}
%
%
%

The rest of this paper is organized as follows. In \cref{notationAndFramework}, we present the notation and framework for GLoME, and BLoME models and their special cases, GLLiM, and BLLiM models. In \cref{mainResult}, we state the main results of this paper: finite-sample oracle inequalities satisfied by the PMLEs. \cref{proofMainTheorem} is devoted to the proofs and main mathematical challenges in establishing our new finite-sample oracle inequalities. We experimentally evaluate our new results in simulated and real datasets in \cref{sec_numerical_Experiment}. Some conclusions and perspectives are provided in \cref{conclusions}. The proofs of technical lemmas can be found in \cref{proofLemma}. 

\section{Notation and framework}\label{notationAndFramework}
From now on, for any $L \in \Ns$, we write $[L]$ to denote the set $\left\{1,\ldots,L\right\}$.  We are interested in estimating the law of the random variable $\Yv =\left(\Yv_j\right)_{j \in [L]}$, conditionally on $\Xv=\left(\Xv_j\right)_{j \in [D]}$.
Here and subsequently, given any $n \in \Ns$, $\left(\Xv_{[n]},\Yv_{[n]}\right) := \left(\Xv_i,\Yv_i\right)_{ i \in [n]}$ denotes a random sample, and $\xv$ and $\yv$ stands for the observed values of the random variables $\Xv$ and $\Yv$, respectively.
The following assumptions will be needed throughout the paper.  We assume that the covariates $\Xv$ are independent but not necessarily identically distributed. The assumptions on the responses $\Yv$ are stronger: conditional on $\Xv_{[n]} $, $\Yv_{[n]}$ are independent, and each $\Yv$ follows a law with true (but unknown) PDF $s_0\left(\cdot \mid \Xv=\xv\right)$, which is approximated via MoE models.
For a matrix $\bfA$, let $m(\bfA)$ and $M(\bfA)$ be, respectively, the modulus of the smallest and largest eigenvalues of $\bfA$.

%


\subsection{GLoME and BLoME models}\label{sec_GLoME_Model}
Motivated by an inverse regression framework, where the role of predictor and response variables will be exchanged such that $\Yv $ becomes the input and $\Xv$ plays the role of a multivariate output, we consider GLoME and BLoME models with inverse conditional probability density functions (PDFs) of the form \eqref{eq_define_GLoME}. This construction goes back to the work of \cite{li1991sliced,xu1995alternative,deleforge2015high,perthame2018inverse}, and is very useful in a high-dimensional regression context, where typically $D \gg L$. In this way, we define PDF of $\Xv$ conditional on $\Yv=\yv$ as follows:
\begin{align}
	s_{\psib_{K,d}}(\xv \mid \yv) &= \sum_{k=1}^K g_{k}\left(\yv;\omegab\right) \phi_D\left(\xv;\upsilonb_{k,d}(\yv),\Sigmab_k\right),\label{eq_define_GLoME}\\
	g_{k}\left(\yv;\omegab\right)&= \frac{ \pi_k\phi_L\left(\yv;\bfc_k,\Gammab_k\right)  }{\sum_{j=1}^K \pi_j \phi_L\left(\yv;\bfc_j,\Gammab_j\right)}\cdot \label{eq_GLoME_GaussianGating}
\end{align}
Here, $g_k(\cdot; \omegab)$ and $\phi_D\left(\cdot;\upsilonb_{k,d}(\cdot),\Sigmab_k\right)$, $k\in[K]$, $K \in \Ns$, $d \in \Ns$, are called normalized Gaussian gating functions (or Gaussian gating networks) and Gaussian experts, respectively. Furthermore, we decompose the parameters of the model as follows: $\psib_{K,d} = \left(\omegab,\upsilonb_{d},\Sigmab\right) \in \Omegab_K\times \Upsilonb_{K,d}\times \bfV_K =:\Psib_{K,d} $, $\omegab = \left(\pib,\cb,\Gammab\right) \in \left(\Pib_{K-1} \times \bfC_K \times \Vb'_K\right) =: \Omegab_K$, $\pib = \left(\pi_k\right)_{k\in[K]}$, $\cb = \left(\bfc_k\right)_{k\in[K]}$, $\Gammab = \left(\Gammab_k\right)_{k\in[K]}$, $\upsilonb_{d}= \left(\upsilonb_{k,d}\right)_{k\in[K]} \in\Upsilonb_{K,d} $, and $\Sigmab= \left(\Sigmab_k\right)_{k\in[K]} \in \bfV_K$. Note that $\Pib_{K-1} =\left\{ \left(\pi_k\right)_{k\in[K]} \in \left(\R^+\right)^K, \sum_{k=1}^K \pi_k = 1\right\}$ is a $K-1$ dimensional probability simplex, $\bfC_K$ is a set of $K$-tuples of mean vectors of size $L \times 1$, $\Vb'_K$ is a set of $K$-tuples of elements in $\cS_L^{++}$, where $\cS_L^{++}$ denotes the collection of symmetric positive definite matrices on $\R^L$, $\Upsilonb_{K,d}$ is a set of $K$-tuples of mean functions from $\R^L$ to $\R^D$ depending on a degree $d$ (\eg~a degree of polynomials), and $\bfV_K$ is a set containing $K$-tuples from $\cS_D^{++}$. 

\begin{figure}[ht!]
	\centering
		\includegraphics[width =.8\linewidth]{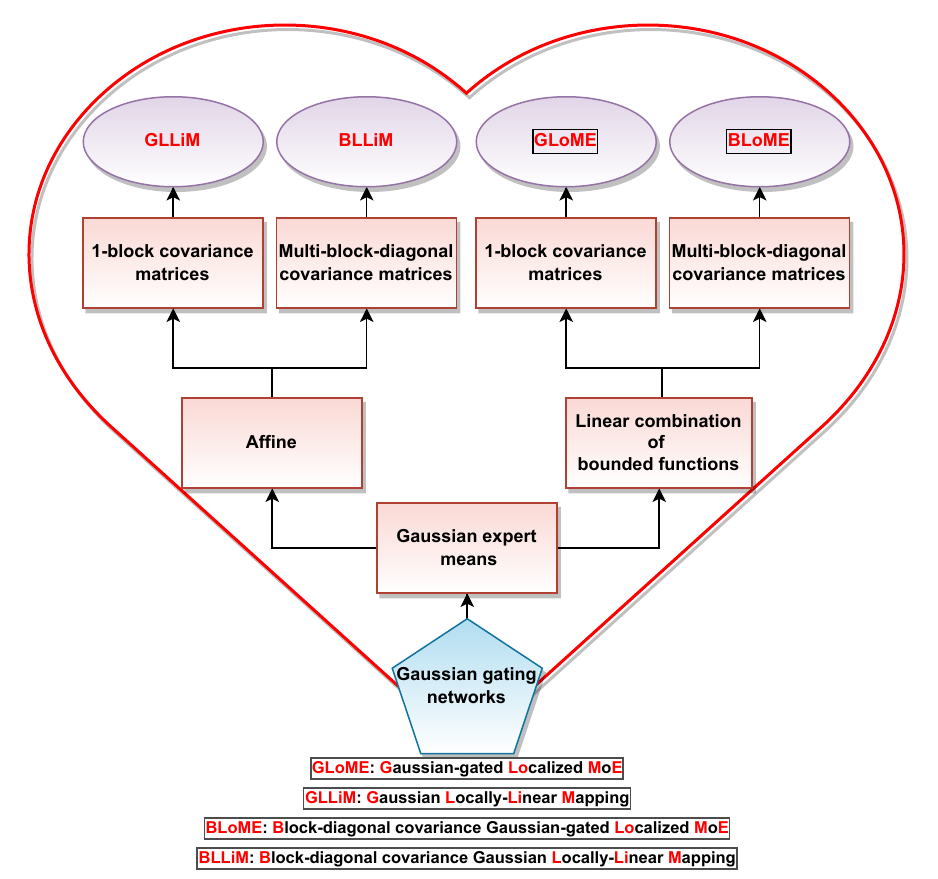}
\caption[A comprehensive classification and nomenclature of MoE models with Gaussian gating networks.]{A comprehensive classification and nomenclature of standard MoE regression models with Gaussian gating networks.}
\label{fig_MoE_GaussianGN}
\end{figure}

Recall that GLLiM and BLLiM models are affine instances of GLoME and BLoME models, and are especially useful for high-dimensional regression data, since there exist link functions between the inverse and forward conditional density; see \cref{fig_MoE_GaussianGN} for comprehensive classification and nomenclature of standard MoE regression models with Gaussian gating networks. Note that the principle of inverse regression is only useful when the functions $\upsilonb_{k,d}(\yv)$ are linear, because there is no explicit way to express the law of $\Yv \mid \Xv$ from that of $\Xv\mid \Yv$ for higher degree of polynomials. However, to have more consistent notations with the previous affine results of GLLiM, BLLiM models from \cite{deleforge2015high,devijver2017nonlinear}, we decide to use the inverse regression frameworks instead of the forward one.

In the BLoME model, we wish to make use of block-diagonal covariance structures by replacing $\Sigmab_k$ and $\bfV_K$ with $\Sigmab_k\left(\bfB_k\right)$ and $\bfV_K\left(\bfB\right)$, defined in \eqref{eq_block_diagonal_structures}, respectively (see, \eg~\citep{devijver2017nonlinear,devijver2018block}).
This block-diagonal structures for covariance matrices are not only used to trade-off between complexity and sparsity but is also motivated by some real applications, where we want to perform prediction on data sets with heterogeneous observations and hidden graph-structured interactions between covariates; for instance, for gene expression data sets in which conditionally on the phenotypic response, genes interact with few other genes only, \ie there are small modules of correlated genes (see	 \citep{devijver2017nonlinear,devijver2018block} for more details).
To be more precise, for $k \in [K]$, we decompose $\Sigmab_{k}\left(\bfB_k\right)$ into $G_k$ blocks, $G_k \in \Ns$, and we denote by $d^{[g]}_k$ the set of variables into the $g$th group, for $g\in \left[G_k\right]$, and by $\card\left(d^{[g]}_k\right)$ the number of variables in the corresponding set. Then, we define $\bfB_k = \left(d^{[g]}_k\right)_{g\in\left[G_k\right]}$ to be a block structure for the cluster $k$, and $\bfB = \left(\bfB_k\right)_{k\in [K]}$ to be the covariate indexes into each group for each cluster. In this way, to construct the block-diagonal covariance matrices, up to a permutation, we make the following definition: $\bfV_K\left(\bfB\right) = \left(\bfV_k\left(\bfB_k\right)\right)_{k\in[K]}$, for every $k\in[K]$,
\begin{align} \label{eq_block_diagonal_structures}
	\bfV_k\left(\bfB_k\right) = \left\{ \Sigmab_{k}\left(\bfB_k\right) \in \cS_D^{++}\left| 
	\begin{array}{l} 
		\Sigmab_{k}\left(\bfB_k\right) = \bfP_k
		\begin{pmatrix}
			\Sigmab_{k}^{[1]}&\zero&\ldots &\zero\\
			\zero&\Sigmab_{k}^{[2]}&\ldots&\zero\\
			\zero&\zero&\ddots&\zero\\
			\zero&\zero&\ldots&\Sigmab_{k}^{[G_k]}
		\end{pmatrix}
		\bfP^{-1}_k,\\
		\Sigmab_{k}^{[g]}\in \cS^{++}_{
			\card\left(d^{[g]}_k\right)}, \forall g \in [G_k]
	\end{array}
	\right.\right\}.
\end{align}
Here, $\bfP_k$ corresponds to the permutation leading to a block-diagonal matrix in cluster $k$. It is worth pointing out that  outside the blocks, all coefficients of the matrix are zeros and we also authorize reordering of the blocks: \eg$\left\{\left(1,3\right);\left(2,4\right)\right\}$ is identical to $\left\{\left(2,4\right);\left(1,3\right)\right\}$, and the permutation inside blocks: \eg the partition of $4$ variables into blocks $\left\{\left(1,3\right);\left(2,4\right)\right\}$ is the same as the partition $\left\{\left(3,1\right);\left(4,2\right)\right\}$.

It is interesting to point out that the BLLiM framework aims to model a sample of high-dimensional regression data issued from a heterogeneous population with hidden graph-structured interaction between covariates. In particular, the BLLiM model is considered as a good candidate for performing model-based clustering and for predicting the response in situations affected by the ``curse of dimensionality" phenomenon, where the number of parameters could be larger than the sample size. Indeed, to deal with high-dimensional regression problems, the BLLiM model is based on an inverse regression strategy, which inverts the role of the high-dimensional predictor and the multivariate response. Therefore, the number of parameters to estimate is drastically reduced.  More precisely, BLLiM utilizes GLLiM, described in \cite{deleforge2015hyper,deleforge2015high}, in conjunction with a block-diagonal structure hypothesis on the residual covariance matrices, to make a trade-off between complexity and sparsity.
This prediction model is fully parametric and highly interpretable. For instance, it can be used for the analysis of transcriptomic data in molecular biology to classify observations or predict phenotypic states; see \eg~\citep{golub1999molecular,nguyen2002tumor,le2008sparse}. Indeed, if the predictor variables are gene expression data measured by microarrays or by the RNA-seq technologies, and the response is a phenotypic variable, 
the BLLiM model not only provides clusters of individuals based on the relation between gene expression data and the phenotype, but also implies a gene regulatory network specific to each cluster of individuals (see \citep{devijver2017nonlinear} for more details).

In order to establish our finite-sample oracle inequalities, \cref{thm_Oracle_Inequality_GLoME,thm_Oracle_Inequality_BLoME}, we need to 
explicitly impose some classical boundedness conditions on the parameter space. 
\subsubsection{Gaussian gating networks}
We shall restrict our study to bounded Gaussian gating parameter vectors. 
Specifically, we assume that there exist deterministic positive constants $a_\pib,A_{\cb},a_{\Gammab},A_{\Gammab}$, such that $\omegab$ belongs to $\widetilde{\Omegab}_K$, where
\begin{align} \label{define_BoundedGatesParameters}
\widetilde{\Omegab}_K &= \bigg\{\omegab \in \Omegab_K : \forall k \in [K], \left\|\bfc_k\right\|_\infty  \le A_{\cb}, \nn\\
& \hspace{2cm}a_{\Gammab} \le m\left({\Gammab}_{k}\right) \le  M\left({\Gammab}_{k}\right) \le A_{\Gammab}, a_\pib \le \pi_k \bigg\}.
\end{align}
%

\subsubsection{Gaussian experts}
Following the same structure for the Gaussian mean experts from \cite{montuelle2014mixture},  the set $\UpsilondK$ will be chosen as a tensor product of compact sets of moderate dimension (\eg~a set of polynomials of degree smaller than $d$, whose coefficients are smaller in absolute values than $T_\Upsilonb$). Then, $\UpsilondK$ is defined as a linear combination of a finite set of bounded functions whose coefficients belong to a compact set. This general setting includes polynomial bases when the covariates are bounded, Fourier bases on an interval, as well as suitably renormalized wavelet dictionaries. More specifically, $ \UpsilondK = \otimes_{k\in[K]}\Upsilonb_{k,d} =: \Upsilonb_{k,d}^K$, where $\Upsilonb_{k,d} = \Upsilonb_{b,d}$, $\forall k \in [K]$, and 
\begin{align} \label{eq_define_meanExperts_linearBounded}
\Upsilonb_{b,d} &= \left\{ \yv \mapsto \left(\sum_{i=1}^{d} \alphab_i^{(j)} \varphi_{\Upsilonb,i}(\yv)\right)_{j\in [D]}=: \left(\upsilonb_{d,j}(\yv)\right)_{j\in [D]}:\norm{\alphab}_\infty \le T_\Upsilonb\right\}.
\end{align}
Here, the subscript $b$ stands for ''bounded functions'', $d \in \Ns, T_\Upsilonb \in \R^+$, and $\left(\varphi_{\Upsilonb,i}\right)_{i \in \left[d\right]}$ is a collection of bounded functions on $\cY$. In particular, we focus on the bounded $\cY$ case and assume that $\cY = [0,1]^L$, without loss of generality. In this case, $\varphi_{\Upsilonb,i}$ can be chosen as monomials with maximum (non-negative) degree $d$: $\yv^{\bfr} = \prod_{l=1}^L \yv^{\bfr_l}_l$. Recall that a multi-index $\bfr = \left(\bfr_l\right)_{l\in[L]}, \bfr_l \in \Ns\cup\left\{0\right\}, \forall l \in [L]$, is an $L$-tuple of nonnegative integers. We define $\left|\bfr\right| = \sum_{l=1}^L \bfr_l$ and the number $|\bfr|$ is called the order or degree of $\yv^{\bfr}$. Then, $ \UpsilondK = \Upsilonb_{p,d}^K$, where
\begin{align}\label{eq_define_meanExperts_polynomial}
\Upsilonb_{p,d} &= \left\{ \yv \mapsto \left(\sum_{|\bfr|=0}^{d} \alphab_\bfr^{(j)} \yv^{\bfr}\right)_{j\in [D]} =: \left(\upsilonb_{d,j}(\yv)\right)_{j\in [D]}:\norm{\alphab}_\infty \le T_\Upsilonb\right\}.
\end{align}
Here, the subscript $p$ stands for ''polynomial". 

For GLoME models, note that any covariance matrix $\Sigmab_k$ can be decomposed into the form $B_k\bfP_k\bfA_k\bfP_k^\top$, such that $B_k = \left|\Sigmab_k\right|^{1/D}$ is a positive scalar corresponding to the volume, $\bfP_k$ is the matrix of eigenvectors of $\Sigmab_k$ and
$\bfA_k$ the diagonal matrix of normalized eigenvalues of $\Sigmab_k$; $B_{-} \in \R^+,B_{+} \in \R^+$,
$\cA\left(\lambda_-,\lambda_+\right)$ is a set of diagonal matrices $\bfA_k$, such that $\left|\bfA_k\right|=1$ and $\forall i \in [D],\lambda_- \le \left(\bfA_k\right)_{i,i} \le \lambda_+$, where $\lambda_-, \lambda_+ \in \R$; and $SO(D)$ is the special orthogonal group of dimension $D$.  In this way, we obtain what is known as the classical covariance matrix sets described by \cite{celeux1995gaussian} for Gaussian parsimonious clustering models, defined by
\begin{align}\label{eq_define_VarianceExperts}
\bfV_K &= \Big\{\left(B_k\bfP_k\bfA_k\bfP_k^\top\right)_{k\in[K]}: \forall k\in[K], B_-\le B_k \le B_+, \nn\\
&\quad \hspace{3cm}  \bfP_k \in SO(D), \bfA_k \in \cA\left(\lambda_-,\lambda_+\right)
\Big\}.
\end{align}

For the block-diagonal covariances of Gaussian experts from BLoME models, we assume that there exist some positive constants $\lambda_m$ and $\lambda_M$ such that, for every $k \in [K]$,
\begin{align}\label{eq_define_blockdiagonal_covariances_boundedness}
	\widetilde{	\bfV}_K\left(\bfB\right) &= \bigg\{ \left(\Sigmab_{k}\left(\bfB_k\right)\right)_{k\in[K]} \in 	\bfV_K\left(\bfB\right): \nn\\
	& \hspace{2cm}0<\lambda_m \le m\left(\Sigmab_{k}\left(\bfB_k\right)\right) \le M\left(\Sigmab_{k}\left(\bfB_k\right)\right) \le \lambda_M  \bigg\}.
\end{align}
Note that this is a typical assumption and is also used in the block-diagonal covariance selection for Gaussian graphical models of \cite{devijver2018block}.


Next, characterizations of GLLiM and BLLiM models are provided in \cref{randomCovModelSelection}.


\subsection{High-dimensional regression via GLLiM and BLLiM models}\label{randomCovModelSelection}
The GLLiM and BLLiM models, as originally introduced in \cite{deleforge2015high,devijver2017nonlinear}, are used to capture the nonlinear relationship between the response and the set of covariates from a high-dimensional regression data,  imposed by a potential hidden graph structured interaction, typically in the case when $D\gg L$, by the $K$ locally affine mappings:
\begin{align}\label{eq_locallyaffine}
\Yv = \sum_{k=1}^K \Indi \left(Z = k\right) \left(\bfA^*_k\Xv +\bfb^*_k + \bfE^*_k\right).
\end{align}
Here, $\Indi$ is an indicator function and $Z$ is a latent variable capturing a cluster relationship, such that $Z=k$ if $\Yv$ originates from cluster $k \in [K]$. Cluster specific affine transformations are defined by matrices $\bfA^*_k \in \R^{L\times D}$ and vectors $\bfb^*_k \in \R^L$. Furthermore, $\bfE^*_k$ are error terms capturing both the reconstruction error due to the local affine approximations and the observation noise in $\R^L$.

Following the common assumption that $\bfE^*_k$ is a zero-mean Gaussian variable with covariance matrix $\Sigmab_k^* \in \R^{L\times L}$, it holds that
\begin{align} \label{eq_conditional_forward}
p\left(\Yv = \yv\mid \Xv = \xv,Z = k;\psib^*_K\right) = \Phi_L\left(\yv;\bfA^*_k\xv +\bfb^*_k,\Sigmab^*_k\right),
\end{align}
where we denote by $\psib^*_K$ the vector of model parameters and $\Phi_L$ is the PDF of a Gaussian distribution of dimension $L$. 
In order to allow the affine transformations to be local, $\Xv$ is defined as a mixture of $K$ Gaussian components as follows:
\begin{align}\label{eq_marginal_forward}
p\left(\Xv=\xv\mid Z = k;\psib^*_K\right) = \Phi_D\left(\xv;\bfc^*_k,\Gammab^*_k\right), p\left(Z=k;\psib^*_k\right) = \pib^*_k,
\end{align}
where $\bfc^*_k \in \R^D, \Gammab^*_k \in \R^{D \times D}$, $\pib^*=\left(\pib^*_k\right)_{k\in[K]} \in \Pib^*_{K-1}$, and $\Pib^*_{K-1}$
is the $K-1$ dimensional probability simplex.
Then, according to formulas for conditional multivariate Gaussian variables and the following hierarchical decomposition 
\begin{align*}
p\left(\Yv = \yv,\Xv = \xv;\psib^*_K\right)
&=  \sum_{k=1}^K \pib^*_k \Phi_D\left(\xv;\bfc^*_k,\Gammab^*_k\right) \Phi_L\left(\yv;\bfA^*_k\xv +\bfb^*_k,\Sigmab^*_k\right),
\end{align*}
we obtain the following \emph{forward conditional density} \citep{deleforge2015high}:
\begin{align} \label{eq_forwardGLLiM}
p\left(\Yv = \yv\mid \Xv = \xv;\psib^*_K\right) = \sum_{k=1}^K\frac{ \pib^*_k \Phi_D\left(\xv;\bfc^*_k,\Gammab^*_k\right)  }{\sum_{j=1}^K \pib^*_j \Phi_D\left(\xv;\bfc^*_j,\Gammab^*_j\right)}\Phi_L\left(\yv;\bfA^*_k\xv +\bfb^*_k,\Sigmab^*_k\right),
\end{align}
where $\psib^*_K =  \left(\pib^*,\thetab^*_K \right) \in \Pi_{K-1} \times \Thetab^*_K =:\Psib^*_K$. Here, $\thetab^*_K=\left( \bfc^*_k,\Gammab^*_k,\bfA^*_k,\bfb^*_k,\Sigmab^*_k\right)_{k\in[K]}$ and 
$$\Thetab^*_K = \left(\R^D\times \cS_D^{++} (\R) \times \R^{L\times D}\times \R^L\times \cS_L^{++} (\R)\right)^K.$$
Without assuming anything on the structure of parameters, the dimension of the model, denoted by $\dim\left(\cdot\right)$, is defined as the total number of parameters that has to be estimated, as follows:
\begin{align*}
\dim\left(\Psib^*_K\right) = K\left(1+D(L+1)+\frac{D(D+1)}{2}+\frac{L(L+1)}{2}+L\right)-1.
\end{align*}

It is worth mentioning that $	\dim\left(\Psib_K\right)$ is typically very large compared to the sample size (see, e.g., \cite{deleforge2015high,devijver2017nonlinear,perthame2018inverse} for more details in their real data sets) whenever $D$ is large and $D\gg L$. Furthermore, it is more realistic to make assumption on the residual covariance matrices $\Sigmab^*_k$ of error vectors $\bfE^*_k$ rather than on $\Gammab^*_k$ (cf. \cite[Section 3]{deleforge2015high}). This justifies the use of the inverse regression trick from \cite{deleforge2015high}, which leads a drastic reduction in the number of parameters to be estimated.

More specifically, in \eqref{eq_forwardGLLiM}, the roles of input and response variables were exchanged such that $\Yv$ becomes the covariates and $\Xv$ plays the role of the multivariate response.  Therefore, its corresponding \emph{inverse conditional density} is defined as a \emph{Gaussian locally-linear mapping} (GLLiM) model, based on the previous hierarchical Gaussian mixture model, as follows:
\begin{align}
p\left(Z=k;\psib_k\right) &= \pi_k,\\
p\left(\Yv=\yv\mid Z = k;\psib_K\right) &= \Phi_L\left(\yv;\bfc_k,\Gammab_k\right),\label{eq_marginal_inverse}\\	
p\left(\Xv = \xv\mid \Yv = \yv,Z = k;\psib_K\right) &= \Phi_D\left(\xv;\bfA_k\yv +\bfb_k,\Sigmab_k\right), \label{eq_conditional_inverse}
\\
p\left(\Xv = \xv\mid \Yv = \yv;\psib_K\right) &= \sum_{k=1}^K\frac{ \pi_k \Phi_L\left(\yv;\bfc_k,\Gammab_k\right)  }{\sum_{j=1}^K \pi_j \Phi_L\left(\yv;\bfc_j,\Gammab_j\right)}\Phi_D\left(\xv;\bfA_k\yv +\bfb_k,\Sigmab_k\right)\label{eq_inverseGLLiM},
\end{align}
where $\Sigmab_k$ is a $D \times D$ covariance structure (usually diagonal in GLLiM models, chosen to reduce the number of parameters) automatically learned from data and $\psib_K$ is the set
of parameters, denoted by $\psib_K =  \left(\pib,\thetab_K \right) \in \Pib_{K-1} \times \Thetab_K =: \Psib_K$. It is important to note that the BLLiM model imposes block-diagonal structures on $\left(\Sigmab_k\right)_{k\in[K]}$ to make a trade-off between complexity and sparsity.

The following interesting feature of both GLLiM and BLLiM models is proved in \cref{proof_1_1_GLLiM}. 

\begin{lem}\label{lem_inverse_regression_mapping}
The parameter $\psib^*_K$ in the forward conditional PDF, defined in \eqref{eq_forwardGLLiM}, can then be deduced from $\psib_K$ in \eqref{eq_inverseGLLiM} via the following one-to-one correspondence:
\begin{align} \label{eq_inverse_regression_mapping}
	\thetab_K=\begin{pmatrix}
		\bfc_k\\
		\Gammab_k \\
		\bfA_k \\
		\bfb_k \\
		\Sigmab_k
	\end{pmatrix}_{k \in [K]}
	\mapsto
	\begin{pmatrix}
		\bfc_k^*\\
		\Gammab_k^* \\
		\bfA^*_k \\
		\bfb^*_k \\
		\Sigmab_k^*
	\end{pmatrix}_{k \in [K]}
	=
	\begin{pmatrix}
		\bfA_k\bfc_k +\bfb_k  \\
		\Sigmab_k+\bfA_k\Gammab_k\bfA_k^\top\\
		\Sigmab_k^*\bfA_k^\top\Sigmab_k^{-1}\\
		\Sigmab_k^* (\Gammab_k^{-1}\bfc_k-\bfA_k^\top\Sigmab_k^{-1}\bfb_k)\\
		\left(\Gammab_k^{-1}+\bfA_k^\top\Sigmab_k^{-1}\bfA_k\right)^{-1}
	\end{pmatrix}_{k \in [K]}
	\in \Thetab^*_K,
\end{align}
with $\pib^* \equiv \pib$.
\end{lem}
\subsection{Collection of GLoME and BLoME models} \label{sec_collectionOfModels}
For GLoME, we only need to choose the degree of polynomials $d$ and the number of components $K$ among finite sets $\cD_\Upsilonb = \left[d_{\max}\right]$ and $\cK = \left[K_{\max}\right]$, respectively, where $d_{\max}\in \Ns$ and $K_{\max}\in \Ns$ may depend on the sample size $n$. 
We wish to estimate the unknown inverse conditional density $s_0$ by inverse conditional densities belonging to the following collection of models $ \left(S_\bfm\right)_{\bfm \in \cM}$. Here, we define $\cM = \left\{\left(K,d\right): K \in \cK,d \in \cD_\Upsilonb\right\}$, and
\begin{align} \label{eq_define_GLoME_Boundedness}
S_\bfm  &=  \Big\{(\xv,\yv) \mapsto
s_{\psib_{K,d}}(\xv\mid  \yv) =:s_\bfm (\xv\mid  \yv): \nn\\
&  \hspace{2cm}\psib_{K,d} = \left(\omegab,\upsilonb_d,\Sigmab\right) \in \widetilde{\Omegab}_K\times \UpsilondK\times \bfV_K 
\Big\},
\end{align}
where $s_{\psib_{K,d}}$, $ \widetilde{\Omegab}_K$, $\UpsilondK$ and $\bfV_K$ are define previously in \eqref{eq_define_GLoME}, \eqref{define_BoundedGatesParameters}, \eqref{eq_define_meanExperts_polynomial} (or more general \eqref{eq_define_meanExperts_linearBounded}) and \eqref{eq_define_VarianceExperts}, respectively.

While for the BLoME model, we also need to select $\bfB$ from a list of candidate structures $\left(\cB_k\right)_{k \in [K]} \equiv \left(\cB\right)_{k \in [K]}$, where $\cB$ denotes the set of all possible partitions of the covariables indexed by $[D]$, for each cluster of individuals.The collection of BLoME models is defined as follows:  $ \left(S_\bfm\right)_{\bfm \in \cM}$, $\cM = \left\{\left(K,d,\bfB\right): K \in \cK,d \in \cD_\Upsilonb, \bfB \in \left(\cB\right)_{k \in [K]}\right\}$,
\begin{align} \label{eq_define_BLoME_Boundedness}
	S_\bfm 
	&=  \bigg\{(\xv,\yv) \mapsto
	s_{\psib_{K,d,\bfB}}(\xv\mid  \yv): =s_\bfm (\xv\mid  \yv): \nn\\
	&  \hspace{2cm} \psib_{K,d,\bfB} = \left(\omegab,\upsilonb_{d},\Sigmab\left(\bfB\right)\right) \in \widetilde{\Omegab}_K\times \UpsilondK \times \widetilde{	\bfV}_K\left(\bfB\right)
	\bigg\},
\end{align}
where $	s_{\psib_{K,d,\bfB}}$, $\widetilde{\Omegab}_K$, $\UpsilondK$, and $\widetilde{	\bfV}_K\left(\bfB\right)$ are defined in \eqref{eq_define_GLoME},  \eqref{define_BoundedGatesParameters}, \eqref{eq_define_meanExperts_polynomial} (or more generally \eqref{eq_define_meanExperts_linearBounded}), and \eqref{eq_define_blockdiagonal_covariances_boundedness}, respectively.
%
%
In theory, we would like to consider the whole collection of models $ \left(S_\bfm\right)_{\bfm \in \cM}$.
However, the cardinality of $\cB$ is large; its size is a Bell number. Even for a moderate number of variables $D$, it is not possible to explore the set $\cB$, exhaustively. We restrict our attention to a random subcollection $\cB^R$ of moderate size. For example, we can consider the BLLiM procedure from \citep[Section 2.2]{devijver2017nonlinear}.

\begin{remark}\label{remark_eq_define_GLoME_Boundedness}
It is worth noting that we can also define the collection of the forward GLoME and BLoME models in the same framework as in \eqref{eq_define_GLoME_Boundedness} and \eqref{eq_define_BLoME_Boundedness}, respectively. For example, in GloME models, the unknown forward conditional density $s_0^*$ is estimated via the following collection of forward models $ \left(S^*_\bfm\right)_{\bfm \in \cM}$, with $\cM = \cK \times \cD_\Upsilonb$, and  
\begin{align} \label{eq_define_GLoME_forward_Boundedness}
	S^*_\bfm &=  \Big\{(\xv,\yv) \mapsto
	s_{\psib^*_K}(\yv\mid  \xv)=:s^*_\bfm(\yv\mid  \xv):\nn\\
	& \hspace{2cm} \psib^*_K = \left(\omegab^*,\upsilonb^*_d,\Sigmab^*\right) \in \widetilde{\Omegab}_K^*\times \UpsilondKstar\times \bfV_K^* =:\widetilde{\Psib}^*_{K,d}\Big\},
\end{align}
where $ \widetilde{\Omegab}_K^*$, $\UpsilondKstar$ and $\bfV_K^*$ are defined similar to \eqref{define_BoundedGatesParameters}, \eqref{eq_define_meanExperts_polynomial} (or more general \eqref{eq_define_meanExperts_linearBounded}) and \eqref{eq_define_VarianceExperts}, respectively.

Motivated by the inverse conditional densities \eqref{eq_inverseGLLiM} and the one-to-one correspondence in \cref{lem_inverse_regression_mapping}, for the sake of simplicity of notation, we only state our main \cref{thm_Oracle_Inequality_GLoME,thm_Oracle_Inequality_BLoME} for the collection of inverse models $\left(S_\bfm\right)_{\bfm \in \cM}$ in \eqref{eq_define_GLoME_Boundedness} and \eqref{eq_define_BLoME_Boundedness}, respectively.
%
%
However, our finite-sample oracle inequalities, \cref{thm_Oracle_Inequality_GLoME,thm_Oracle_Inequality_BLoME} can be applied to the forward model $ \left(S^*_\bfm\right)_{\bfm \in \cM}$, established in \eqref{eq_define_GLoME_forward_Boundedness}, if we consider $\yv$ and $\xv$ as realizations of predictors and response variables, respectively. 

\end{remark}

\subsection{Loss functions and penalized maximum likelihood estimator}
In the maximum likelihood approach, the Kullback--Leibler divergence is the most natural loss function, which is defined for two densities $s$ and $t$ by
\begin{align*}
\kl(s,t) = \begin{cases}
	\int_{\R^D} \ln\left(\frac{s(y)}{t(y)}\right)s(y)dy& \text{ if $sdy$ is absolutely continuous \wrt~$tdy$},\\
	+\infty & \text{ otherwise}.
\end{cases}
\end{align*}
However, to take into account the structure of conditional
densities and the random covariates $\left(\Yv_{[n]}\right)$, we consider the \emph{tensorized Kullback--Leibler divergence} $\tkl$, defined as:
\begin{align} \label{tensorizedKLD}
\tkl(s,t) = \E{\Yv_{[n]}}{\frac{1}{n} \sum_{i=1}^n \kl\left(s\left(\cdot\mid  \Yv_i\right),t\left(\cdot\mid  \Yv_i\right)\right)},
\end{align}
if $sdy$ is absolutely continuous \wrt~$tdy$, and $+\infty$ otherwise. Note that if the predictors are fixed, this divergence is the classical fixed design type divergence in which there is no expectation. We refer to our result as a \emph{weak oracle inequality}, because its statement is based on a smaller divergence, when compared to $\tkl$, namely the \emph{tensorized Jensen-Kullback--Leibler divergence}:
\begin{align*}
\jtkl(s,t) = \E{\Yv_{[n]}}{\frac{1}{n} \sum_{i=1}^n \frac{1}{\rho}\kl\left(s\left(\cdot\mid  \Yv_i\right),\left(1-\rho\right)s\left(\cdot\mid  \Yv_i\right)+\rho t\left(\cdot\mid  \Yv_i\right)\right)},
\end{align*}
with $\rho \in \left(0,1\right)$. We note that $\jtkl$ was first used in \cite{cohen2011conditional,cohen2013partition}. However, a version of this divergence appears explicitly with $\rho = \frac{1}{2}$ in \cite{massart2007concentration}, and it is also found implicitly in \cite{birge1998minimum}. This loss is always bounded by $\frac{1}{\rho}\ln \frac{1}{1-\rho}$ but behaves like $\tkl$, when $t $ is close to $s$. The main tools in the proof of such a weak oracle inequality are deviation inequalities for sums of random variables and their suprema. These tools require a boundedness assumption on the controlled functions which is not satisfied by $-\ln\frac{s_\bfm }{s_0}$, and thus also not satisfied by $\tkl$. Therefore, we consider instead the use of $\jtkl$. In particular, in general, it holds that  $C_\rho \thel \le \jtkl \le \tkl$, where  $C_\rho = \frac{1}{\rho}\min\left(\frac{1-\rho}{\rho},1\right)\left(\ln\left(1+\frac{\rho}{1-\rho}\right)-\rho\right)$ (see \cite[Prop. 1]{cohen2011conditional}) and $\thel$ is a tensorized extension of the squared Hellinger distance $\thel$, defined by
\begin{align}
\thel(s,t) = \E{\Yv_{[n]}}{\frac{1}{n} \sum_{i=1}^n \hel\left(s\left(\cdot\mid  \Yv_i\right), t\left(\cdot\mid  \Yv_i\right)\right)}. \label{eq_tensorized_Hellinger}
\end{align}
Moreover, if we assume that, for any $\bfm \in \cM$ and any $s_\bfm  \in S_\bfm , s_0 d\lambda \ll s_\bfm  d\lambda$, then (see \cite{montuelle2014mixture,cohen2011conditional})
\begin{align} \label{montuelle.KLandJKL.1}
\frac{C_\rho}{2+\ln \norm{s_0/s_\bfm }}_\infty \tkl(s_0,s_\bfm ) \le \jtkl(s_0, s_\bfm ).
\end{align} 

In the context of MLE, given the collection of conditional densities $S_\bfm $,
we aim to estimate $s_0$ by the conditional density $\widehat{s}_\bfm$ that maximizes the likelihood (conditionally to $\left(\yv_i\right)_{i\in[n]}$) or equivalently that minimizes the negative log-likelihood (NLL), which we can write as:
\begin{align*}
	\widehat{s}_\bfm = \argmin_{s_\bfm  \in S_\bfm } \sum_{i=1}^n -\ln \left[s_\bfm  \left(\xv_i \mid  \yv_i \right)\right].
\end{align*}
We should work with almost minimizer of this quantity and define a $\eta$-log-likelihood minimizer as any $\widehat{s}_\bfm$ that satisfies:
\begin{align}\label{eq_define_NLL}
	\sum_{i=1}^n -\ln \left[\widehat{s}_\bfm  \left(\xv_i \mid  \yv_i \right)\right]  \le \inf_{s_\bfm  \in S_\bfm } \sum_{i=1}^n -\ln \left[s_\bfm  \left(\xv_i \mid  \yv_i \right)\right] + \eta,
\end{align}
where the error term $\eta$ is necessary to avoid any existence issue, \eg~the infimum may not be reached. Roughly speaking, the Ekeland variational principle asserts that, for any extended-valued lower semicontinuous function which is bounded below, one can add a small perturbation to ensure the existence of the minimum, see \eg~\cite[Chapter 2]{borwein_techniques_2004}. This framework is also used in \cite{cohen2011conditional,cohen2013partition,montuelle2014mixture}.


As always, using the NLL of the estimate in each model as a criterion is not sufficient. It is an underestimation of the risk of the estimate and this leads to choosing models that are too complex. In the context of PMLE, by adding a suitable penalty $\pen(\bfm)$, one hopes to create a trade-off between good data fit and model complexity.
For a given choice of $\pen(\bfm)$, the \emph{selected model} $S_{\widehat{\bfm}}$
is chosen as the one whose index is an $\eta'$-almost minimizer of the sum of the NLL and this penalty:
%
%
\begin{align}\label{eq_define_penalizedNLL}
	\sum_{i=1}^n -\ln \left[\widehat{s}_{\widehat{\bfm}} \left(\xv_i \mid  \yv_i\right)\right]+ \pen\left(\widehat{\bfm}\right) \le \inf_{\bfm \in \cM} \left\{\sum_{i=1}^n -\ln \left[\widehat{s}_\bfm \left(\xv_i \mid  \yv_i\right)\right]+ \pen(\bfm)\right\} + \eta'.
\end{align}
Note that $\widehat{s}_{\widehat{\bfm}}$ is then called the $\eta'$-penalized likelihood estimate and depends on both the error terms $\eta$ and $\eta'$. From hereon in, the term \emph{selected model (estimate) or best data-driven model (estimate)} is used to indicate that it satisfies the definition in \eqref{eq_define_penalizedNLL}.

\section{Main results of finite-sample oracle inequalities} \label{mainResult}
\subsection{Deterministic collection of GLoME models}
\cref{thm_Oracle_Inequality_GLoME} provides a lower bound on the penalty function, $\pen(\bfm)$, which guarantees that the PMLE selects a best data-driven model that performs almost as well as the best model. 
We highlight our main contribution and briefly establish the proof of \cref{thm_Oracle_Inequality_GLoME} in \cref{sec_contribution_proof_GLoME}.  A more detail technical proof can be found in \cref{sec_proof_lem_Bracketing_Entropy_Gates_intro_mine,sec_proof_lem_Bracketing_Entropy_Experts}.
\begin{theorem}[Finite-sample oracle inequality for GLoME models] \label{thm_Oracle_Inequality_GLoME}
Assume that we observe $\left(\xv_{[n]},\yv_{[n]}\right)$, arising from an unknown conditional density $s_0$. Given a collection of GLoME models, $ \left(S_\bfm \right)_{\bfm \in \cM}$, there is a constant $C$ such that for any $\rho \in (0,1)$, for any $\bfm\in \cM$, $ z_\bfm \in \R^+$, $\Xi = \sum_{\bfm \in \cM}e^{-z_\bfm} < \infty$ and any $C_1 > 1$, there is a constant $\kappa$ depending only on $\rho$ and $C_1$, such that if for every index $\bfm\in \cM$,
	\begin{align}
		\pen(\bfm) > \kappa\left[\left(C+\ln n\right)\dim \left(S_\bfm \right)+ z_\bfm\right], \label{eq_lowerBoundPenalty}
\end{align}
\noindent
then the $\eta'$-penalized likelihood estimate $\widehat{s}_{\widehat{\bfm}}$, defined in \eqref{eq_define_NLL} and \eqref{eq_define_penalizedNLL},
satisfies
\begin{align} \label{montuelle.OracleIneq.Gaussian}
	\E{\Xv_{[n]},\Yv_{[n]}}{\jtkl\left(s_0,\widehat{s}_{\widehat{\bfm}}\right)} &\le C_1 \inf_{\bfm \in \cM} \left(\inf_{s_\bfm  \in S_\bfm } \tkl\left(s_0,s_\bfm \right)+ \frac{\pen(\bfm)}{n}\right)  \nn\\
	& \quad + \frac{\kappa C_1 \Xi }{n}+ \frac{ \eta +\eta'}{n}.
\end{align}
\end{theorem}
It is worth noting that \cref{thm_Oracle_Inequality_GLoME} extends a corollary of \citep[Theorem 1]{montuelle2014mixture}, which can be verified via Lemma 1 from \cite{nguyen2020approximationMoE}, which makes explicit the relationship between softmax and Gaussian gating classes.

	Note that we only aim to construct upper bounds and do not focus on the important question of the existence of the corresponding lower bounds or minimax adaptive estimation. 
	However, as a special case of GLoME models, \cite{maugis2011non,maugis2013adaptive} consider a collection of univariate Gaussian mixture models having the same and known component variance. They show that the PMLE $\widehat{s}_{\widehat{\bfm}}$ is minimax adaptive to the regularity $\beta$, $\beta > 0$, of univariate density classes $\cH_{\beta}$ whose logarithm of the elements  is locally $\beta$-H\"older, with convergence rate $n^{-\frac{2\beta}{2\beta +1}}$ up to a power of $\ln(n)$. Although this result is stated for the Hellinger risk, it remains valid for the Kullback–Leibler divergence if we further assume that $\ln\left(\|s/t\|_\infty\right)$ is uniformly bounded on $\cup_{\bfm \in \cM} S_{\bfm}$, see \eg~Lemma 7.23 in \cite{massart2007concentration}. Note that this assumption guarantees that the Kullback–Leibler divergence and the Hellinger distance are equivalent.

A special case of GLoME models, namely model-selection (nearly-D-sparse) aggregation in mixture models, is considered by \cite{bunea_spades_2010,bertin_adaptive_2011,rigollet_kullbackleibler_2012,butucea_optimal_2017,dalalyan_optimal_2018} and is related to our model selection results. 
More precisely, the authors of \cite{rigollet_kullbackleibler_2012} did not consider PMLEs with Kullback-Leibler loss but only with $L_2$-loss or aggregation of a finite number of densities. Similarly, the results from \cite{bunea_spades_2010,bertin_adaptive_2011} dealt with the $L_2$-loss and investigate the Lasso and the Dantzig estimators, respectively, suitably adapted to the problem of density estimation. With respect to the Kullback–Leibler loss, \cite{butucea_optimal_2017} established model selection type oracle inequalities with high probability rather than in expectation. In particular, a bound in expectation with Kullback–Leibler loss in \cite{dalalyan_optimal_2018} is perhaps the most relevant reference to our result.
They established exact oracle inequalities with a rate-optimal remainder term $\left(\left(\ln K\right)/n\right)^{1/2}$, up to a possible logarithm correction, in the problem of convex aggregation when the number $K$ of components is larger than $n^{1/2}$. It is worth noting that they did not consider PMLEs as we do or as \cite{maugis2011non,maugis2013adaptive}, see \eg~\eqref{eq_define_penalizedNLL}. Instead, the constraint that the weight vector in the mixture model belongs to the probabilistic simplex acts as a sparsity-inducing penalty. However, in their collection of mixture models, the mixture components are deterministic and chosen from a dictionary obtained on the basis of previous experiments or expert knowledge rather estimated from data. The adjective exact refers to the fact that the ``bias term" $\tkl\left(s_0,s_\bfm \right)$ is not multiplied by factor strictly larger than one as  in our \cref{thm_Oracle_Inequality_GLoME}.

	To the best of our knowledge, providing a minimax analysis for the PMLEs or the problem of model-selection (nearly-D-sparse) aggregation in the context of standard MoE models is still an open question. In particular, it is not trivial to extend such optimal risk bounds regarding Gaussian mixtures, see \eg~\cite[Theorem 2.8]{maugis2013adaptive} or \cite[Theorems 2.3 and 3.1]{dalalyan_optimal_2018}, to standard MoE models. However, we wish to provide such a minimax analysis in future research.

%
%
%

In contract to \cref{thm_Oracle_Inequality_GLoME}, in \cref{thm_Oracle_Inequality_BLoME}, we consider a random subcollection of models $\widetilde{\cM}$, included in the whole collection $\cM$. This is particularly useful in a high-dimensional context where we cannot test all the models. Therefore, we have to restrict ourselves to a smaller subcollection of models, which is then potentially random.

%
\subsection{Random subcollection of BLoME models}
Note that the constructed collection of models with block-diagonal structures for each cluster of individuals is designed, for example, by the BLLiM procedure from \cite{devijver2017nonlinear}, where each collection of partitions is sorted by sparsity level. Nevertheless, our finite-sample oracle inequality, \cref{thm_Oracle_Inequality_BLoME}, still holds for any random subcollection of $\cM$, which is constructed by some suitable tools in the framework of BLoME regression models.
%
We highlight our main contribution and briefly establish the proof of \cref{thm_Oracle_Inequality_BLoME} in \cref{sec_contribution_collection_random_MoE}.  A more detailed technical proof can be found in \cref{sec_lem_Inequality_Hellinger_Supx_dy_dk_intro,proofLemma_nguyen2021nonBLoME_intro,sec_proof_lem_bracketingEntropyGaussianBlock_intro}.
\begin{theorem}[Finite-sample oracle inequality for BLoME models]\label{thm_Oracle_Inequality_BLoME}
Let $(\xv_{[n]},\yv_{[n]})$ be observations coming from an unknown conditional density $s_0$. For each $\bfm = \left(K,d,\bfB\right) \in \left(\cK \times \cD_\Upsilonb \times \cB\right)\equiv \cM$, let $S_\bfm$ be define by \eqref{eq_define_BLoME_Boundedness}. Assume that there exists $\tau > 0$ and $\epsilon_{KL} > 0$ such that, for all $\bfm \in \cM$, one can find $\bar{s}_\bfm \in S_\bfm$, such that
\begin{align}\label{eq_inequ_random_subcollection}
	\tkl\left(s_0,\bar{s}_\bfm\right) \le \inf_{t \in  S_\bfm} \tkl\left(s_0,t\right) + \frac{\epsilon_{KL}}{n}, \text{ and } \bar{s}_\bfm\ge e^{-\tau} s_0.
\end{align}
Next, we construct some random subcollections $\left(S_\bfm\right)_{\bfm \in \widetilde{\cM}}$ of $\left(S_\bfm\right)_{\bfm \in \cM}$ by letting $\widetilde{\cM}\equiv\left(\cK \times \cD_\Upsilonb \times \cB^R\right) \subset \cM$, where $\cB^R$ is a random subcollection $\cB$, of moderate size.
Consider the collection $\left(\widehat{s}_\bfm\right)_{\bfm \in \widetilde{\cM}}$ of $\eta$-log likelihood minimizers satisfying \eqref{eq_define_NLL} for all $\bfm \in \widetilde{\cM}$.
Then, there is a constant $C$ such that for any $\rho \in (0,1)$, and any $C_1 > 1$, there are two constants $\kappa$ and $C_2$ depending only on $\rho$ and $C_1$ such that, for every index, $\bfm \in \cM$, $ z_\bfm \in \R^+$, $ \Xi = \sum_{\bfm \in \cM}e^{-z_\bfm} < \infty$ and 
\begin{align}
	\pen(\bfm) \ge \kappa \left[\left(C+\ln n\right)\dim(S_\bfm) + (1 \vee \tau)z_\bfm\right], \label{eq_Oracle_Inequality_BLoME_pentalty}
\end{align}
the $\eta'$-penalized likelihood estimator $\widehat{s}_{\widehat{\bfm}}$, defined as in \eqref{eq_define_penalizedNLL} on the subset $\widetilde{\cM} \subset \cM$,
satisfies
\begin{align}
	\E{\Xv_{[n]},\Yv_{[n]}}{\jtkl\left(s_0,\widehat{s}_{\widehat{\bfm}}\right)}&\le C_1\E{\Xv_{[n]},\Yv_{[n]}}{\inf_{\bfm \in \widetilde{\cM}}\left(\inf_{t \in S_\bfm} \tkl\left(s_0,t\right)+2 \frac{\pen(\bfm)}{n}\right)} \nn\\
	&\quad 
	+ C_2(1 \vee \tau)\frac{\Xi^2}{n} + \frac{\eta' + \eta}{n}. \label{eq_Oracle_Inequality_BLoME}
\end{align}
\end{theorem}
It is worth mentioning that the comparison of the two inequalities in \cref{thm_Oracle_Inequality_GLoME,thm_Oracle_Inequality_BLoME} involves the following major differences. First, to control the random subcollection, in \cref{thm_Oracle_Inequality_BLoME},  we further assume condition \eqref{eq_inequ_random_subcollection}.
However, this is not a strong assumption and is satisfied, for example, if $s_0$ is bounded, with a compact support.
Assumption \eqref{eq_inequ_random_subcollection} is needed because we consider a random subcollection from the whole collection. Thanks to this assumption, we can utilize Bernstein's inequality to control the additional randomness. It is important to stress that the parameter $\tau$ depends on the true unknown density $s_0$ and cannot be explicitly determined for this reason. It appears not only in the oracle type inequality \eqref{eq_Oracle_Inequality_BLoME}, but also in the penalty term \eqref{eq_Oracle_Inequality_BLoME_pentalty}. However, in some special cases of BLoME models, for instance finite mixture regression models, we can construct a larger penalty independent of $\tau$ to obtain an oracle type inequality but the price to pay for getting rid of $\tau$ in the risk bound is the increased value of error upper bound, see Appendix C in \cite{devijver2017joint} for more details.
Second, the constant $\Xi$ associated to the Kraft-type inequality for the collection appears squared in the upper bound of the oracle inequality in \cref{thm_Oracle_Inequality_BLoME}. This is because of the random subcollection $\widetilde{\cM}$ of $\cM$, if the model collection is fixed, we get a linear bound as in \cref{thm_Oracle_Inequality_GLoME}. 

\subsection{Practical application of finite-sample oracle inequalities} \label{sec_practical_application}
It is important to mention that \cref{thm_Oracle_Inequality_GLoME,thm_Oracle_Inequality_BLoME} provide some theoretical justification regarding the penalty shapes when using the slope heuristic for our collection of MoE models, including GLLiM, GLoME, BLLiM, BLoME models. 

More precisely, our oracle inequalities show that the performance in $\jtkl$ loss of our PMLEs are roughly comparable to that of oracle models if we take large enough constants in front of the penalties, whose forms are only known up to multiplicative constants and proportional to the dimensions of models. This motivate us to make use of the slope heuristic criterion to select our data-driven hyperparameters, including the number of mixture components, the degree of polynomial Gaussian expert's mean functions, and the potential hidden block-diagonal structures of the covariance matrices of the multivariate predictor, see more in \cref{sec_numerical_Experiment}.  

To be more precise, we consider the condition \eqref{penaltyLowerBound} for GLoME models. As shown in the proof of Theorem 3.1, 
in fact we can replace the assumption on $\pen(\bfm)$ by a milder one. More precisely, given a constant $\fC$, which we specify later, there is a constant $\kappa$ depending only on $\rho$ and $C_1$, such that
for every index $\bfm \in \cM$, $\pen(\bfm)$ is bounded from below by
\begin{align} \label{penaltyLowerBound}
	\kappa \left(\dim(S_\bfm ) \left(2 \left(\sqrt{\fC} + \sqrt{\pi}\right)^2 + \left(\ln \frac{n}{\left(\sqrt{\fC} + \sqrt{\pi}\right)^2\dim\left(S_\bfm \right)}\right)_+\right)+ z_\bfm\right).
\end{align}
In our numerical experiment, namely, \cref{fig_Slope_Heuristic_2000_DDSE_Slope_CNLL,fig_Slope_Heuristic_10000_DDSE_Slope_CNLL}, we confirm the validity of the linear penalty shape assumption, which supports the use of penalties proportional to the dimension. This implies that the logarithm terms are not detected in practice as shown in \cref{fig_Slope_Heuristic_2000_DDSE_Slope_CNLL,fig_Slope_Heuristic_10000_DDSE_Slope_CNLL} and thus only the preponderant term in $\dim\left(S_\bfm\right)$ is retained in the penalty form.

In particular, based on Appendix B.4 from \cite{cohen2011conditional}, we can make explicit the dependence of the theoretical constant $\kappa$, with respect to $\rho$ and $C_1$ as follows. For instance, in \cref{thm_Oracle_Inequality_GLoME}, given any $\rho \in (0,1)$ and $C_1 > 1$, define $\epsilon_{\pen} = 1-\frac{1}{C_1}$. Then $\kappa$ is determined by 
\begin{align*}
	\kappa = \frac{\kappa'_0 \left(\kappa'_1 + \kappa'_2\right)^2 \left(\sqrt{1+ \frac{72 C_{\rho} \epsilon_{\pen}}{\rho \kappa'_0 \left(\kappa'_1 + \kappa'_2\right)^2}}+1\right)}{2 C_{\rho} \epsilon_{\pen}} + \frac{18}{\rho},
\end{align*}
where $\epsilon_d$ is a given positive constant and
\begin{align*}
	\kappa'_0 &= \frac{2\left(2+ \epsilon_d\right)}{1+ \epsilon_d},
	\kappa'_1 =  \frac{1}{\sqrt{\rho \left(1- \rho\right)}} \left(3 \kappa'_3 \sqrt{2}+12+16\sqrt{\frac{1-\rho}{\rho}}\right), \\
	\kappa'_3 &\le 27,
	\kappa'_2 = \frac{1}{\sqrt{\rho \left(1-\rho\right)}} \left(42+ \frac{3}{4\sqrt{\kappa'_0}}\right) .
\end{align*}

For example, if $\rho = \frac{1}{2}$, $C_1 = 2$, $\epsilon_d = 1$, $\kappa_3' = 27$, then $\epsilon_{\pen} = 1-\frac{1}{C_1} = \frac{1}{2}$, $\kappa_0'=3$, $\kappa_1' = 56+162\sqrt{2}$, $\kappa_2' =  84+\frac{\sqrt{3}}{2}$, and
\begin{align} \label{eq_kappa_practice}
	C_\rho &= \frac{1}{\rho}\min\left(\frac{1-\rho}{\rho},1\right)\left(\ln\left(1+\frac{\rho}{1-\rho}\right)-\rho\right) = 2 \ln2 -1 , \nn\\
	%
	%
	\kappa  & \approx 2126069.
	%
	%
\end{align}
According to the previous example, we can see that the theoretical penalty is lower bounded by $\kappa$, which can be large in practice. This result is not surprising since, according to \citep[Appendix B.4, page 40, line 7]{cohen2011conditional}, if we choose $\epsilon_d$ small enough then $\kappa$ scales proportionally to 
\begin{align*}
	\frac{1}{C_\rho \rho (1-\rho)\epsilon_{\pen}} 
	= \frac{\rho}{(1-\rho)^2 \left(\ln\left(1+ \frac{\rho}{1-\rho}\right)-\rho\right)}
\end{align*}
and thus explodes to $+\infty$ when $\rho$ goes to $1$ and $C_1$ goes to $1$.	
Therefore, it is important to study a natural question as to whether the constant $\kappa$ appearing in the penalty can be estimated from the data without losing theoretical guaranties on the performance? No definite answer exists so far for MoE regression models, however at least our numerical experiment in \cref{sec_numerical_Experiment} shows that the slope heuristic proposed by \cite{birge2007minimal} may lead to a good practical solution. In particular, we seek to mathematically and fully justify the slope heuristic in MoE regression models as in least-squares regression on a random (or fixed) design with regressogram (projection) estimators, respectively, see \cite{birge2007minimal,arlot_data_driven_2009,arlot_data_driven_2009_NIPS,arlot2019minimal} for greater details. Furthermore, as is often the case in empirical processes theory, the constant $\kappa$ appearing in the bound is pessimistic as shown in \eqref{eq_kappa_practice}. Numerical experiments in \cref{sec_numerical_Experiment} show that there is a hope that this is only a technical issue.
 For instance, \cref{fig_Slope_Heuristic_2000_Djump_Slope_CNLL,fig_Slope_Heuristic_10000_Djump_Slope_CNLL} show that the practical values of $\kappa$, selecting by slope heuristic, are generally very small  compared to the theoretical values.


%

\section{Proofs of  finite-sample oracle inequalities and main mathematical challenges} \label{proofMainTheorem}
\subsection{Proof for deterministic collections of GLoME models}
The deterministic collection of MoE models include GLLiM, GLoME, SGaME and LinBoSGaME models, where finite-sample oracle inequalities have only been well studied for the latter two models \citep[Theorem 2]{cohen2011conditional}, see also \citep[Theorem 2.2.]{cohen2013partition} and \citep[Theorem 2]{montuelle2014mixture}. We first summarize this general model selection theorem and the techniques that \cite{montuelle2014mixture} used to control the bracketing entropy of LinBoSGaME models with softmax gating networks in \cref{sec_general_model_selection_intro,sec_contribution_proof_LinBoSGaME}, respectively. Then, we explain why such techniques can not be directly applied to our collection of GLoME models via \cref{sec_contribution_proof_GLoME} and propose our approach to highlight the main issues and our contributions.

\subsubsection{A general conditional density model selection theorem for deterministic collection of models }\label{sec_general_model_selection_intro}
Before stating a general model selection theorem for conditional densities, we have to present some regularity assumptions.

First, we need an information theory type assumption to control the complexity of our collection. We assume the existence of a Kraft-type inequality for the collection \citep{massart2007concentration,barron2008mdl}.
\begin{assumption}[K] \label{assumption_K_intro}
There is a family $\left(z_\bfm\right)_{\bfm \in \cM}$ of non-negative numbers and a real number $\Xi$ such that
\begin{align*}
	\Xi = \sum_{\bfm \in \cM}e^{-z_\bfm} < + \infty.
\end{align*}
\end{assumption}

For technical reasons, a separability assumption, always satisfied in the setting of this paper, is also required. It is a mild condition, which is classical in empirical process theory \citep{van1996weak,geer2000empirical}. \cref{assumption_Sep_intro} allows us to work with a countable subset. 
\begin{assumption}[Sep]\label{assumption_Sep_intro}
For every model $S_\bfm $ in the collection $\left(S_\bfm \right)_{\bfm\in\cM}$, there exists some countable subset $S'_\bfm$ of $S_\bfm $ and a set $\cX'_\bfm$ with $\iota \left(\cX \setminus \cX'_\bfm\right) = 0$, where $\iota$ denotes Lebesgue measure, such that for every $t \in S_\bfm $, there exist some sequences $\left(t_k\right)_{k \geq 1}$ of elements of $S'_\bfm$, such that for every $\yv \in \cY$ and every $\xv \in \cX'_\bfm, \ln\left(t_k\left(\xv \mid  \yv\right)\right) \xrightarrow{k \rightarrow + \infty}  \ln\left(t\left(\xv \mid  \yv\right)\right)$.
\end{assumption} 

Next, recall that the bracketing entropy of a set $S$ with respect to any distance $d$, denoted by $\cH_{\left[\cdot\right],d}(\left(\delta,S\right))$, is defined as the logarithm of the minimal number $\cN_{\left[\cdot\right],d}\left(\delta,S\right)$ of brackets $\left[t^{-},t^{+}\right]$ covering $S$, such that $d(t^-,t^+) \le \delta$. That is,
\begin{align}\label{eq_define_BracketingEntropy_intro}
N_{\left[\cdot\right],d}\left(\delta,S\right):= \min\left\{n \in \Ns : \exists\left[ t^-_k,t^+_k\right]_{k\in[n]}\text{ s.t }d(t^-_k,t^+_k) \le \delta,S \subset \bigcup_{k=1}^n \left[t^-_k,t^+_k\right] \right\}.
\end{align}
Here, $s \in \left[t^-_k,t^+_k\right]$ means that $t^{-}_k(\xv,\yv) \le s(\xv,\yv)\le t^{+}_k(\xv,\yv)$, $\forall (\xv,\yv) \in \cX \times \cY$.

We also need the following important \cref{assumption_H_intro} on Dudley-type integral of these bracketing entropies, which is often utilized in empirical process theory \citep{van1996weak,geer2000empirical,kosorok2007introduction}.
\begin{assumption}[H]\label{assumption_H_intro}
For every model $S_\bfm $ in the collection $\cS$, there is a non-decreasing function $\phi_\bfm$ such that $\delta \mapsto \frac{1}{\delta}\phi_\bfm(\delta)$ is non-increasing on $\left(0,\infty\right)$ and for every $\delta \in \R^+$,
\begin{align*}
	\int_{0}^\delta \sqrt{\cH_{\left[\cdot\right],\thell}\left(\delta,S_\bfm \left(\widetilde{s},\delta\right)\right)}d\delta \le \phi_\bfm(\delta),
\end{align*}
where $S_\bfm \left(\widetilde{s},\delta\right) = \left\{s_\bfm  \in S_\bfm  : \thell \left(\widetilde{s},s_\bfm \right) \le \delta\right\}$.
The model complexity  of $S_\bfm $ is then defined as $\cD_\bfm$ =$n \delta^2_\bfm$, where $\delta_\bfm$ is the unique root of $\frac{1}{\delta} \phi_\bfm(\delta) = \sqrt{n} \delta$.

\end{assumption} 
Observe that the model complexity does not depend on the bracketing entropies of the global models $S_\bfm $, but rather on those of smaller localized sets $S_\bfm \left(\widetilde{s},\delta\right)$. We are now able to state an important weak oracle inequality, \cref{thm_weakOracleInequality_PennecEL_intro}, originally from \cite[Theorem 2]{cohen2011conditional,cohen2013partition}.
\begin{theorem} \label{thm_weakOracleInequality_PennecEL_intro}
Assume that we observe $\left(\xv_{[n]},\yv_{[n]}\right)$, arising from an unknown conditional density $s_0$. Let $\cS = \left(S_\bfm \right)_{\bfm \in \cM}$ be an at most countable conditional density model collection. Assume that \cref{assumption_K_intro} (K), \cref{assumption_Sep_intro} (Sep), and \cref{assumption_H_intro} (H) hold for every model $S_\bfm  \in \cS$.
Then, for any $\rho \in (0,1)$ and any $C_1 > 1$,  there is a constant $\kappa$ depending only on $\rho$ and $C_1$, such that for every index $\bfm\in \cM$,
\begin{align*}
	\pen(\bfm) \geq \kappa\left(n \delta_\bfm^2+ z_\bfm\right)
\end{align*}
with $\delta_\bfm$ is the unique root of $\frac{1}{\delta} \phi_\bfm(\delta) = \sqrt{n} \delta$, such that the $\eta'$-penalized likelihood estimator $\widehat{s}_{\widehat{\bfm}}$ 
satisfies
\begin{align} \label{montuelle.OracleIneqGenConDensity_intro}
	\E{\Xv_{[n]},\Yv_{[n]}}{\jtkl\left(s_0,\widehat{s}_{\widehat{\bfm}}\right)} &\le C_1 \inf_{\bfm \in \cM} \left(\inf_{s_\bfm  \in S_\bfm } \tkl\left(s_0,s_\bfm \right)+ \frac{\pen(\bfm)}{n}\right) \nn\\
	& \quad  + \frac{\kappa C_1 \Xi }{n}+ \frac{ \eta +\eta'}{n}.
\end{align}
\end{theorem}

For the sake of generality, \cref{thm_weakOracleInequality_PennecEL_intro} is relatively abstract. Since the assumptions of \cref{thm_weakOracleInequality_PennecEL_intro} are as general as possible. But from the practical point of view, a natural question is the existence of interesting model collections that satisfy these assumptions. We will sketch the proof for collection of LinBoSGaME models from \citep[Theorem 1]{montuelle2014mixture} and show that their result is not directly applicable to the GLoME setting. The main reason is that the technique for handling the linear combination of bounded functions for the weight functions of logistic schemes of \cite{montuelle2014mixture} is not valid for the Gaussian gating parameters in GLoME models. 
Therefore, we propose a \emph{reparameterization trick}
%
%
to bound the metric entropy of the Gaussian gating parameters space; see \cref{sec_contribution_proof_GLoME} for more details. 


\subsubsection{Sketch of the proof for LinBoSGaME models}\label{sec_contribution_proof_LinBoSGaME}
To prove the main conditional density model selection theorem for LinBoSGaME models, \citep[Theorem 1]{montuelle2014mixture}, the authors have to make use of \cref{thm_weakOracleInequality_PennecEL_intro}. Then, they need to prove that their collection of LinBoSGaME models have to satisfy \cref{assumption_K_intro} (K), \cref{assumption_Sep_intro} (Sep), and \cref{assumption_H_intro} (H). However, they did not verify \cref{assumption_K_intro} (K) and \cref{assumption_Sep_intro} (Sep) and considered them as primative assumptions on their LinBoSGaME models because of the complexity of LinBoSGaME models and technical reasons. Therefore, the main difficulty remains on verifying \cref{assumption_H_intro} (H) via bracketing entropy controls of the linear combination of bounded functions for the weight functions of logistic schemes. 

Firstly, they define the following distance over conditional densities:
\begin{align*}
\sup_\yv d_\xv(s,t) = \sup_{\yv \in \cY} \left(\int_\cX \left(\sqrt{s(\xv\mid  \yv)} - \sqrt{t(\xv\mid  \yv)}\right)^2 d\xv \right)^{1/2}.
\end{align*}
This leads straightforwardly to $\thel(s,t) \le \sup_\yv d^2_\xv(s,t)$. Then, they also define
\begin{align*}
\sup_\yv d_k\left(g,g'\right) = \sup_{\yv\in \cY} \left(\sum_{k=1}^K \left(\sqrt{g_k(\yv)}-\sqrt{g'_k(\yv)}\right)^2\right)^{1/2},
\end{align*}
for any gating functions $g$ and $g'$.
To this end, given any densities $s$ and $t$ over $\cX$, the following distance, depending on $\yv$, is constructed as follows:
\begin{align*}
\sup_\yv \max_k d_\xv(s,t) &= \sup_{\yv \in \cY} \max_{k \in [K]} d_\xv\left(s_k(\cdot,\yv),t_k(\cdot,\yv)\right)
\nn\\
&
=\sup_{\yv \in \cY} \max_{k \in [K]} \left(\int_\cX \left(\sqrt{s_k(\xv,\yv)}-\sqrt{t_k(\xv,\yv)}\right)^2d\xv\right)^{1/2}.
\end{align*}

Then, they prove that definition of complexity of model $S_\bfm $ in \cref{assumption_H_intro} (H) is related to an classical entropy dimension with respect to a Hellinger type divergence $\thell$, due to \cref{prop_montuelle_Proposition1_intro}.
\begin{prop}[Proposition 2 from \citep{cohen2011conditional}]\label{prop_montuelle_Proposition1_intro}
For any $\delta \in (0,\sqrt{2}]$, such that $ \cH_{\left[\cdot\right],\thell}\left(\delta,S_\bfm \right) \le \dim(S_\bfm )\left(C_\bfm+\ln\left(\frac{1}{\delta}\right)\right)$, the function 
\begin{align*}
	\phi_\bfm\left(\delta\right) = \delta \sqrt{\dim\left(S_\bfm \right)}\left(\sqrt{C_\bfm}+\sqrt{\pi}+\sqrt{\ln\left(\frac{1}{\min\left(\delta,1\right)}\right)}\right)
\end{align*} satisfies \cref{assumption_H_intro} (H). Furthermore, the unique solution $\delta_\bfm$ of $\frac{1}{\delta} \phi_\bfm\left(\delta\right) =  \sqrt{n} \delta$, satisfies
\begin{align*}
	n \delta^2_\bfm \le \dim(S_\bfm ) \left(2 \left(\sqrt{C_\bfm} + \sqrt{\pi}\right)^2 + \left(\ln \frac{n}{\left(\sqrt{C_\bfm} + \sqrt{\pi}\right)^2\dim\left(S_\bfm \right)}\right)_+\right).
\end{align*}
\end{prop}
Therefore, \cref{prop_montuelle_Proposition1_intro} implies that \cref{assumption_H_intro} (H) can be proved via \cref{lem_model_complexity_Hellinger}.
\begin{lem}\label{lem_model_complexity_Hellinger}
For any $\delta \in (0,\sqrt{2}]$, the collection of LinBoSGaME models, $\cS = \left(S_\bfm \right)_{\bfm \in \cM}$, satisfies
\begin{align*}
	\cH_{\left[\cdot\right],\thell}\left(\delta,S_\bfm \right) \le \dim(S_\bfm )\left(C_\bfm+\ln\left(\frac{1}{\delta}\right)\right).
\end{align*}
\end{lem}
\cref{lem_model_complexity_Hellinger} is then obtained by decomposing the entropy terms between the softmax gating functions and the Gaussian experts.
Note that both LinBoSGaME and GLoME models share the same structures of Gaussian experts mean, see \cref{fig_MoE_GaussianGN} again
for more details.
Therefore, in \cref{sec_contribution_proof_GLoME}, we only highlight our contributions regarding the control of bracketing entropy for the parameter of gating network compared to \cite{montuelle2014mixture}.

More precisely, the authors of \cite{montuelle2014mixture} rewrite the softmax gating parameters' space
as follows: 
\begin{align}
%
\bfW_{K,d_{\bfW}} &= \left\{0\right\}\otimes \bfW^{K-1},
\bfW =  \left\{ \yv \mapsto
\sum_{d = 1}^{d_{\bfW}} \omega_d \theta_{\bfW,d}\left(\yv\right) \in \R:\max_{d \in \left[d_{\bfW}\right]}\left|\omega_d\right| \le T_{\bfW}\right\},\nn\\
\cP_K&= \left\{ \yv \mapsto\left(\frac{e^{\bfw_{k}(\yv)}}{\sum_{l=1}^K e^{\bfw_{l}(\yv)}}\right)_{k\in[K]} =: \left(g_{\bfw,k}\left(\yv\right)\right)_{k \in [K]}, \bfw \in \bfW_{K,d_{\bfW}}\right\}.\nn
\end{align}
Then, they also require the definition of metric entropy of the set $\cW_K$: $\entropy(\delta,\cW_K)$, which measures the logarithm of the minimal number of balls of radius at most $\delta$, according to a distance $d_{\norm{\sup}_\infty}$, needed to cover $\cW_K$ where
\begin{align} \label{eq.def.metricEntropy-Distance_intro}
d_{\norm{\sup}_\infty} \left(\left(s_k\right)_{k\in[K]},\left(t_k\right)_{k\in[K]}\right) = \max_{k \in [K]}\sup_{\yv \in \cY} \norm{s_k(\yv) - t_k(\yv)}_2,
\end{align}
for any $K$-tuples of functions $\left(s_k\right)_{k\in[K]}$, $\left(t_k\right)_{k\in[K]}$ and $\norm{s_k(\yv) - t_k(\yv)}_2$ is the Euclidean distance in $\R^L$.
By using Lemma 5 and Proposition 2 from \citep{montuelle2014mixture},  \cref{lem_model_complexity_Hellinger} holds true if we can prove  \cref{lem_Bracketing_Entropy_Gates_intro_mine_intro}. Note that the first inequality of \cref{lem_Bracketing_Entropy_Gates_intro_mine_intro} comes from \citep[Lemma 4]{montuelle2014mixture} and describes relationship between the bracketing entropy of $\cP_K$ and the entropy of $\cW_K$.
\begin{lem}\label{lem_Bracketing_Entropy_Gates_intro_mine_intro}
For all $\delta \in (0,\sqrt{2}]$, there exists a constant $C_{\cW_K}$ such that
\begin{align*}
	\cH_{[\cdot],\sup_\yv d_k}\left(\frac{\delta}{5},\cP_K\right)
	&\le 
	\entropy\left(\frac{ 3 \sqrt{3}\delta}{20 \sqrt{K-1}},\cW_K\right)\nn\\
	&\le \dim\left(\cW_K\right) \left(C_{\cW_K} + \ln \left(\frac{ 20 \sqrt{K-1}}{3 \sqrt{3}\delta}\right)\right).
\end{align*}
\end{lem}
By the linearity from the construction of linear combination of a finite set of bounded functions whose coefficients belong to a compact set, in the argument from \citep[Proof of Part 1 of Lemma 1, Page 1689]{montuelle2014mixture}, the second inequality of \cref{lem_Bracketing_Entropy_Gates_intro_mine_intro} is then easily established as follows. 
\begin{proof}[Proof of the second inequality of \cref{lem_Bracketing_Entropy_Gates_intro_mine_intro}]
Note that for all 
$$\bfw = \left(0,w_k\right)_{k\in[K-1]} \in \bfW_{K,d_{\bfW}}, \quad \bfv= \left(0,v_k\right)_{k\in[K-1]} \in \bfW_{K,d_{\bfW}},$$ it holds that 
\begin{align}
	d_{\left\|\sup\right\|_\infty} \left(\bfw - \bfv\right) 
	&= \max_{k \in [K-1]} \left\|\cW_K-\bfv_k\right\|_\infty\nn\\
	&=
	\max_{k \in [K-1]} \sup_{\yv \in \Yv} \left|\sum_{i = 1}^{d_{\bfW}} \omega_{k,i}^{\bfw} \theta_{\bfW,i}(\yv)-\sum_{i = 1}^{d_{\bfW}} \omega_{k,i}^{\bfv} \theta_{\bfW,i}(\yv)\right|\nn \\
	&\le \max_{k \in [K-1]} \sum_{i = 1}^{d_{\bfW}} \left|\omega_{k,i}^{\bfw}- \omega_{k,i}^{\bfv}\right|  \underbrace{\sup_{\yv \in \Yv} \left|\theta_{\bfW,i}(\yv)\right|}_{\le 1}\nn\\ 
	&\le d_{\bfW} \max_{k \in [K-1], i \in \left[ d_{\bfW}\right]} \left|\omega_{k,i}^{\bfw}- \omega_{k,i}^{\bfv}\right|.\nn
\end{align}
Therefore, we obtain
\begin{align}
	&\entropy\left(\frac{ 3 \sqrt{3}\delta}{20 \sqrt{K-1}},\cW_K\right)\nn\\
	&\le \entropy\left(\frac{ 3 \sqrt{3}\delta}{20 \sqrt{K-1}d_{\bfW}},\left\{ \omegab \in \R^{(K-1)d_{\bfW}}:\left\|\omegab\right\|_\infty \le T_{\bfW}\right\}\right)\nn\\
	&\le 
	(K-1) d_{\bfW}\ln\left(1+\frac{20 \sqrt{K-1}d_{\bfW} T_{\bfW}}{ 3 \sqrt{3}\delta}\right)\nn\\
	&\le (K-1) d_{\bfW}\ln\left(\sqrt{2}+\frac{20 \sqrt{K-1}d_{\bfW} T_{\bfW}}{ 3 \sqrt{3}}+\ln\left(\frac{1}{\delta}\right)\right)\nn.
\end{align}

\end{proof}
However, proving \cref{lem_Bracketing_Entropy_Gates_intro_mine}, our equivalent of \cref{lem_Bracketing_Entropy_Gates_intro_mine_intro}, will be much harder and could be used for controlling the bracketing entropy for many standard MoE regression models with Gaussian gating networks, see \cref{sec_contribution_proof_GLoME} for more details.


\subsubsection{Our contributions on the proof for GLoME models}\label{sec_contribution_proof_GLoME}
To prove \cref{thm_Oracle_Inequality_GLoME}, we also need to make use of \cref{thm_weakOracleInequality_PennecEL_intro}.
Then, our model collection has to satisfy \cref{assumption_K_intro} (K), \cref{assumption_Sep_intro} (Sep), and \cref{assumption_H_intro} (H). 

Note that in the proof of LinBoSGaME models, the authors of \cite{montuelle2014mixture} did not prove \cref{assumption_K_intro} (K) and \cref{assumption_Sep_intro} (Sep) and considered them as assumptions on their LinBoSGaME models because of the complexity of LinBoSGaME models and technical reasons. However, in our proof for GLoME models, we consider an explicit example where the model is defined by $ \cM = \cK \times \cD_\Upsilonb = \left[K_{\max}\right] \times \left[d_{\max}\right]$, $K_{\max},d_{\max} \in \Ns$, leads to the \cref{assumption_K_intro} (K) is always satisfied. It is interesting to find the optimal family $\left(z_\bfm\right)_{\bfm \in \cM}$ satisfying \cref{assumption_K_intro} (K). To the best of our knowledge, this question is only partially answered in some special cases of MoE regression models, \eg~Gaussian graphical model \cite{devijver2018block}, Gaussian finite mixture models \cite{maugis2011non}, finite mixture of Gaussian regression models \cite{devijver2015finite,devijver2017joint}. However, for the standard MoE regression models such as LinBoSGaME and GLoME models, such question still remains open due to the complexity of models. We hope to resolve this important and interesting problem in our future work. 
Furthermore, note that the \cref{assumption_Sep_intro} (Sep) is true when we consider Gaussian densities \cite{massart2007concentration}. 

Therefore, our model has only to satisfy the remaining \cref{assumption_H_intro} (H). Following the same strategy as in the proof of LinBoSGaME models, the main task for GLoME models is to establish \cref{lem_Bracketing_Entropy_Gates_intro_mine}, which is much more difficult compared to \cref{lem_Bracketing_Entropy_Gates_intro_mine_intro} and is proved in \cref{sec_proof_lem_Bracketing_Entropy_Gates_intro_mine}. 

For the Gaussian gating parameters, to make use of the first inequality from \cref{lem_Bracketing_Entropy_Gates_intro_mine_intro} of \cite{montuelle2014mixture}, we propose the following \emph{reparameterization trick} of the Gaussian gating space, which is defined in \eqref{eq_GLoME_GaussianGating}, via the logistics scheme $\cP_K$ and the nonlinear space $\cW_K$ as follows:
\begin{align}\label{eq_reparameterizationTrick_intro}
\cW_K &= \left\{ \yv \mapsto \left(\ln\left(\pi_k\phi_L\left(\yv;\bfc_k,\Gammab_k\right)\right)\right)_{k \in [K]} =: \left(\bfw_{k}(\yv;\omegab)\right)_{k \in [K]}= \bfw\left(\yv;\omegab\right) : \omegab \in \widetilde{\Omegab}_K\right\},\nn\\
\cP_K&= \left\{ \yv \mapsto\left(\frac{e^{\bfw_{k}(\yv)}}{\sum_{l=1}^K e^{\bfw_{l}(\yv)}}\right)_{k\in[K]} =: \left(g_{\bfw,k}\left(\yv\right)\right)_{k \in [K]}, \bfw \in \cW_K\right\}.
\end{align}
We aim to provide the following important upper bound for metric  entropy of nonlinear space  \cref{lem_Bracketing_Entropy_Gates_intro_mine}, which play a key step for controlling the bracketing entropy not only for GLoME models but also for any standard MoE regression models with Gaussian gating networks, \eg~BLLiM and BLoME models, see again \cref{fig_MoE_GaussianGN} for comprehensive descriptions of this general class.
\begin{lem}\label{lem_Bracketing_Entropy_Gates_intro_mine}
For all $\delta \in (0,\sqrt{2}]$, there exists a constant $C_{\cW_K}$ such that
\begin{align*}
	\cH_{[\cdot],\sup_\yv d_k}\left(\frac{\delta}{5},\cP_K\right)
	&	\le 
	\entropy\left(\frac{ 3 \sqrt{3}\delta}{20 \sqrt{K-1}},\cW_K\right) \nn\\
	&\le 
	\dim\left(\cW_K\right) \left(C_{\cW_K} + \ln \left(\frac{ 20 \sqrt{K-1}}{3 \sqrt{3}\delta}\right)\right).
\end{align*}
\end{lem}
More precisely, \cref{assumption_H_intro} (H)
is proved due to \cref{lem_model_complexity_Hellinger} and the following \cref{lem_montuelle_lemma5new} by using the fact that $\sum_{k=1}^K g_{\bfw,k}\left(\yv\right)= 1, \forall \yv \in \cY, \forall \bfw \in \cW_K$ and $\cH_{\left[\cdot\right],\thell}\left(\delta,S_\bfm \right) \le 	\cH_{[\cdot],\sup_\yv d_\xv }\left(\delta,S_\bfm \right)$, which is obtained by definition of bracketing entropy and $\thell(s,t) \le \sup_\yv d_\xv(s,t)$.

\begin{lem}[Lemma 5 from \cite{montuelle2014mixture}]\label{lem_montuelle_lemma5new}
Let
\begin{align*}
	\gkd &= \left\{\cX \times \cY \ni (\xv,\yv) \mapsto \left(\Phi_D\left(\xv;\vb_{k,d}(\yv),\Sigmab_k\right)
	\right)_{k \in [K]}: \upsilonb_d \in \UpsilondK, \Sigmab \in \bfV_K \right\}.
\end{align*}
For all $\delta \in (0,\sqrt{2}]$ and $\bfm \in \cM$,
\begin{align*}
	\cH_{[\cdot],\sup_\yv d_\xv }\left(\delta,S_\bfm \right) \le 	\cH_{[\cdot],\sup_\yv d_k }\left(\frac{\delta}{5},\cP_K\right) + 	\cH_{[\cdot],\sup_\yv\max_k d_\xv }\left(\frac{\delta}{5},\gkd\right).
\end{align*}
\end{lem}
By making use of \cref{lem_montuelle_lemma5new}, the remaining task is to control the bracketing entropy of the Gaussian gating functions and experts separately via \cref{lem_Bracketing_Entropy_Gates_intro_mine,lem_Bracketing_Entropy_Experts}, which are proved in \cref{sec_proof_lem_Bracketing_Entropy_Gates_intro_mine,sec_proof_lem_Bracketing_Entropy_Experts}, respectively.


\begin{lem}\label{lem_Bracketing_Entropy_Experts} For all $\delta \in (0,\sqrt{2}]$, there exists a constant $C_{\gkd}$ such that
\begin{align} \label{eq_Bracketing_Entropy_Experts}
	\cH_{[\cdot],\sup_\yv\max_k d_\xv }\left(\frac{\delta}{5},\gkd\right) \le \dim\left(\gkd\right) \left(C_{\gkd} + \ln \left(\frac{1}{\delta}\right)\right).
\end{align}
\end{lem}

To this end, \cref{lem_montuelle_lemma5new} allows us to conclude that given $\fC=C_{\cW_K} + \ln \left(\frac{5K_{\max}\sqrt{K_{\max}}}{a_W}\right) + C_{\gkd}$,
\begin{align*}
\cH_{[\cdot],\sup_\yv d_\xv }\left(\delta,S_\bfm \right) &\le \cH_{[\cdot],\sup_\yv d_k}\left(\frac{\delta}{5},\cP_K\right) + \cH_{[\cdot],\sup_\yv\max_k d_\xv }\left(\frac{\delta}{5},\gkd\right)\nn\\
&\le \dim(S_\bfm )\left(\fC + \ln \left(\frac{1}{\delta}\right)\right).
\end{align*}
Then, \cref{prop_montuelle_Proposition1_intro} leads to
\begin{align*}
n \delta^2_\bfm \le   \dim(S_\bfm ) \left(2 \left(\sqrt{\fC} + \sqrt{\pi}\right)^2 + \left(\ln \frac{n}{\left(\sqrt{\fC} + \sqrt{\pi}\right)^2\dim\left(S_\bfm \right)}\right)_+\right).
\end{align*}
Finally, \cref{thm_weakOracleInequality_PennecEL_intro} implies that for any given collection of GLoME models $\left(S_\bfm \right)_{\bfm \in \cM}$,
the oracle inequality of \cref{thm_Oracle_Inequality_GLoME} is satisfied if $\pen(\bfm)$ is bounded from below by
\begin{align*}
\kappa \left(\dim(S_\bfm ) \left(2 \left(\sqrt{\fC} + \sqrt{\pi}\right)^2 + \left(\ln \frac{n}{\left(\sqrt{\fC} + \sqrt{\pi}\right)^2\dim\left(S_\bfm \right)}\right)_+\right)+ z_\bfm\right).
\end{align*}

\subsection{Proof of random subcollection of BLoME models  }\label{sec_contribution_collection_random_MoE}
The random subcollection of MoE models include BLLiM and BLoME models where finite-sample oracle inequalities were only well studied for finite mixture of Gaussian regression models via a model selection theorem for MLE among a random subcollection of models in regression framework of \citep[Theorem 5.1]{devijver2015finite}, see also \citep[Theorem 7.3]{devijver2018block}. This is an extension of a whole collection of conditional densities from \citep[Theorem 2]{cohen2011conditional}, and of \citep[Theorem 7.11]{massart2007concentration}, working only for density estimation. 
In \cref{sec_general_model_selection_random_intro}, we first summarize this theorem and the techniques that \cite{devijver2017joint} used to control the bracketing entropy of finite mixture of Gaussian regression models with joint rank and variable selection for parsimonious estimation in a high-dimensional framework. Then, we explain why such techniques can not be applied to our collection of BLLiM and BLoME models to highlight the main challenges and our contributions.

\subsubsection{A model selection theorem for MLE among a random subcollection}\label{sec_general_model_selection_random_intro}
We can now state the main result of \citep[Theorem 5.1]{devijver2015finite} for the model selection theorem for MLE among a random subcollection.
\begin{theorem}[Theorem 5.1 from \cite{devijver2015finite}]\label{thm_5_1_devijver2015finite_intro}
Let $\left(\xv_{[n]},\yv_{[n]}\right)$ be observations coming from an unknown conditional density $s_0$. Let the model collection $\left(S_\bfm\right)_{\bfm \in \cM}$ be an at most countable collection of conditional density sets. Assume that \cref{assumption_K_intro} (K), \cref{assumption_Sep_intro} (Sep), and \cref{assumption_H_intro} (H) hold for every $\bfm \in \cM$. Let $\epsilon_{KL} > 0$, and $\bar{s}_\bfm \in S_\bfm$, such that
\begin{align*}
	\tkl\left(s_0,\bar{s}_\bfm\right) \le \inf_{t \in  S_\bfm} \tkl\left(s_0,t\right) + \frac{\epsilon_{KL}}{n};
\end{align*}
and let $\tau > 0$, such that
\begin{align}\label{eq_assumptionRS_intro}
	\bar{s}_\bfm\ge e^{-\tau} s_0.
\end{align}
Introduce $\left(S_\bfm\right)_{\bfm \in \widetilde{\cM}}$, a random subcollection of $\left(S_\bfm\right)_{\bfm \in \cM}$.
Consider the collection $\left(\widehat{s}_\bfm\right)_{\bfm \in \widetilde{\cM}}$ of $\eta$-log likelihood minimizer satisfying \eqref{eq_define_NLL} for all $\bfm \in \widetilde{\cM}$.
Then, for any $\rho \in (0,1)$, and any $C_1 > 1$, there are two constants $\kappa$ and $C_2$ depending only on $\rho$ and $C_1$, such that, for every index $\bfm \in \cM$,
\begin{align*}
	\pen(\bfm) \ge \kappa \left(\cD_\bfm + (1 \vee \tau)z_\bfm\right),
\end{align*}
and where the model complexity $\cD_\bfm$ is defined in \cref{assumption_H_intro} (H),
the $\eta'$-penalized likelihood estimator $\widehat{s}_{\widehat{\bfm}}$, defined as in \eqref{eq_define_penalizedNLL} on the subset $\widetilde{\cM} \subset \cM$,
satisfies
\begin{align*}
	\E{\Xv_{[n]},\Yv_{[n]}}{\jtkl\left(s_0,\widehat{s}_{\widehat{\bfm}}\right)}&\le C_1\E{\Xv_{[n]},\Yv_{[n]}}{\inf_{\bfm \in \widetilde{\cM}}\left(\inf_{t \in S_\bfm} \tkl\left(s_0,t\right)+2 \frac{\pen(\bfm)}{n}\right)} \nn\\
	& \quad + C_2(1 \vee \tau)\frac{\Xi^2}{n} + \frac{\eta' + \eta}{n}.
\end{align*}
\end{theorem}
\subsubsection{Sketch of the proof for BLoME models}\label{sec_contribution_proof_PSGaBloME}
To work with conditional density estimation in the BLLiM and BLoME models, it is natural to make use of \cref{thm_weakOracleInequality_PennecEL_intro}. However, it is worth mentioning that because the model collection constructed by the BLLiM \cite{devijver2017nonlinear} or by some suitable procedures for BLoME models in practice is usually random, we have to use a model selection theorem for MLE among a random subcollection (cf. \citep[Theorem 5.1]{devijver2015finite} and \citep[Theorem 7.3]{devijver2018block}). 


More precisely, we explain how \cref{thm_5_1_devijver2015finite_intro} implies the finite-sample oracle inequality, \cref{thm_Oracle_Inequality_BLoME}. To this end, our collections of BLoME models has to satisfy some regularity assumptions, 
see \cref{proofLemma_nguyen2021nonBLoME_intro} for more details. For BLoME models, the main difficulty in proving our oracle inequality lies in bounding the bracketing entropy of the Gaussian gating functions and Gaussian experts with block-diagonal covariance matrices. To overcome the former issue, we follow a reparameterization trick of the Gaussian gating parameters space in \eqref{eq_reparameterizationTrick_intro} and \cref{lem_Bracketing_Entropy_Gates_intro_mine}. For the second one, based on some ideas of Gaussian mixture models from \cite{genovese2000rates,maugis2011non},  we contribute a novel extension for standard MoE models with Gaussian gating networks. Note that our contributions extend the recent novel result on block-diagonal covariance matrices in \cite{devijver2018block}, which is only developed for Gaussian graphical models.

\subsubsection{Our contributions on BLoME models}\label{proofLemma_nguyen2021nonBLoME_intro}
It should be stressed that all we need is to verify that \cref{assumption_H_intro} (H), \cref{assumption_Sep_intro} (Sep) and \cref{assumption_K_intro} (K) hold for every $\bfm \in \cM$. According to the result from \citep[Section 5.3]{devijver2015finite}, \cref{assumption_Sep_intro} (Sep) holds when we consider Gaussian densities and the assumption defined by \eqref{eq_assumptionRS_intro} is true if we assume further that the true conditional density $s_0$ is bounded and compactly supported. Furthermore, since we restricted $d$ and $K$ to $\cD_\Upsilonb = \left[d_{\max}\right]$ and $\cK = \left[K_{\max}\right]$, respectively, it is true that there exists a family $\left(z_\bfm\right)_{\bfm \in \cM}$ and $\Xi >0$ such that, \cref{assumption_K_intro} (K) is satisfied.
Therefore, the proof for the remaining \cref{assumption_H_intro} (H) is our novel contribution. In particular, to the best of our knowledge, there are no results that can be directly applied to \cref{assumption_H_intro} (H) for BLoME models due to their complexity with the Gaussian gating functions and Gaussian experts with block-diagonal covariance matrices. This highlights our contributions to this challenging problem concerning the work of \cite{genovese2000rates,maugis2011non,devijver2015finite,devijver2017joint,devijver2018block,montuelle2014mixture}.

By using \cref{lem_model_complexity_Hellinger}, \cref{assumption_H_intro} (H) holds true if we can prove that for any $\delta \in (0,\sqrt{2}]$, the collection of LinBoSGaME models, $\cS = \left(S_\bfm \right)_{\bfm \in \cM}$, satisfies
\begin{align}
\cH_{\left[\cdot\right],\thell}\left(\delta,S_\bfm \right) \le \dim(S_\bfm )\left(C_\bfm+\ln\left(\frac{1}{\delta}\right)\right). \label{eq_bractketingEntropyModelSm_intro}
\end{align}
\begin{proof}[Proof of \ref{eq_bractketingEntropyModelSm_intro}]
Note that \eqref{eq_bractketingEntropyModelSm_intro} can be established by first decomposing the entropy term between the Gaussian gating functions and the Gaussian experts. 
Motivated by our \emph{reparameterization trick} of the Gaussian gating space in \cref{sec_contribution_proof_GLoME} and \eqref{eq_reparameterizationTrick_intro}, we define Gaussian experts $\gkdb$ as follows.
\begin{align*}
	%
	\gkdb &= \left\{ (\xv,\yv) \mapsto \left(\phi\left(\xv;\upsilonb_{k,d}(\yv),\Sigmab_k\left(\bfB_k\right)\right)
	\right)_{k \in [K]}: \upsilonb_{d} \in \UpsilondK, \Sigmab(\bfB) \in \bfV_K(\bfB) \right\}.
\end{align*}
There are two possible ways to decompose the bracketing entropy of $S_\bfm$ based on different distances.
For the first approach, we can use  \cref{lem.bracketingEntropyDecomposition_intro} \citep[Lemma 5]{montuelle2014mixture}:
\begin{lem}\label{lem.bracketingEntropyDecomposition_intro}
	For all $\delta \in (0,\sqrt{2}]$ and $\bfm \in \cM$,
	\begin{align*}
		\cH_{[\cdot],\sup_\yv d_\xv }\left(\delta,S_\bfm\right) \le 	\cH_{[\cdot],\sup_\yv d_k }\left(\frac{\delta}{5},\cP_K\right) + 	\cH_{[\cdot],\sup_\yv\max_k d_\xv }\left(\frac{\delta}{5},\gkdb\right).
	\end{align*}
\end{lem}
As mentioning in Appendix B.2.1 from \cite{montuelle2014mixture}, \cref{lem.bracketingEntropyDecomposition_intro} boils down to assuming that $\Yv$ is bounded. Furthermore, they also claim that this boundedness assumption can be relaxed when using smaller distance $\thell$ but bounding the corresponding bracketing entropy becomes much more challenging. 
We successfully weaken such boundedness assumption via utilizing the smaller distance: $\thell$, for the bracketing entropy of $S_\bfm$ although bounding such bracketing entropy for $\cW_K$ and $\cG_{K,\cB}$ becomes much more challenging. This reinforces our new contributions concerning the control of bracketing entropy of BLoME models.
Consequently, this leads to the second approach via \cref{lem_bracketingEntropyDecomposition2_intro} \citep[Lemma 6]{montuelle2014mixture}.
\begin{lem}\label{lem_bracketingEntropyDecomposition2_intro}
	For all $\delta \in (0,\sqrt{2}]$,
	\begin{align*}
		\cH_{[\cdot],\thell }\left(\delta,S_\bfm\right) \le \cH_{[\cdot],d_{\cP_K} }\left(\frac{\delta}{2},\cP_K\right) + 	\cH_{[\cdot], d_{\gkdb} }\left(\frac{\delta}{2},\gkdb\right),
	\end{align*}
	where
	\begin{align*}
		d^2_{\cP_K}\left(g^+,g^-\right) &= \E{\Yv_{[n]}}{\frac{1}{n}\sum_{i=1}^n d^2_k\left(g^+ \left(\Yv_i\right),g^-(\Yv_i)\right)} \nn\\
		&= \E{\Yv_{[n]}}{\frac{1}{n}\sum_{i=1}^n \sum_{k=1}^K \left(\sqrt{g^+_k \left(\Yv_i\right)}-\sqrt{g^-_k \left(\Yv_i\right)}\right)^2},\\
		d^2_{\gkdb}\left(\phi^+,\phi^-\right) &= \E{\Yv_{[n]}}{\frac{1}{n}\sum_{i=1}^n \sum_{k=1}^K d^2_\xv\left(\phi^+_k \left(\cdot,\Yv_i\right),\phi^-_k\left(\cdot,\Yv_i\right)\right)}
		\nn\\
		&= \E{\Yv_{[n]}}{\frac{1}{n}\sum_{i=1}^n \sum_{k=1}^K \int_{\cX}\left(\sqrt{\phi^+_k \left(\xv,\Yv_i\right)}-\sqrt{\phi^+_k \left(\xv,\Yv_i\right)}\right)^2 d\xv}.
	\end{align*}
\end{lem}
Next, we make use of \cref{lem_Inequality_Hellinger_Supx_dy_dk_intro}, which is proved in \cref{sec_lem_Inequality_Hellinger_Supx_dy_dk_intro}, to provide an upper bound on the bracketing entropy of $S_\bfm$ and $\cP_K$ on the corresponding distances $\thell$ and $d_{\cP_K}$, respectively.
\begin{lem}\label{lem_Inequality_Hellinger_Supx_dy_dk_intro}
	It holds that
	\begin{align}
		\thell(s,t) &\le \sup_\yv d_\xv(s,t) \text{, and }
		\cH_{[\cdot],\thell}\left(\delta,S_\bfm\right) \le \cH_{[\cdot],\sup_\yv d_\xv }\left(\delta,S_\bfm\right), \label{eq_Inequality_Hellinger_Supx_dy_intro}\\
		d_{\cP_K}\left(g^+,g^-\right)&\le \sup_\yv d_k(g^+,g^-) \text{, and } \cH_{[\cdot],d_{\cP_K} }\left(\frac{\delta}{2},\cP_K\right) \le \cH_{[\cdot],\sup_\yv d_k }\left(\frac{\delta}{2},\cP_K\right). \label{eq_Inequality_Hellinger_Supx_dk_intro}
	\end{align}
\end{lem}
%
%
\cref{lem_bracketingEntropyDecomposition2_intro,lem_Inequality_Hellinger_Supx_dy_dk_intro} imply that
\begin{align*}
	\cH_{[\cdot],\thell}\left(\delta,S_\bfm\right) \le \cH_{[\cdot],\sup_\yv d_k }\left(\frac{\delta}{2},\cP_K\right) + 	\cH_{[\cdot],d_{\gkdb} }\left(\frac{\delta}{2},\gkdb\right).
\end{align*}
Based on this metric, one can first relate the bracketing entropy of $\cP_{K}$ to $\entropy(\delta,\cW_K)$, and then obtain the upper bound for its entropy via \cref{lem_Bracketing_Entropy_Gates_intro_mine}.
Then, we present our main contribution for BLoME models via
\cref{lem_bracketingEntropyGaussianBlock_intro}. This lemma allows us to construct the Gaussian brackets to handle the metric entropy for Gaussian experts, which is established in \cref{sec_proof_lem_bracketingEntropyGaussianBlock_intro}.
\begin{lem}\label{lem_bracketingEntropyGaussianBlock_intro}
	\begin{align}\label{eq.bracketingEntropyGaussianBlock_intro}
		\cH_{[\cdot],d_{\gkdb} }\left(\frac{\delta}{2},\gkdb\right) \le \dim\left(\gkdb\right) \left(C_{\gkdb} + \ln \left(\frac{1}{\delta}\right)\right).
	\end{align}
\end{lem}
Finally, \eqref{eq_bractketingEntropyModelSm_intro} can easily proved via \cref{lem_Bracketing_Entropy_Gates_intro_mine,lem_bracketingEntropyGaussianBlock_intro}.
\end{proof}
%

\section{Numerical experiments} \label{sec_numerical_Experiment}
Note that our numerical experiments in \cref{sec_numerical_Experiment}, and the codes written in the {\bf R} programming language \citep{team2020r} are available in the following link:\\ \href{https://github.com/Trung-TinNGUYEN/NamsGLoME-Simulation}{https://github.com/Trung-TinNGUYEN/NamsGLoME-Simulation}.

\subsection{The procedure}
We illustrate our theoretical results in settings similar to those considered by \cite{chamroukhi2010hidden} and \cite{montuelle2014mixture}, including simulated as well as real data sets. We first observe $n$ random samples $\left(\xv_i,\yv_i\right)_{i\in[n]}$ from an forward conditional density $s^*_0$, and look for the best estimate among $s^*_\bfm \in S^*_\bfm, \bfm \in \cM$, defined in \eqref{eq_define_GLoME_forward_Boundedness}. We considered the simple case where the mean experts are linear functions, which leads to GLoME and supervised GLLiM are identical models.
Our aim is to estimate the best number of components $K$, as well as the model parameters. As described in more detail in \cite{deleforge2015high}, we use a GLLiM-EM algorithm to estimate the model parameters for each $K$, and select the optimal model using the penalized approach that was described earlier. More precisely, in the following numerical experiments, the GLoME model is learned using functions from the package \href{https://cran.r-project.org/web/packages/xLLiM/index.html}{xLLiM}, available on CRAN. It solves the inverse regression problem, defined in \eqref{eq_inverseGLLiM}, and obtain the inverse maximum likelihood estimates (MLE) $\left(\widehat{s}_\bfm \left(\xv_i\mid \yv_i\right)\right)_{i \in [N]}$, $\bfm \in \cM$, then via \eqref{eq_inverse_regression_mapping}, we obtain the forward MLE $\left(\widehat{s}_\bfm^* \left(\yv_i\mid \xv_i\right)\right)_{i \in [N]}$, $\bfm \in \cM$.

According to the general procedure for model selection, we first compute the forward MLE for each model $\bfm \in \cM$, where $\cM = \cK$. Then, we select the model that satisfies the definition \eqref{eq_define_penalizedNLL} with $\pen(\bfm) = \kappa \dim(S_\bfm ^*)$, where $\kappa$ is a positive hyperparameter.
The point we want to stress here is that the theoretical result stated in \cref{thm_Oracle_Inequality_GLoME} 
gives the general form of penalty functions but it
does not provide explicit penalties since the lower bounds on penalty functions in \eqref{eq_lowerBoundPenalty} 
are defined up to an unknown multiplicative
constant $\kappa$ and mixture parameters are not bounded in practice. Nevertheless, \cref{thm_Oracle_Inequality_GLoME} is required to justify the shape of penalties using for the slope heuristic and allows the number of variables $D$ as well as responses $L$ to be large compared to the fixed sample size $n$. In particular, our \cref{thm_Oracle_Inequality_GLoME} and \cref{remark_eq_define_GLoME_Boundedness} guarantee that there exists a $\kappa$ large enough for which the estimate has the desired properties. Therefore, we need a data-driven method to choose $\kappa$.
%
%
According to the AIC or the BIC, we can select $\kappa =1$ or $\kappa = \frac{\ln n}{2}$. 
An important limitation of these criteria, however, is that they are only valid asymptotically. This implies that there are no finite-sample guarantees when using AIC or BIC for choosing between different levels of complexity. Their use in small sample settings is thus ad hoc.
%
To overcome such difficulties, Birg\'e~and Massart \cite{birge2007minimal} proposed a novel approach, called slope heuristics,
supported by a non-asymptotic oracle inequality. This method leads to an optimal data-driven choice of multiplicative constants for penalties, \eg~$\kappa$ in our framework.
%
Thus, we shall concentrate our attention on the slope heuristic for choosing the number of mixture components in our numerical experiments.

%

\subsection{Simulated data sets}
Note that our main objective here is to investigate how well the empirical tensorized Kullback--Leibler divergence between the true model ($s_0^*$) and the selected model $\widehat{s}_{\widehat{\bfm}}^*$ follows the finite-sample oracle inequality of \cref{thm_Oracle_Inequality_GLoME}, as well as the rate of convergence of the error term. 
Therefore, we focus on  $1$-dimensional data sets, that is, with $L = D =1$.
Beyond the statistical estimation and model selection objectives considered here, the dimensionality reduction capability of GLLiM in high-dimensional regression data, typically $D \gg L$, can be found in \citep[Section 6]{deleforge2015high}.

\begin{figure}
\centering
\begin{subfigure}{.25\textheight}
	\centering
	\includegraphics[scale = 0.3]{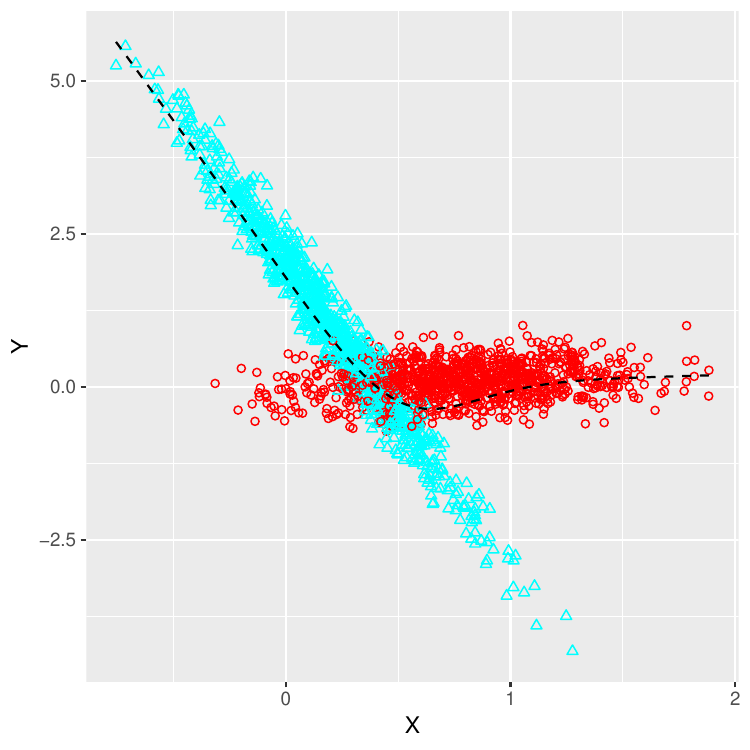}
	\caption{Typical realization of example WS}
	\label{subfig_WS_2000_TypicalRealizations_NS}
\end{subfigure}
\hspace{1cm}
\begin{subfigure}{.25\textheight}
	\centering
	\includegraphics[scale = 0.3]{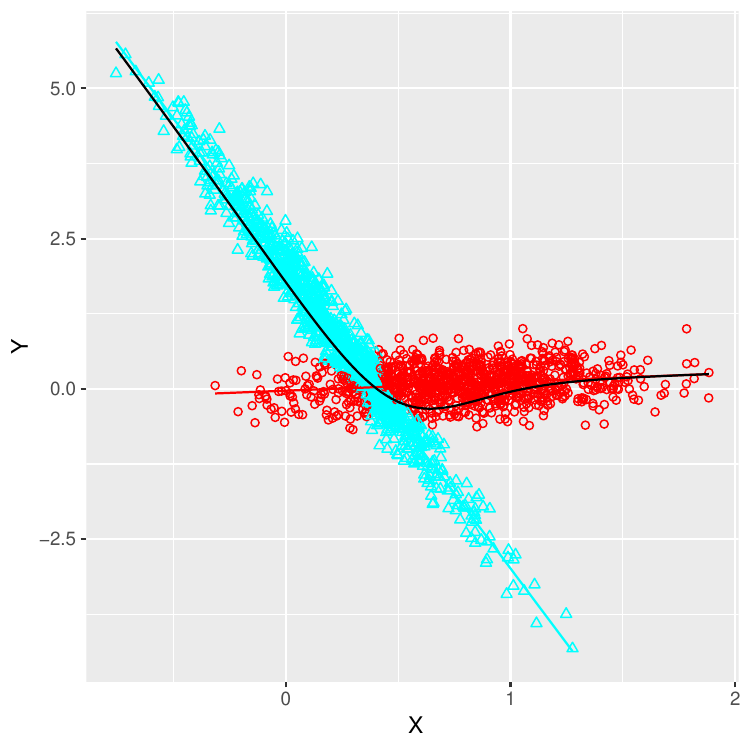}
	\caption{Clustering by GLoME}
	\label{subfig_Clustering_Regression_WS}
\end{subfigure}

\begin{subfigure}{.25\textheight}
	\centering
	\includegraphics[scale = 0.3]{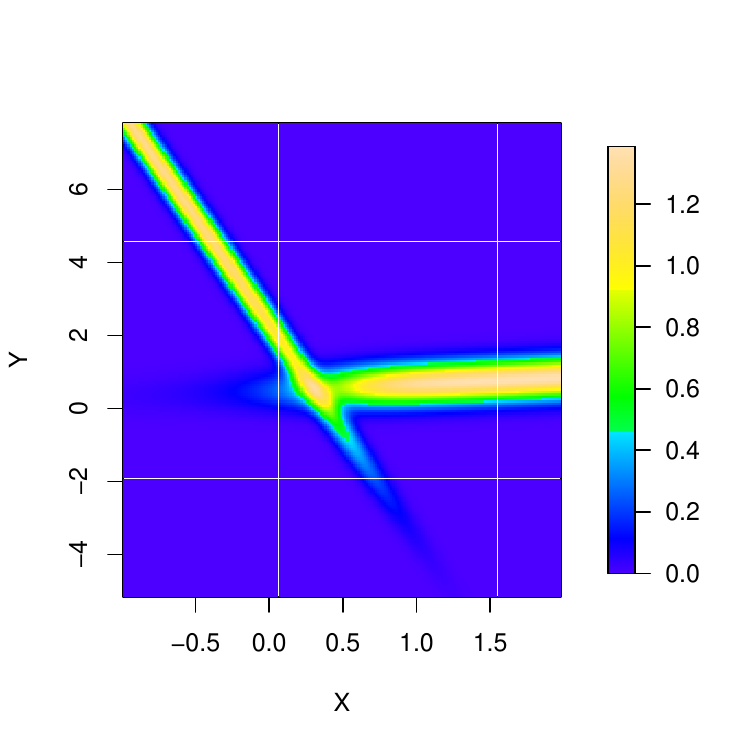}
	\caption{2D view of the resulting conditional density with the 2 regression components}
	\label{subfig_Conditional_density_2D_WS}
\end{subfigure}
\hspace{1cm}
\begin{subfigure}{.25\textheight}
	\centering
	\includegraphics[scale = 0.3]{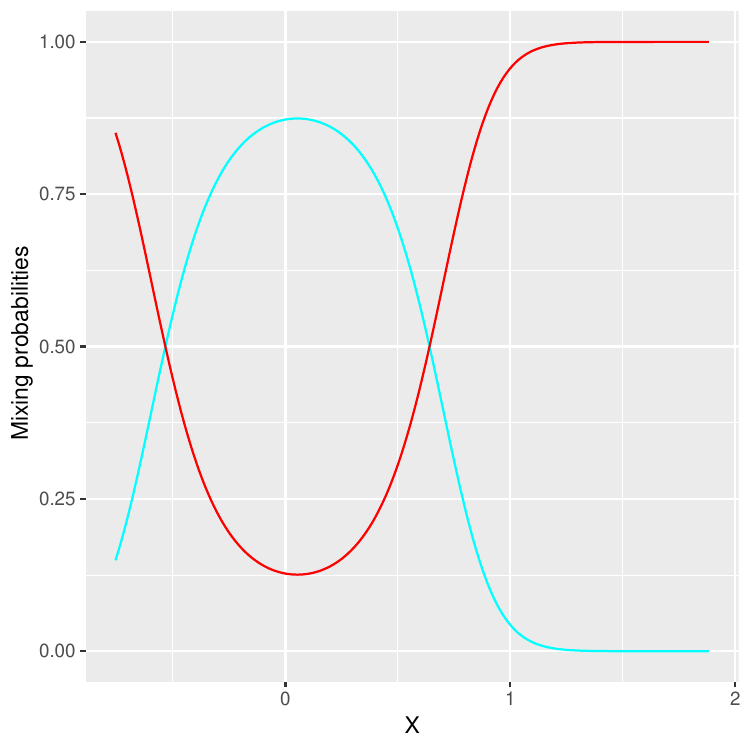}
	\caption{Gating network probabilities}
	\label{subfig_Clustering_Posterior_WS}
\end{subfigure}%

\begin{subfigure}{.25\textheight}
	\centering
	\includegraphics[scale = 0.3]{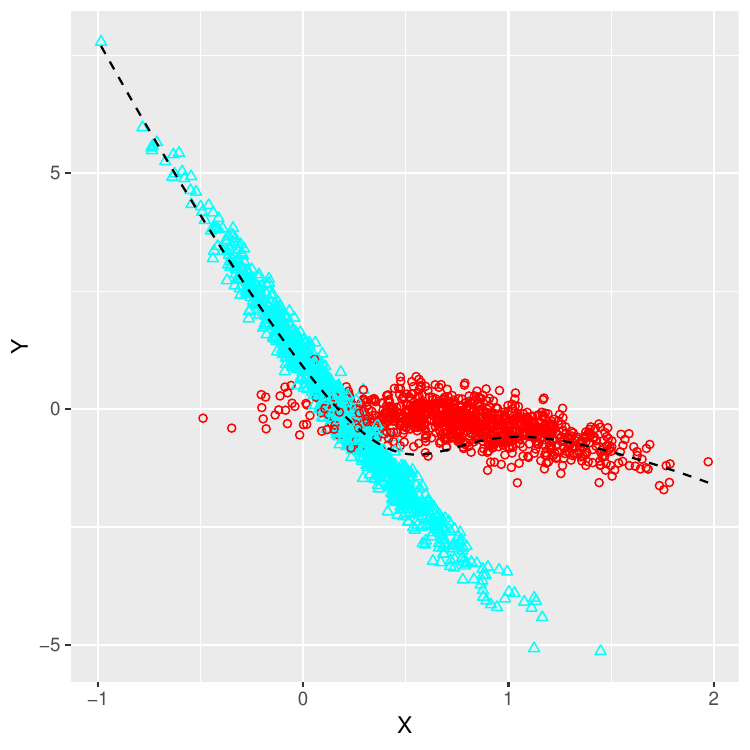}
	\caption{Typical realization of example MS}
	\label{subfig_MS_2000_TypicalRealizations_NS}
\end{subfigure}
\hspace{1cm}
\begin{subfigure}{.25\textheight}
	\centering
	\includegraphics[scale = 0.3]{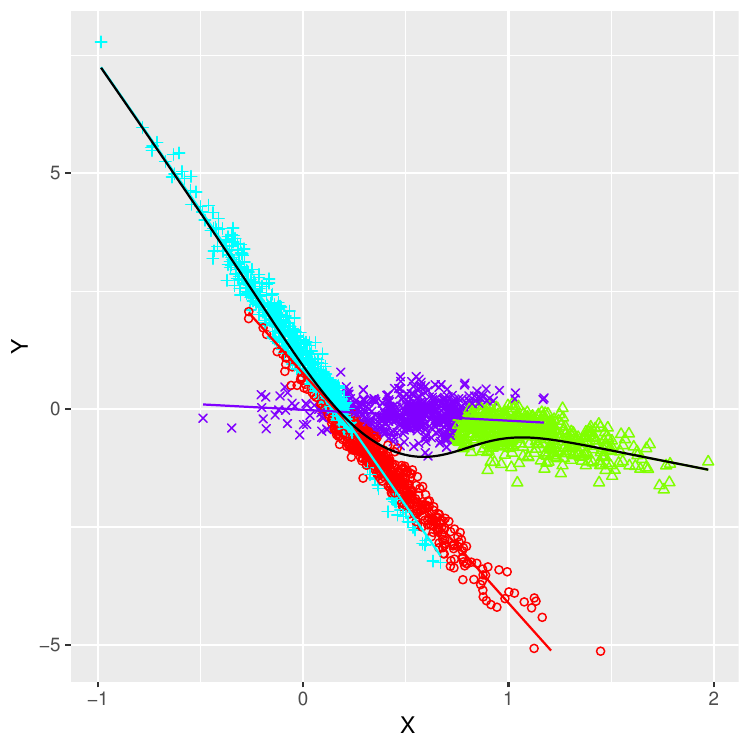}
	\caption{Clustering by GLoME}
	\label{subfig_Clustering_Regression_MS}
\end{subfigure}

\begin{subfigure}{.25\textheight}
	\centering
	\includegraphics[scale = 0.3]{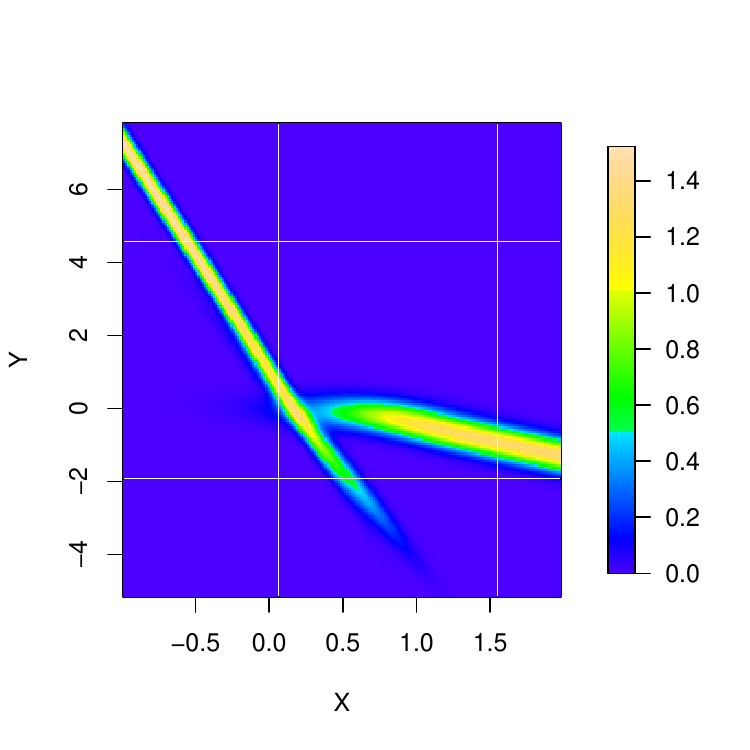}
	\caption{2D view of the resulting conditional density with the 4 regression components}
	\label{subfig_Conditional_density_2D_MS}
\end{subfigure}
\hspace{1cm}
\begin{subfigure}{.25\textheight}
	\centering
	\includegraphics[scale = 0.3]{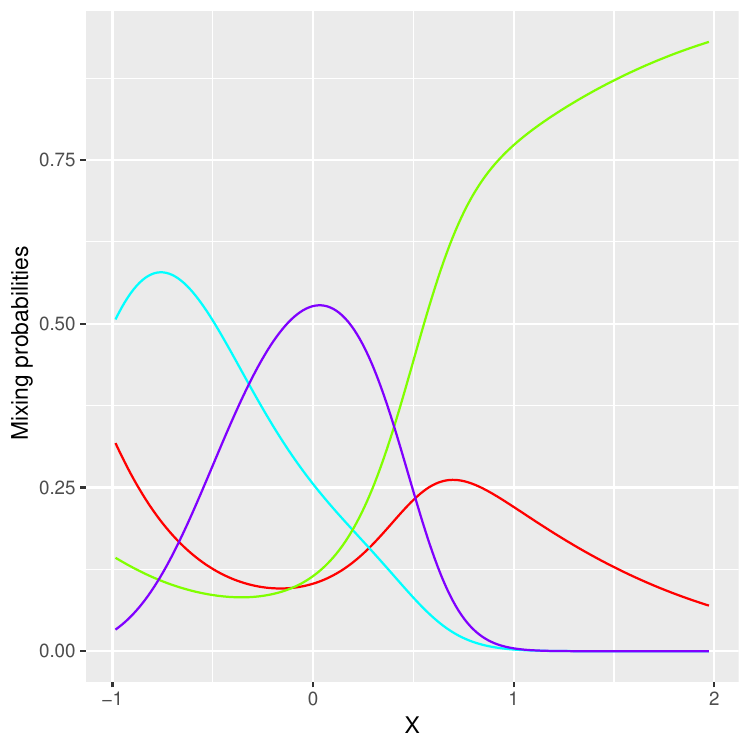}
	\caption{Gating network probabilities}
	\label{subfig_Clustering_Posterior_MS}
\end{subfigure}%

\caption{Clustering deduced from the estimated conditional density of GLoME by a MAP principle with 2000 data points of example WS and MS. The dash and solid black curves present the true and estimated mean functions. }
\label{fig_Estimation_WS_MS}
\end{figure}

We construct simulated data sets following two scenarios: a \emph{well-specified} (WS) case in which the true forward conditional density belongs to the class of proposed models:
\begin{align*}
s_0^*(y\mid  x) &= \frac{\Phi(x;0.2,0.1)}{\Phi(x;0.2,0.1)+\Phi(x;0.8,0.15)}\Phi(y;-5x+2,0.09) \nn\\
& \quad + \frac{\Phi(x;0.8,0.15)}{\Phi(x;0.2,0.1)+\Phi(x;0.8,0.15)}\Phi(y;0.1x,0.09),
\end{align*}
and a \emph{misspecified} (MS) case, whereupon such an assumption is not true:
\begin{align*}
s_0^*(y\mid  x) &= \frac{\Phi(x;0.2,0.1)}{\Phi(x;0.2,0.1)+\Phi(x;0.8,0.15)}\Phi\left(y;x^2-6x+1,0.09\right) \nn\\
& \quad + \frac{\Phi(x;0.8,0.15)}{\Phi(x;0.2,0.1)+\Phi(x;0.8,0.15)}\Phi\left(y;-0.4x^2,0.09\right).
\end{align*}
\cref{subfig_WS_2000_TypicalRealizations_NS,subfig_MS_2000_TypicalRealizations_NS} show some typical realizations of 2000 data points arising from the WS and MS scenarios. Note that by using GLoME, our estimator performs well in the WS setting (\crefrange{subfig_Clustering_Regression_WS}{subfig_Clustering_Posterior_WS}).
In the MS case, we expect our algorithm to automatically balance the model bias and its variance (\crefrange{subfig_Clustering_Regression_MS}{subfig_Clustering_Posterior_MS}), which leads to the choice of a complex model, with $4$ mixture components. This observation will be elaborated upon in the subsequent experiments.

\begin{figure}
\begin{subfigure}{.5\textheight}
	\centering
	\includegraphics[trim={0 0 11cm 0},clip,width=1.1\linewidth]{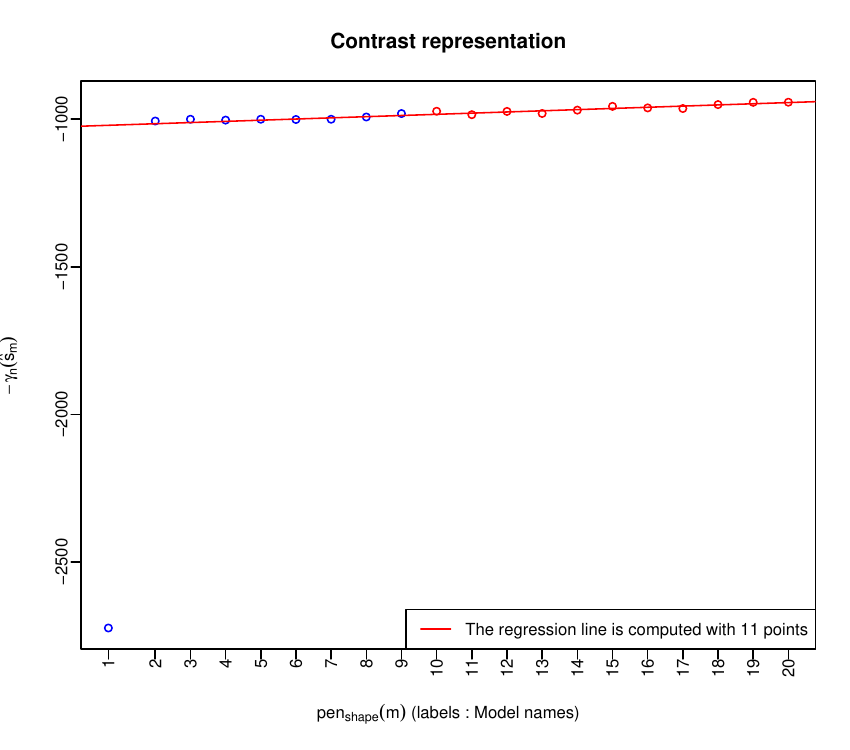}
	\caption{Example WS.}
	\label{fig:sub_WS_2000_DDSE_Slope_CNLL}
\end{subfigure}
\begin{subfigure}{.5\textheight}
	\centering
\includegraphics[trim={0 0 11cm 0},clip,width=1.1\linewidth]{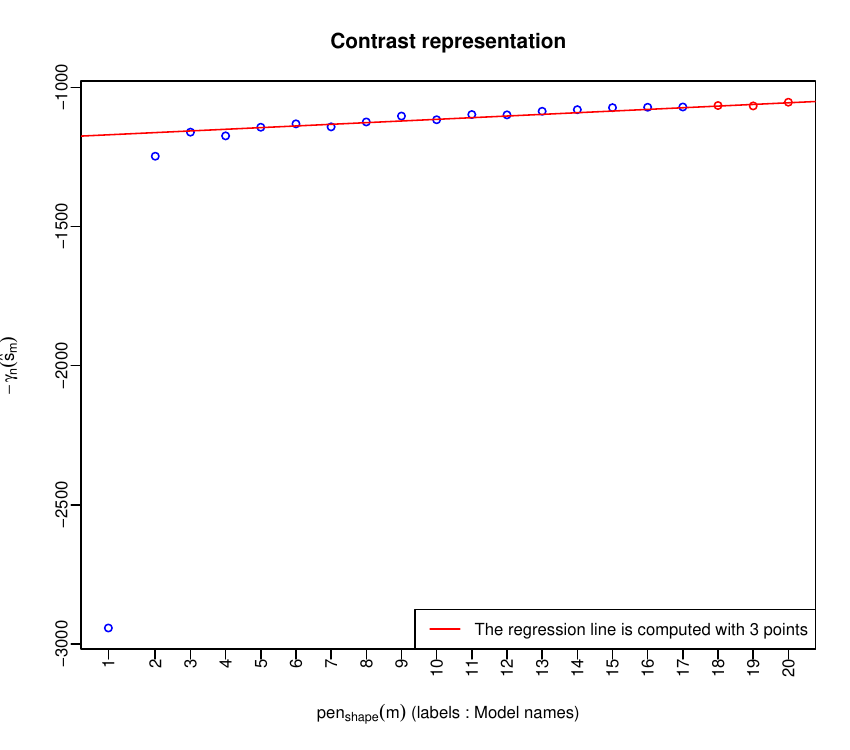}
	\caption{Example MS.}
	\label{fig:sub_MS_2000_DDSE_Slope_CNLL}
\end{subfigure}

\caption{Plot of the selected model dimension using the slope criterion with 2000 data points.}
\label{fig_Slope_Heuristic_2000_DDSE_Slope_CNLL}
\end{figure}

\begin{figure}
	\begin{subfigure}{.5\textheight}
		\centering
		\includegraphics[trim={0 0 11cm 0},clip,width=1.1\linewidth]{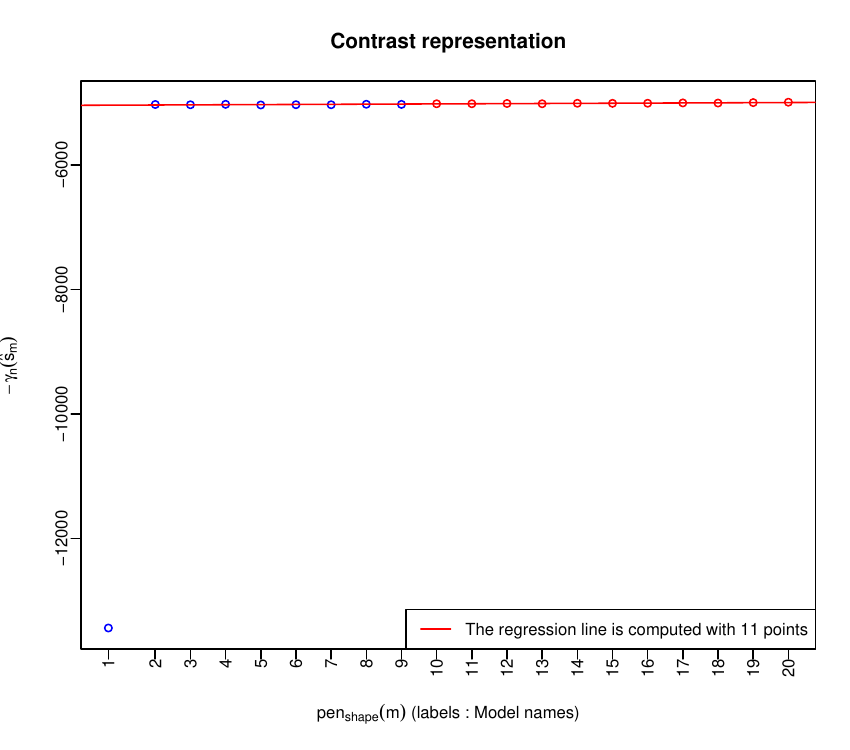}
		\caption{Example WS.}
		\label{fig:sub_WS_10000_DDSE_Slope_CNLL}
	\end{subfigure}
	\begin{subfigure}{.5\textheight}
		\centering
		\includegraphics[trim={0 0 11cm 0},clip,width=1.1\linewidth]{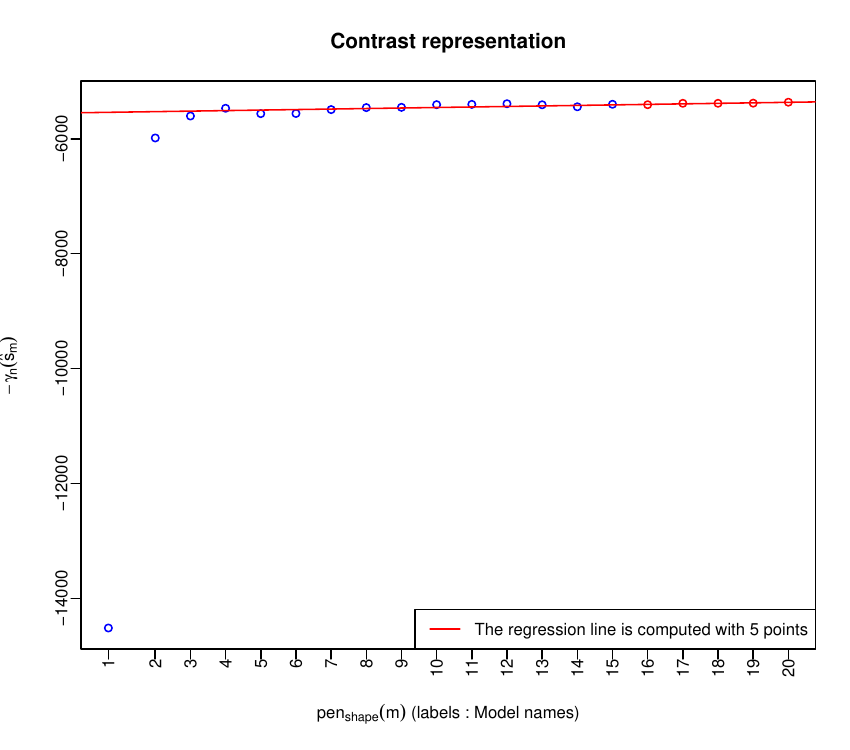}
		\caption{Example MS.}
		\label{fig:sub_MS_10000_DDSE_Slope_CNLL}
	\end{subfigure}
	\caption{Plot of the selected model dimension using the slope criterion with 10000 data points.}
	\label{fig_Slope_Heuristic_10000_DDSE_Slope_CNLL}
\end{figure}

\begin{figure}
\centering
\begin{tabular}{c}
	\includegraphics[height=0.55\textheight]{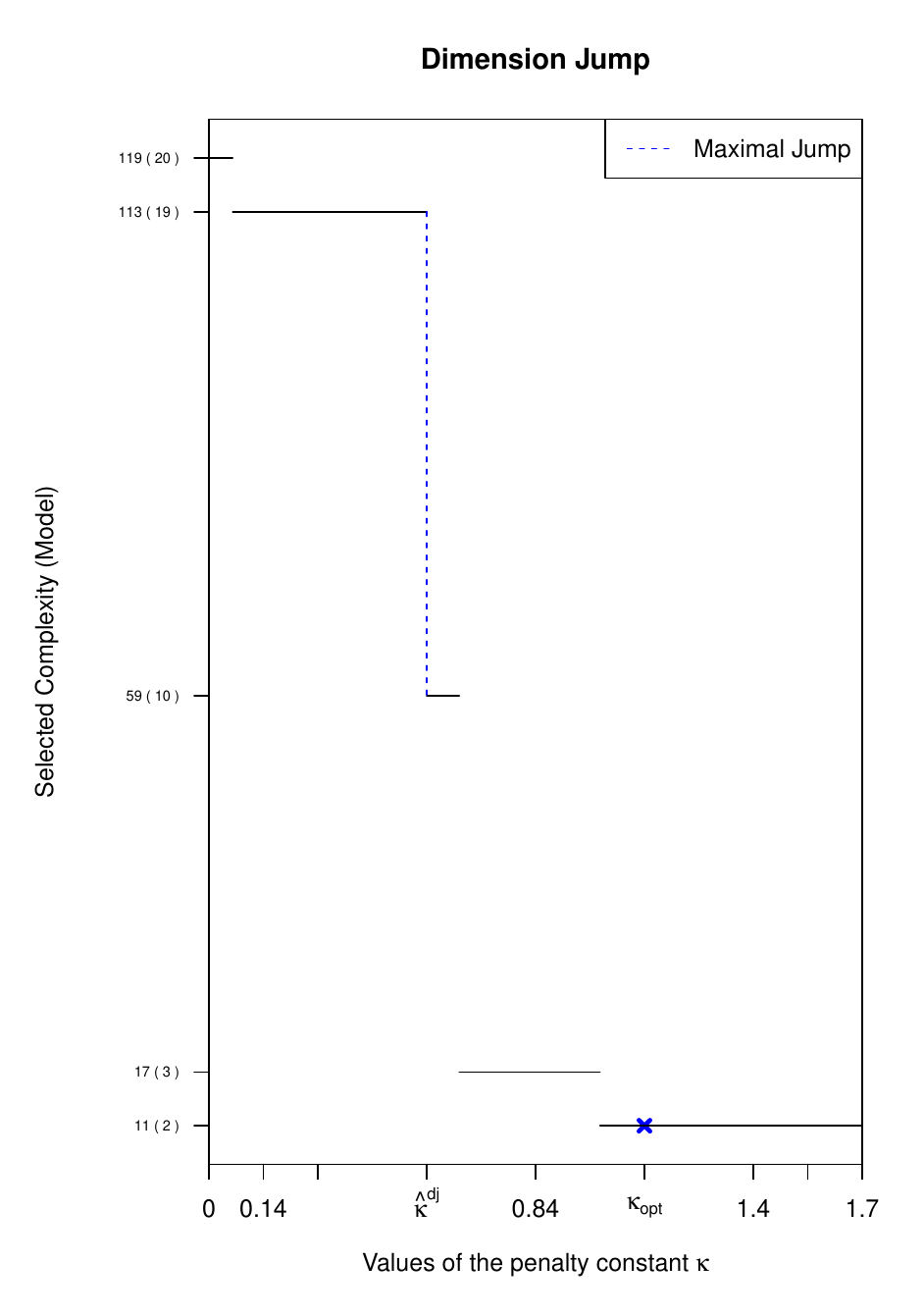}\\
	(a) Example WS.\\
	\includegraphics[height=0.55\textheight]{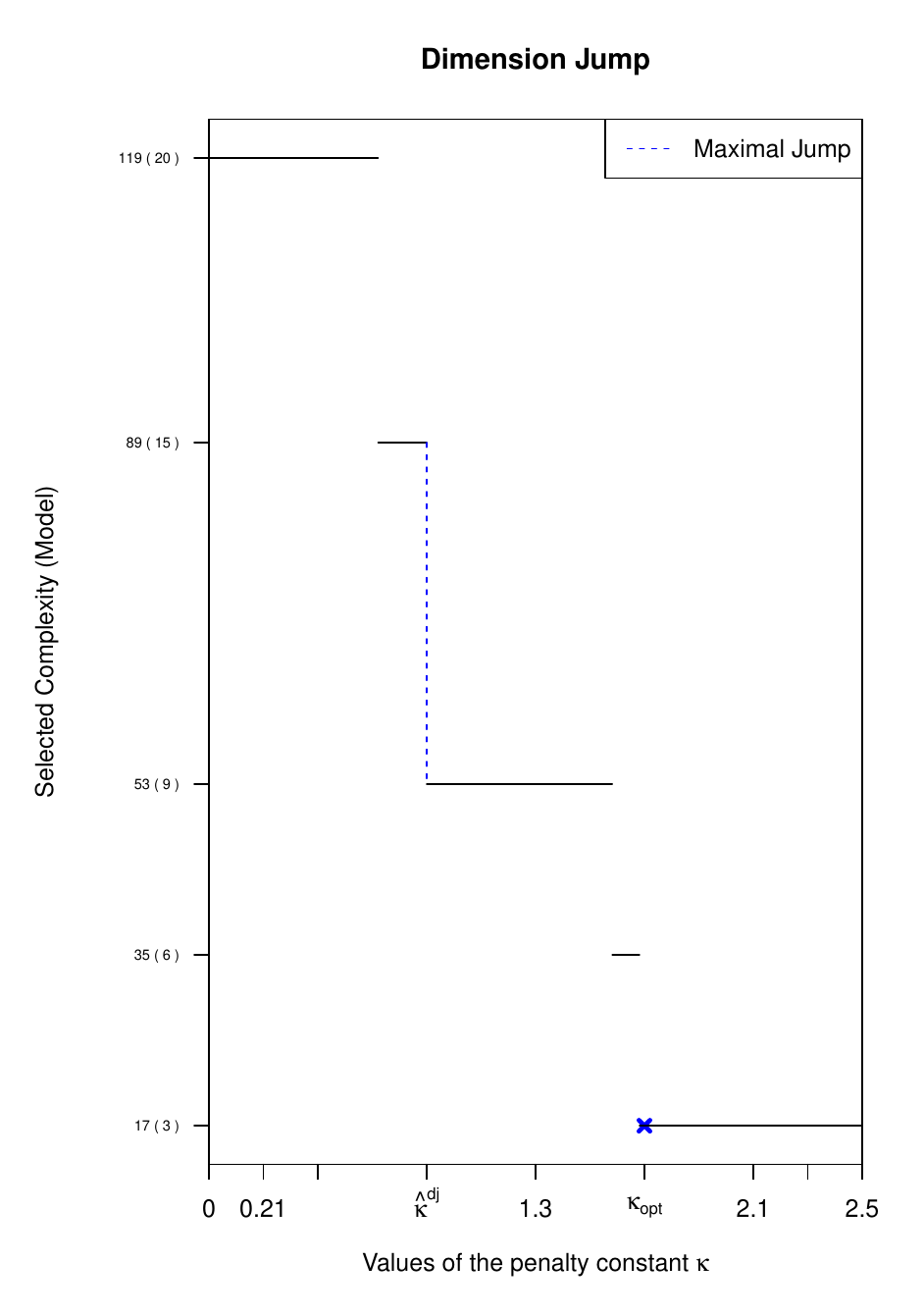} \\
	(b) Example MS.
	%
\end{tabular}
\caption{Plot of the selected model dimension using the jump criterion with 2000 data points.}
\label{fig_Slope_Heuristic_2000_Djump_Slope_CNLL}
\end{figure}

\begin{figure}
	\centering
	\begin{tabular}{c}
		\includegraphics[height=0.55\textheight]{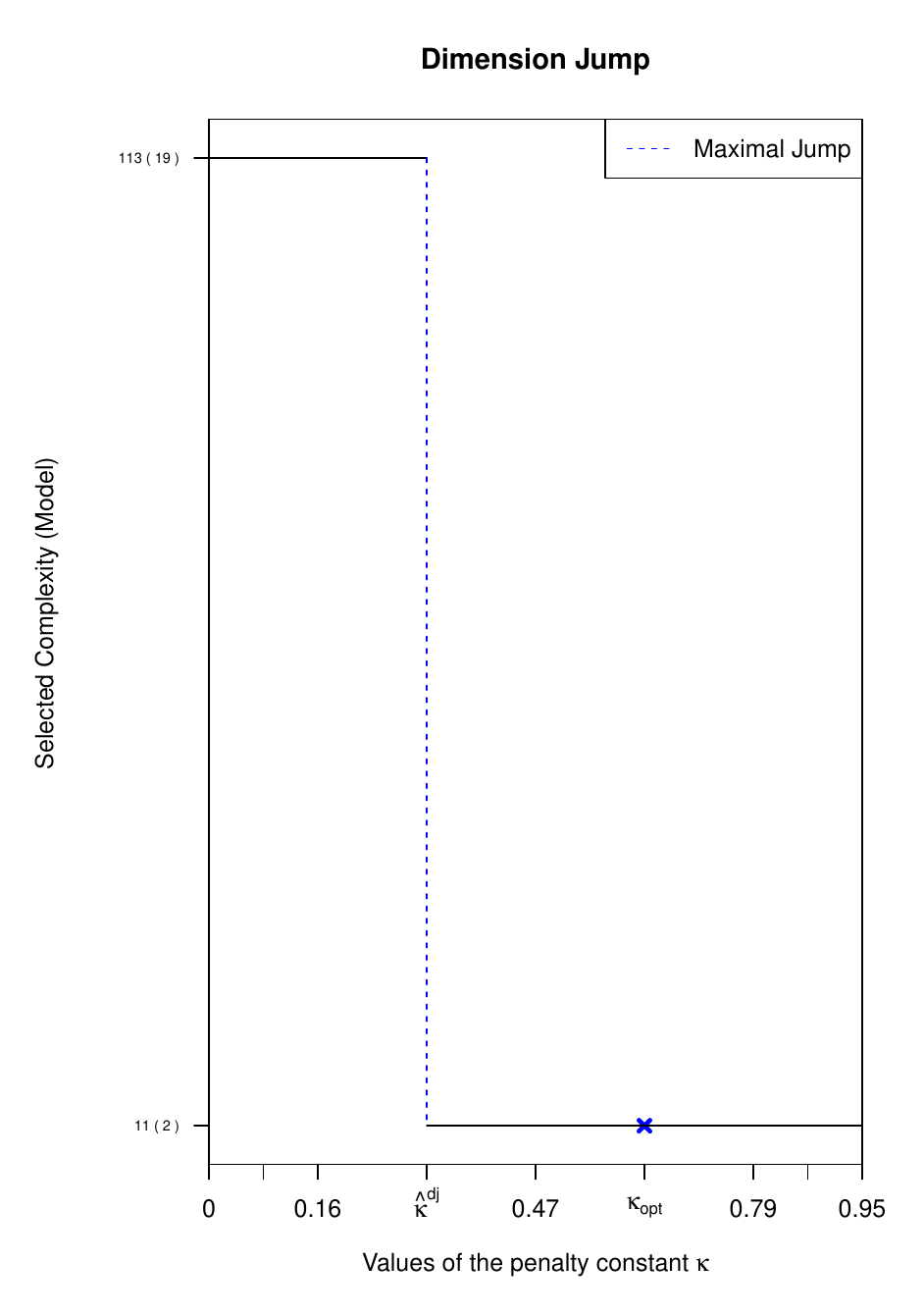}\\
		(a) Example WS.\\
		\includegraphics[height=0.55\textheight]{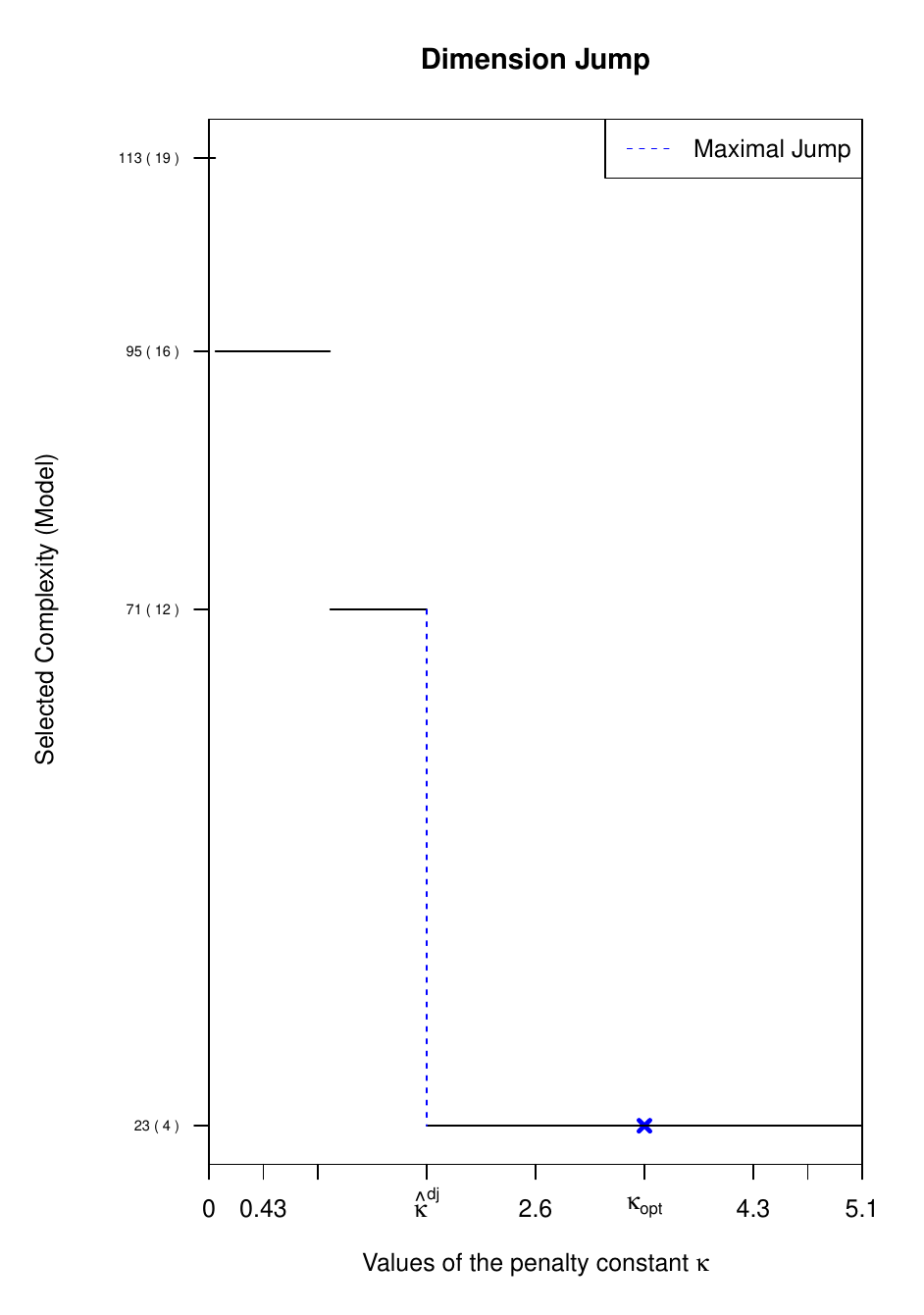} \\
		(b) Example MS.
		%
	\end{tabular}
	\caption{Plot of the selected model dimension using the jump criterion with 10000 data points.}
	\label{fig_Slope_Heuristic_10000_Djump_Slope_CNLL}
\end{figure}

\begin{figure}
	\centering
	\begin{subfigure}{.5\textheight}
		\centering
		\includegraphics[height=0.5\textheight]{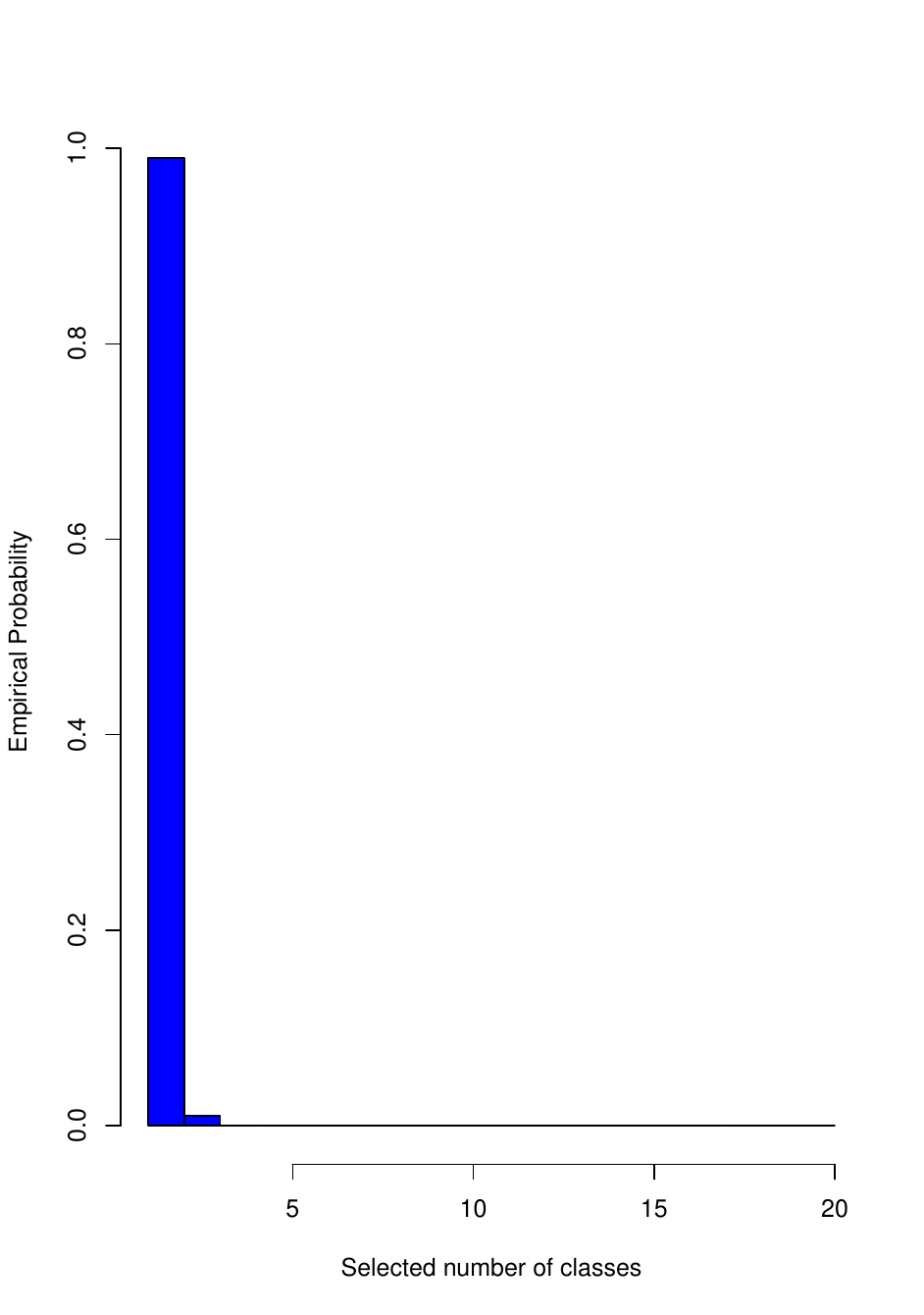}
		\caption{2000 data points.}
		\label{subfig_WS_2000_Djump_CNLL_NS}
	\end{subfigure}
	%
	\begin{subfigure}{.5\textheight}
		\centering
		\includegraphics[height=0.5\textheight]{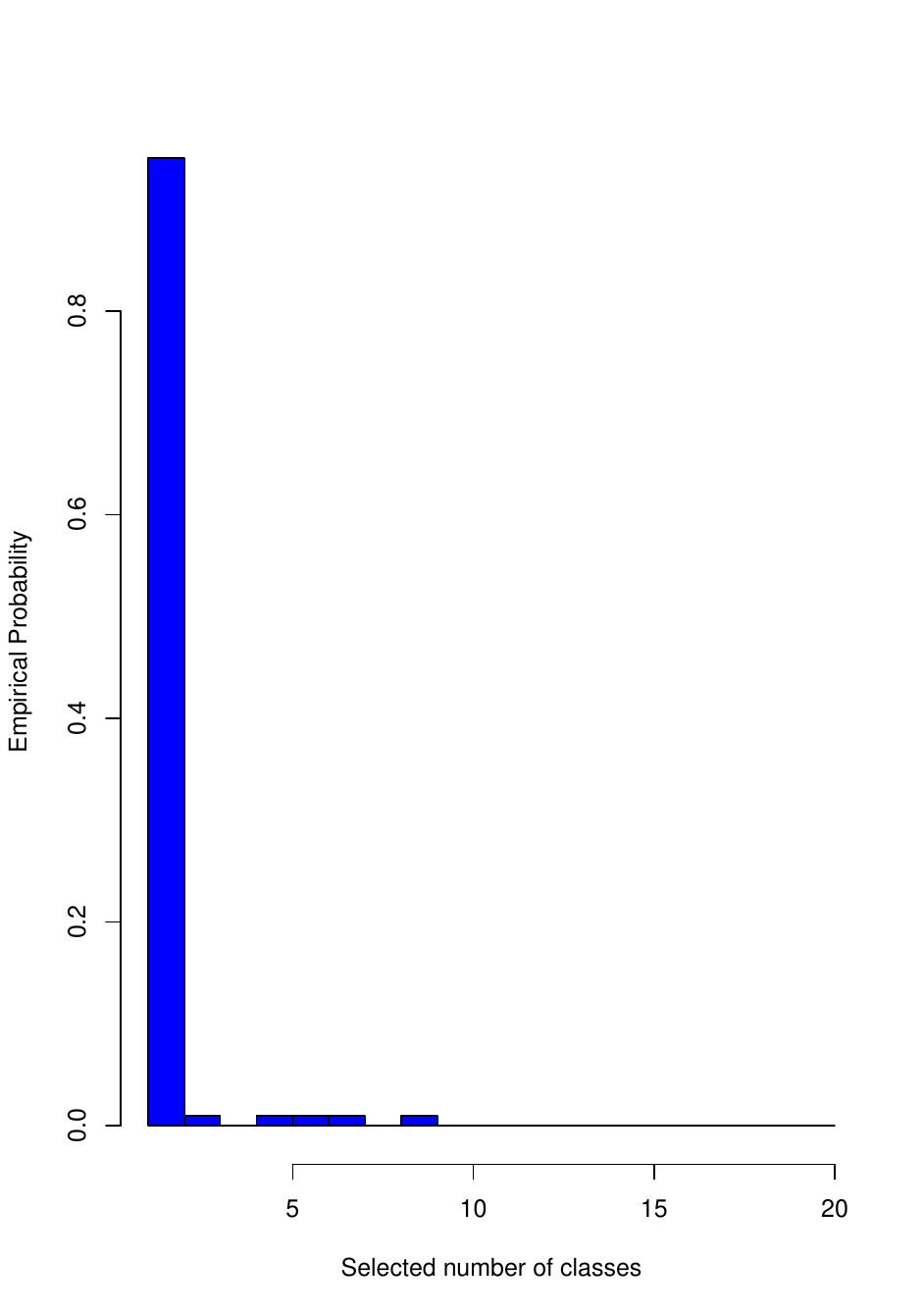}
		\caption{10000 data points.}
		\label{subfig_WS_10000_Djump_CNLL_NS}
	\end{subfigure}
	\caption{Comparison histograms of selected $K$ in WS case using jump criterion over 100 trials.}
	\label{fig_histogramK_Djump_CNLL_WS}
\end{figure}

\begin{figure}
	\centering	\begin{subfigure}{.5\textheight}
		\centering
		\includegraphics[height=0.5\textheight]{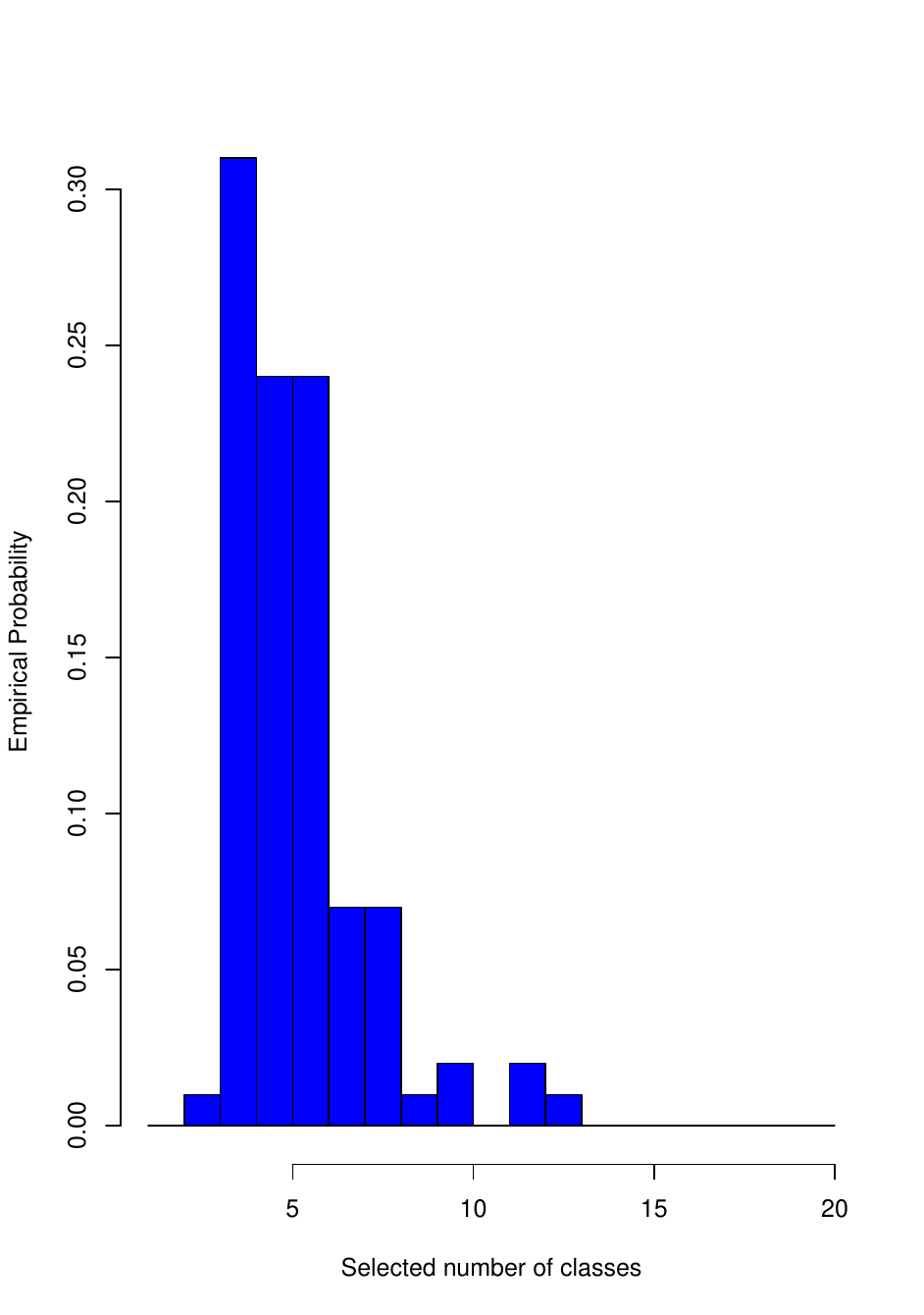}
		\caption{2000 data points.}
		\label{subfig_MS_2000_Djump_CNLL_NS}
	\end{subfigure}
	%
	\begin{subfigure}{.5\textheight}
		\centering
		\includegraphics[height=0.5\textheight]{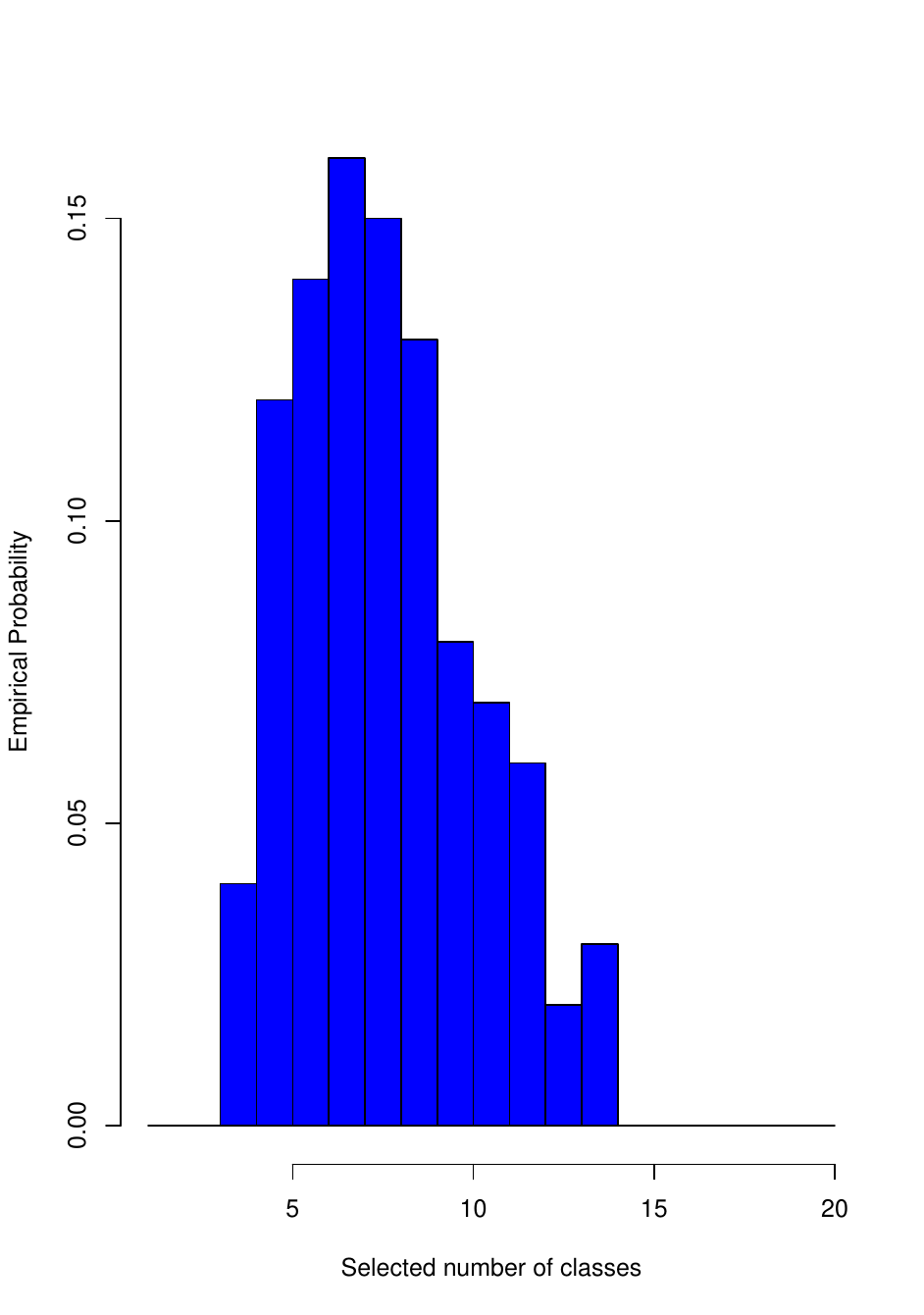}
		\caption{10000 data points.}
		\label{subfig_MS_10000_Djump_CNLL_NS}
	\end{subfigure}%
	\caption{Comparison histograms of selected $K$ in MS case using jump criterion over 100 trials.}
	\label{fig_histogramK_Djump_CNLL_MS}
\end{figure}

\begin{figure}
	\centering
	\begin{subfigure}{.5\textheight}
		\centering
		\includegraphics[height=0.5\textheight]{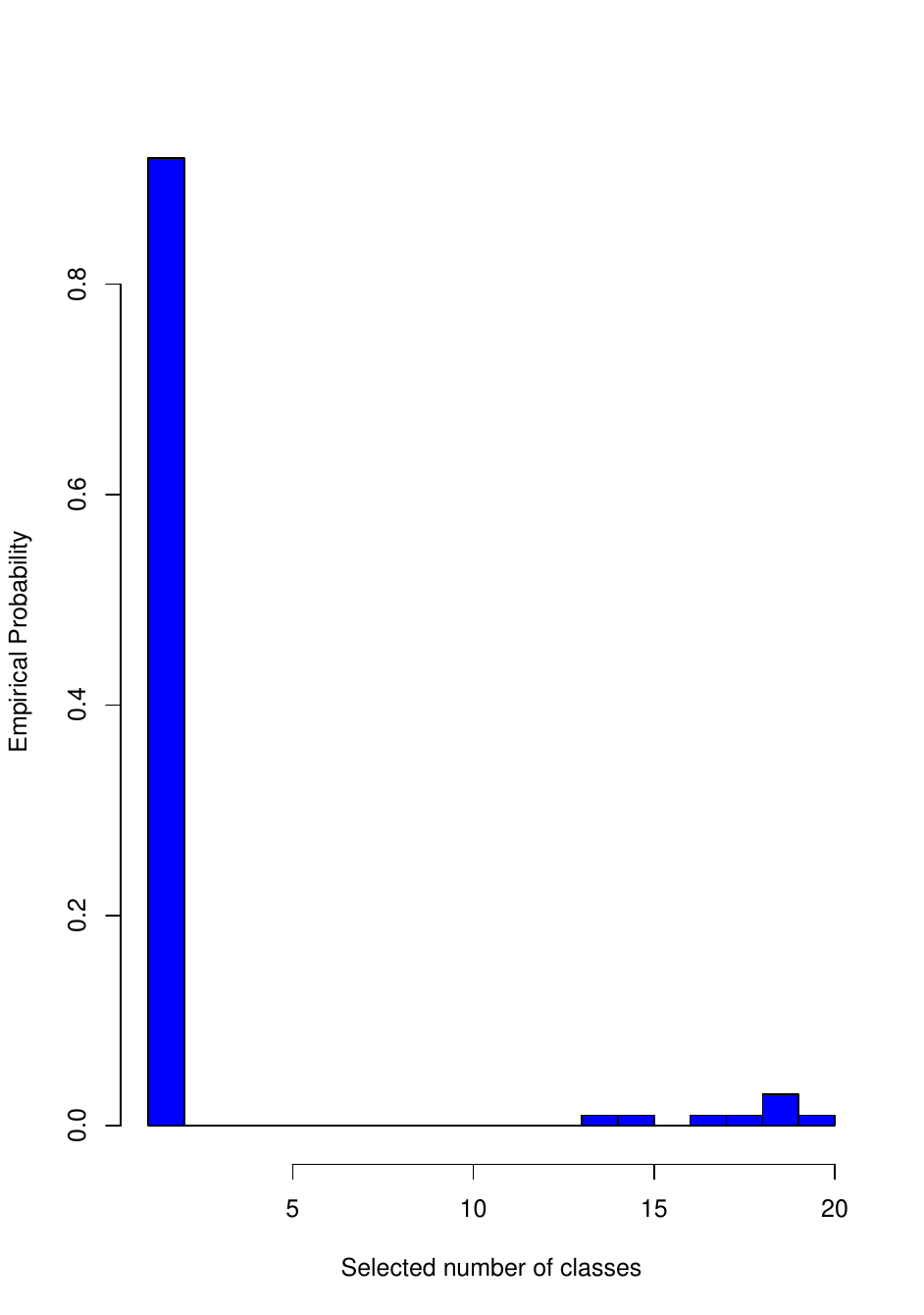}
		\caption{2000 data points.}
		\label{subfig_WS_2000_DDSE_CNLL_NS}
	\end{subfigure}
	%
	\begin{subfigure}{.5\textheight}
		\centering
		\includegraphics[height=0.5\textheight]{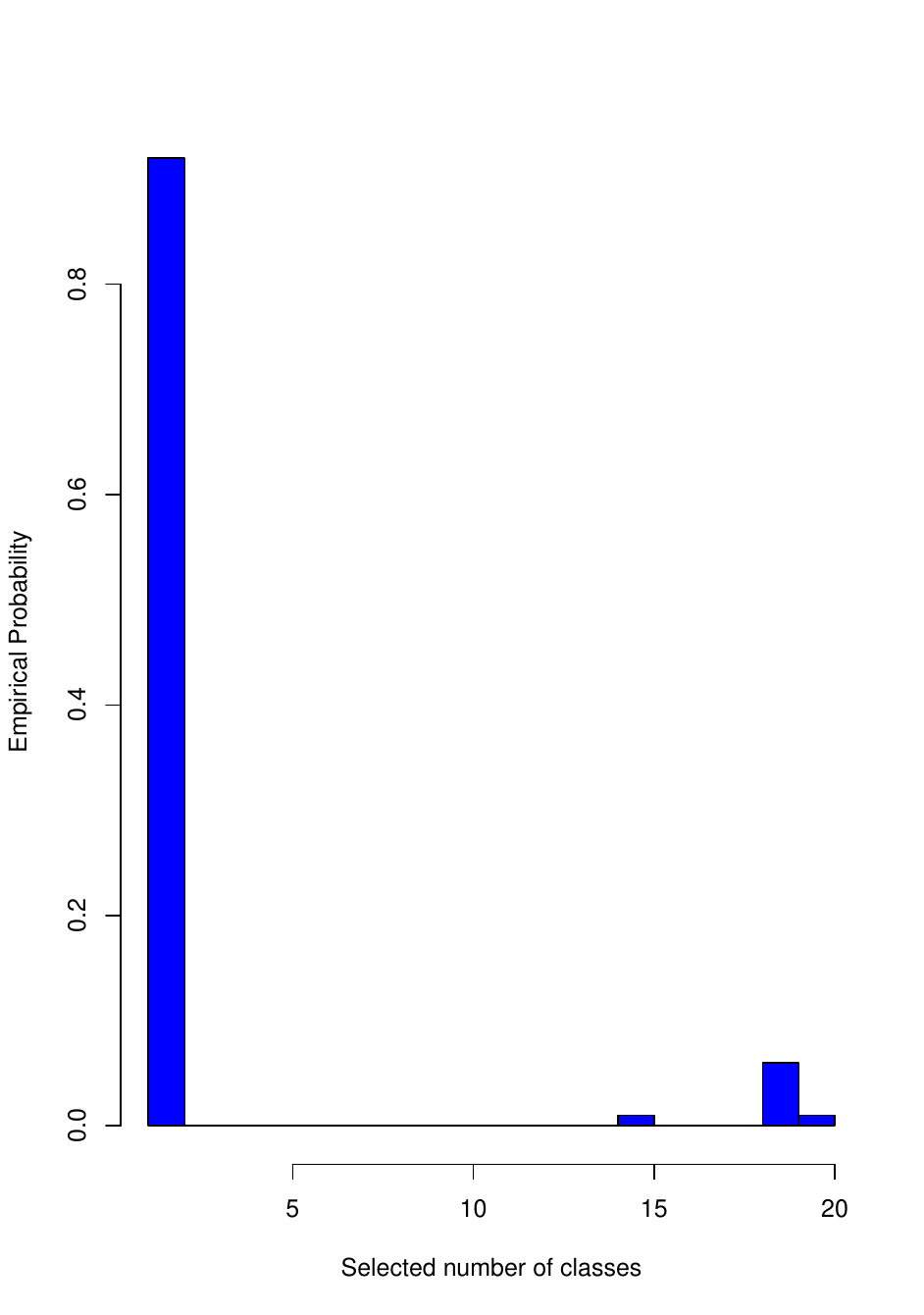}
		\caption{10000 data points.}
		\label{subfig_WS_10000_DDSE_CNLL_NS}
	\end{subfigure}
	\caption{Comparison histograms of selected $K$ in WS case using slope criterion over 100 trials.}
	\label{fig_histogramK_DDSE_CNLL_WS}
\end{figure}

\begin{figure}
	\centering	\begin{subfigure}{.5\textheight}
		\centering
		\includegraphics[height=0.5\textheight]{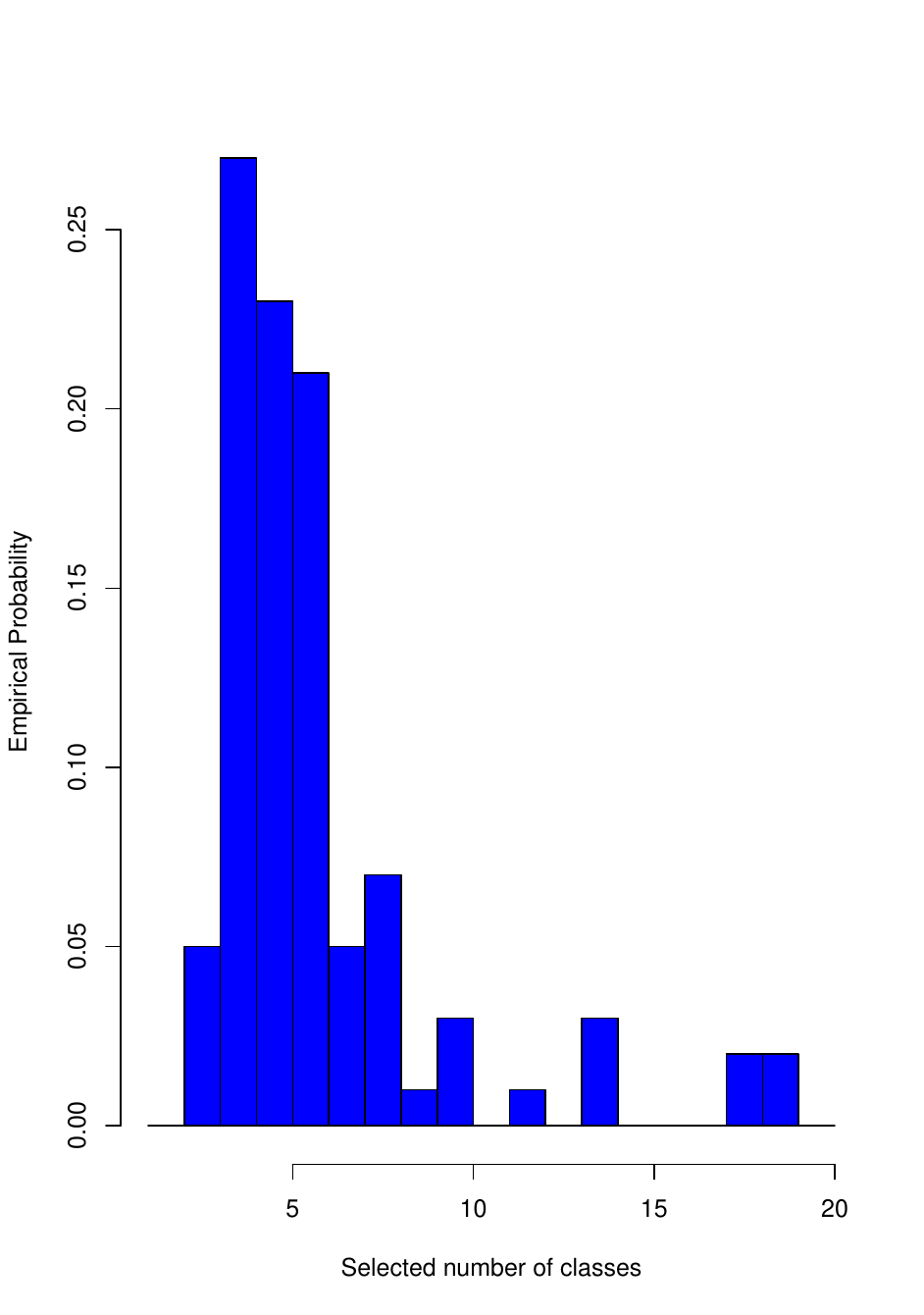}
		\caption{2000 data points}
		\label{subfig_MS_2000_DDSE_CNLL_NS}
	\end{subfigure}
	\hspace{0.5cm}
	\begin{subfigure}{.5\textheight}
		\centering
		\includegraphics[height=0.5\textheight]{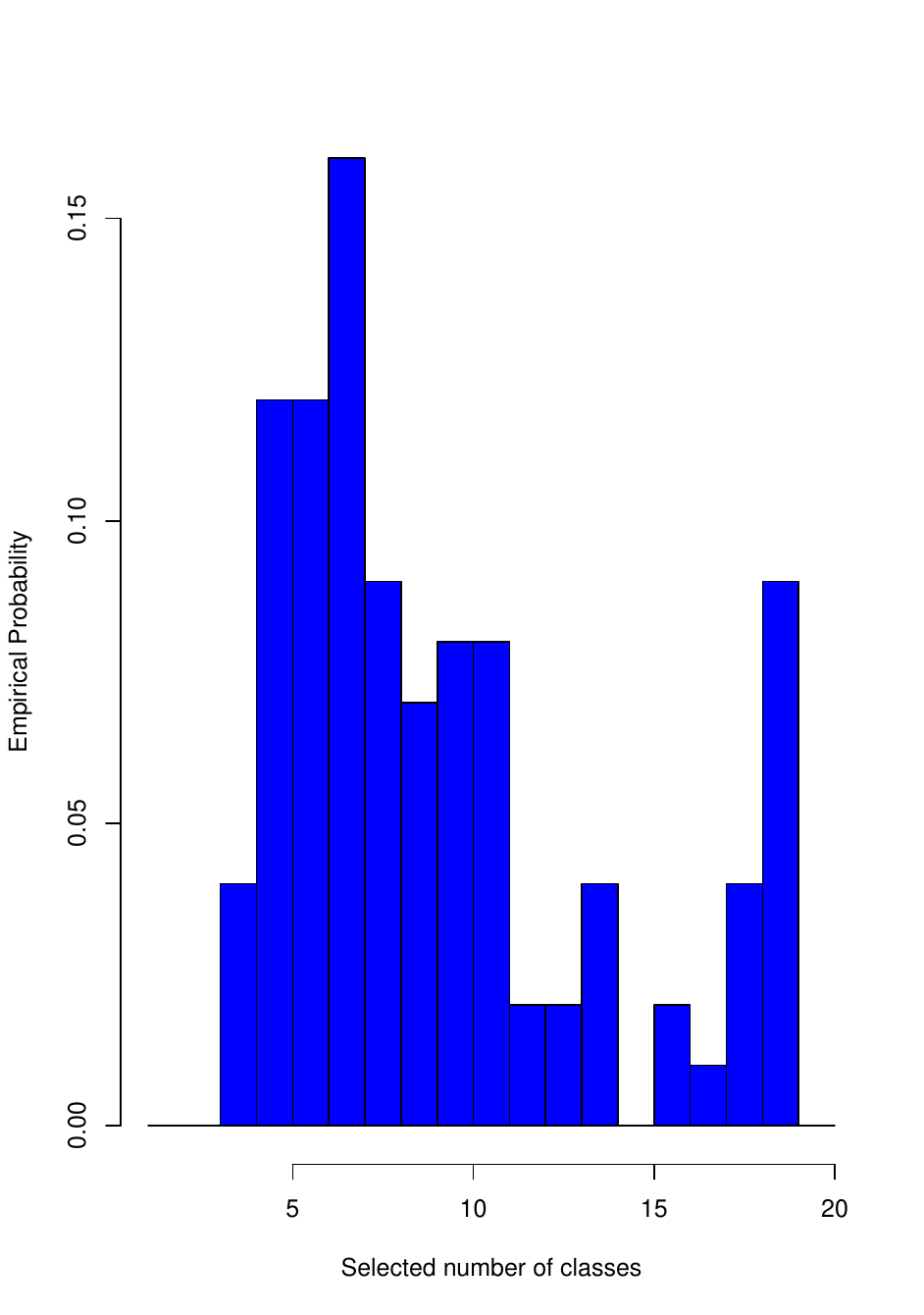}
		\caption{10000 data points.}
		\label{subfig_MS_10000_DDSE_CNLL_NS}
	\end{subfigure}%
	\caption{Comparison histograms of selected $K$ in MS case using slope criterion over 100 trials.}
	\label{fig_histogramK_DDSE_CNLL_MS}
\end{figure}

Firstly, by using the \href{https://cran.rstudio.com/web/packages/capushe/index.html}{capushe} (CAlibrating Penalities Using Slope HEuristics) package in {\bf R}  \citep{baudry2012slope}, we can select the penalty coefficient, along with the number of mixture components $K$. This heuristic comprises two possible criteria: the slope criterion and the jump criterion.  The first criterion consists of computing the \emph{asymptotic} slope of the log-likelihood (\cref{fig_Slope_Heuristic_2000_DDSE_Slope_CNLL,fig_Slope_Heuristic_10000_DDSE_Slope_CNLL}), drawn according to the model dimension, and then penalizing the log-likelihood by twice the slope times the model dimension. Regarding the second criterion, one aims to represent the dimension of the selected model according to $\kappa$ (\cref{fig_Slope_Heuristic_2000_Djump_Slope_CNLL,fig_Slope_Heuristic_10000_Djump_Slope_CNLL}), and find $\widehat{\kappa}$, such that if $\kappa < \widehat{\kappa}$, then the
dimension of the selected model is large, and of reasonable size, otherwise. The slope
heuristic prescribes then the use of $\kappa = 2 \widehat{\kappa}$. In our simulated data sets, \cref{fig_histogramK_Djump_CNLL_WS,fig_histogramK_Djump_CNLL_MS,fig_histogramK_DDSE_CNLL_WS,fig_histogramK_DDSE_CNLL_MS} show that the jump criterion appears to work better. The slope criterion sometimes chooses highly complex models in the WS case, with the problem exacerbated in the MS case.



Next, a close inspection shows that the bias-variance trade-off differs between the two examples. We run our experiment over $100$ trials with $K \in \cK = [20]$, using both the jump and slope criteria. The first remark is that the best choice of $K=2$ appears to be selected with very high probability, even for large samples ($n=10000$) in the WS case. This can be observed in \cref{subfig_WS_2000_Djump_CNLL_NS,subfig_WS_10000_Djump_CNLL_NS,subfig_WS_2000_DDSE_CNLL_NS,subfig_WS_10000_DDSE_CNLL_NS}. In the MS case, the best choice for $K$ should balance between the model approximation error term and the variance one, which is observed in \cref{subfig_MS_2000_Djump_CNLL_NS,subfig_MS_10000_Djump_CNLL_NS,subfig_MS_2000_DDSE_CNLL_NS,subfig_MS_10000_DDSE_CNLL_NS}. Here, the larger the number of samples $n$, the larger the value of $K$ that is selected as optimal.

From hereon in, we only focus on the jump criterion, due to its stability regarding the choice of $K$. We wish to measure the performances of our chosen GLoME models in term of the tensorized Kullback--Leibler divergence, $\tkl$, which can not be calculated exactly in the case of Gaussian mixtures. Therefore, we evaluate the divergence using a Monte Carlo simulation, since we know the true density. We note that the variability of this randomized approximation has been demonstrated to be negligible in practice, which is also supported in the numerical experiments by \cite{montuelle2014mixture}. More precisely, we compute the Monte Carlo approximation for the tensorized Kullback--Leibler divergence as follows. First, note that the Monte Carlo approximation for tensorized Kullback--Leibler divergence between the true model, $s_0^*$, and the selected model $\widehat{s}_{\widehat{\bfm}}^*$ can be approximated as
\begin{align*}
\tkl\left(s_0^*,\widehat{s}_{\widehat{\bfm}}^*\right) &= \E{\Xv}{\frac{1}{n} \sum_{i=1}^n \kl\left(s_0^*\left(\cdot\mid  \xv_i\right),\widehat{s}_{\widehat{\bfm}}^*\left(\cdot\mid  \xv_i\right)\right)}\nn\\
& \approx \frac{1}{n} \sum_{i=1}^n \kl\left(s_0^*\left(\cdot\mid  \xv_i\right),\widehat{s}_{\widehat{\bfm}}^*\left(\cdot\mid  \xv_i\right)\right)
\nn\\
& = \frac{1}{n} \sum_{i=1}^n \frac{1}{n_y} \sum_{j=1}^{n_y} \ln\left(\frac{s_0^*\left(\yv_{i,j}\mid  \xv_i\right)}{\widehat{s}_{\widehat{\bfm}}^*\left(\yv_{i,j}\mid  \xv_i\right)}\right),
\end{align*}
where the data $\xv_i, i \in [n]$, and $\left(\yv_{i,j}\right)_{j \in \left[n_y\right]}$ are drawn from $s_0^*\left(\cdot\mid  \xv_i\right)$. Then, $\E{}{\tkl\left(s_0^*,\widehat{s}_{\widehat{\bfm}}^*\right)}$ is approximated again by averaging over $N_t = 100$ Monte Carlo trials. Therefore, the simulated data used for approximation can be written as $\left(\xv_i,\yv_{i,j}\right)_{t}$ with $i\in [n], j \in \left[n_y\right], t\in \left[N_t\right]$.

\begin{figure}
	\centering		
	\begin{subfigure}{.65\textheight}
		\centering
		\includegraphics[width=\textwidth]{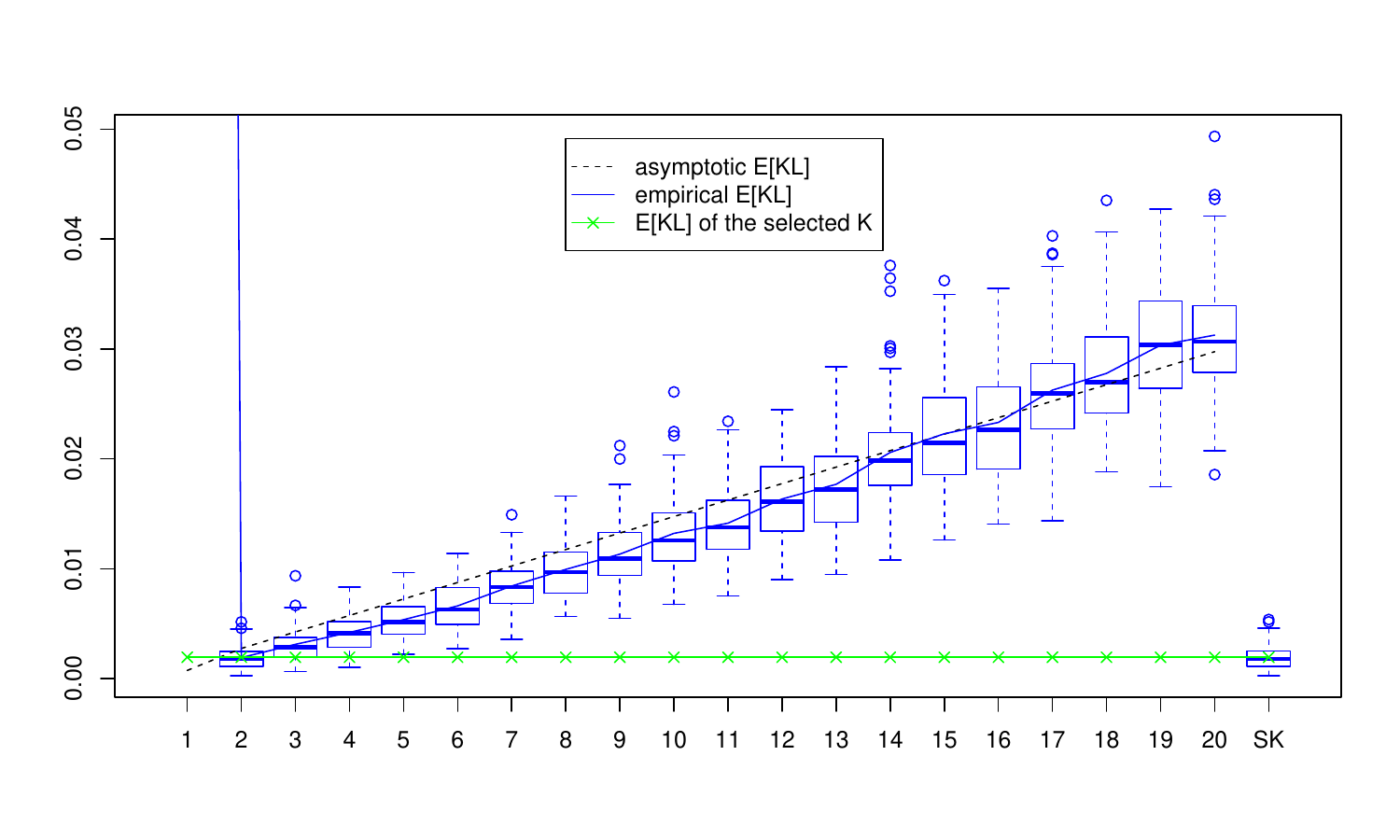}
		\caption{2000 data points.}
		\label{Boxplot_KL_WS_2000_Djump_CNLL}
	\end{subfigure}%
	
	\begin{subfigure}{.65\textheight}
		\centering
		\includegraphics[width=\textwidth]{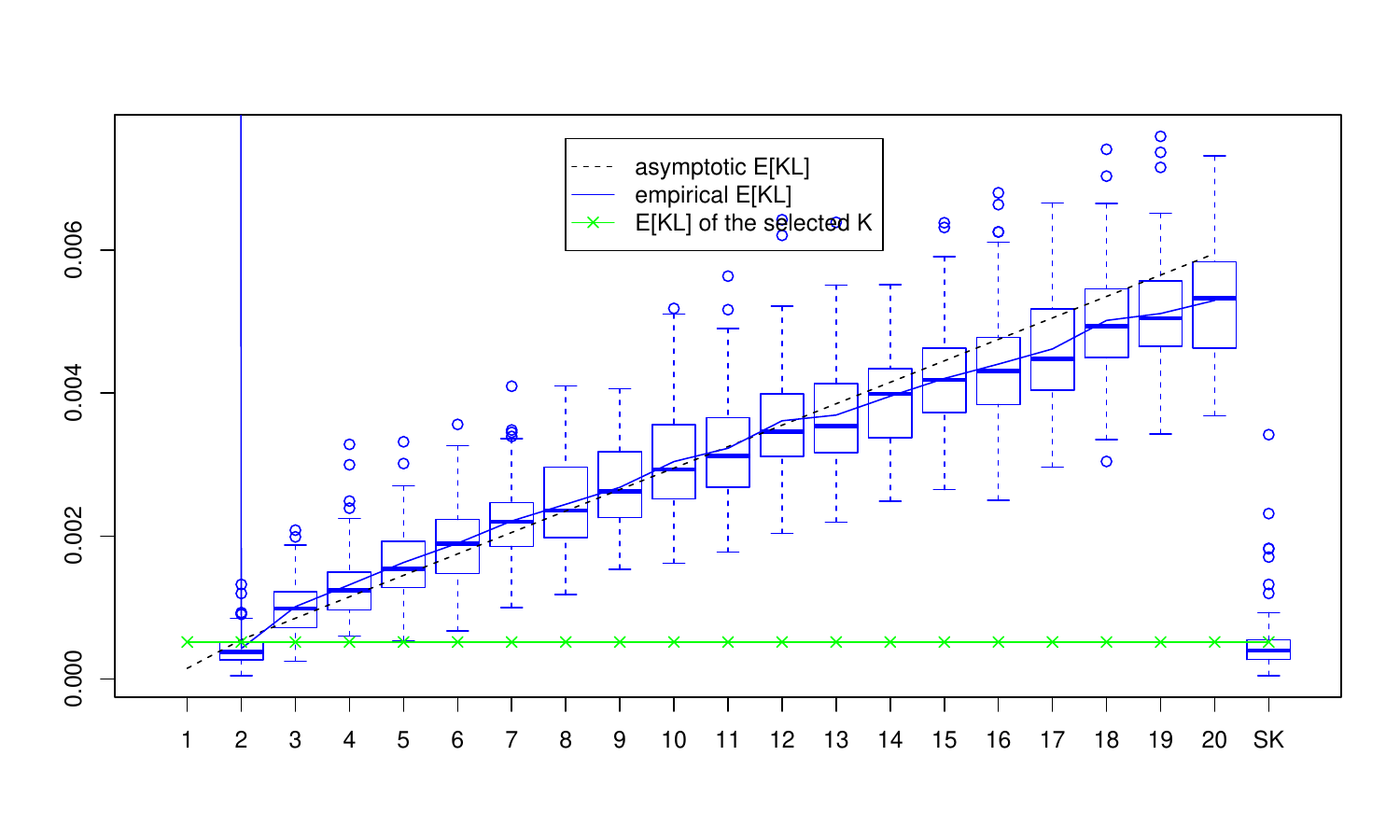}
		\caption{10000 data points.}
		\label{Boxplot_KL_WS_10000_Djump_CNLL}
	\end{subfigure}
	
	\caption{Box-plot of the tensorized Kullback--Leibler divergence according to the number of mixture components in WS case using the jump criterion over 100 trials. The tensorized Kullback--Leibler divergence of the penalized estimator $\widehat{s}_{\widehat{K}}$ is shown in the right-most box-plot of each graph.}
	\label{Boxplot_KL_Djump_CNLL_WS}
\end{figure}

\begin{figure}
	\centering		
		\begin{subfigure}{.65\textheight}
			\centering
			\includegraphics[width=\textwidth]{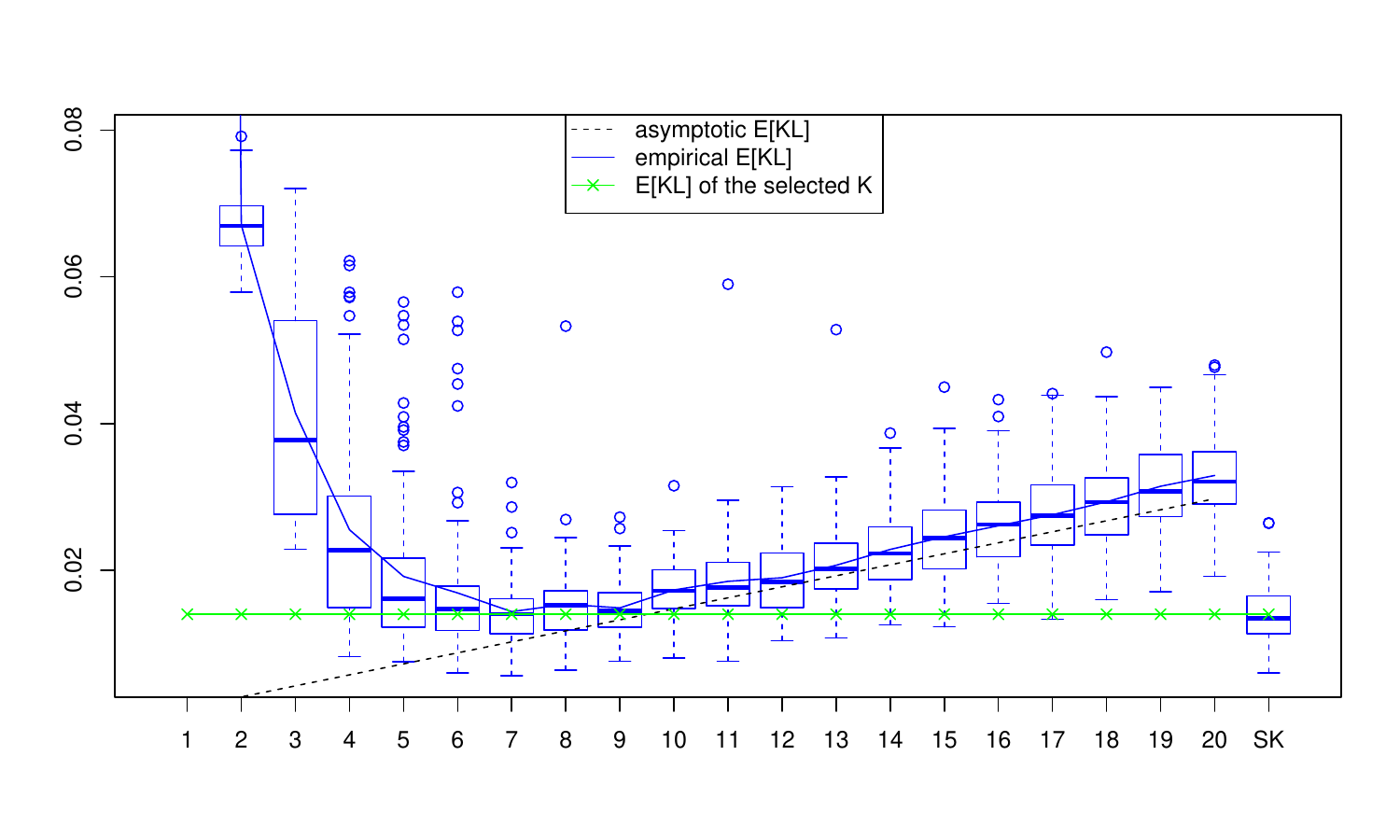}
				\caption{2000 data points.}
				\label{Boxplot_KL_MS_2000_Djump_CNLL}
			\end{subfigure}
		%
		
		\begin{subfigure}{.65\textheight}
			\centering
			\includegraphics[width=\textwidth]{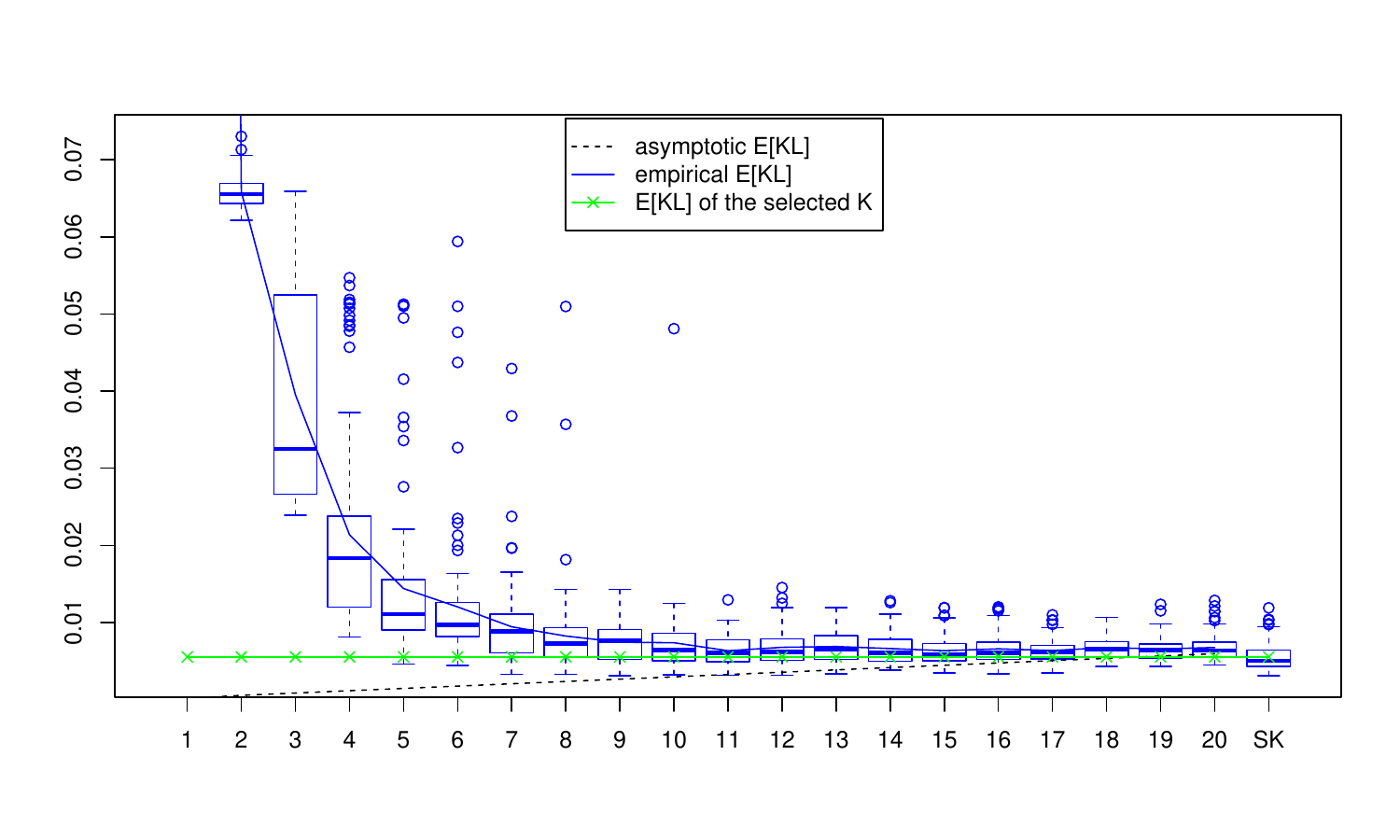}
				\caption{10000 data points.}
				\label{Boxplot_KL_MS_10000_Djump_CNLL}
			\end{subfigure}
	\caption{Box-plot of the tensorized Kullback--Leibler divergence according to the number of mixture components in MS case using the jump criterion over 100 trials. The tensorized Kullback--Leibler divergence of the penalized estimator $\widehat{s}_{\widehat{K}}$ is shown in the right-most box-plot of each graph.}
	\label{Boxplot_KL_Djump_CNLL_MS}
\end{figure}

Based on the approximation, \cref{Boxplot_KL_Djump_CNLL_WS,Boxplot_KL_Djump_CNLL_MS} show the box plots and the mean of the tensorized Kullback--Leibler divergence over $100$ trials, based on the jump criterion. Our boxplots confirm that the mean tensorized Kullback--Leibler divergence between $\widehat{s}^*_K$ and $s_0^*$, over $K \in \left\{1,\ldots,20\right\}$ number of mixture components, is always larger than the mean of tensorized Kullback--Leibler divergence between the penalized estimator $\widehat{s}_{\widehat{K}}^*$ and $s_0^*$, which is consistent with \cref{thm_Oracle_Inequality_GLoME}. In particular, if the true model belongs to our nested collection, the mean tensorized Kullback--Leibler divergence seems to behave like $\frac{\dim\left(S_\bfm ^*\right)}{2n}$ (shown by a dotted line), which can be explained by the AIC heuristic. More precisely, we firstly assume that \begin{align*}
S_\bfm ^*=\left\{\cX \times \cY \ni (\xv,\yv) \mapsto s_\bfm^*:=s_{\psib_\bfm^*}(\yv\mid \xv):\psib_\bfm^* \in \Psib_\bfm^* \subset \R^{\dim\left(S_\bfm ^*\right)}\right\}
\end{align*}
is identifiable and make some strong regularity assumptions on $\psib_\bfm^* \mapsto s_{\psib_m}^*$. Further, we assume the existence of $\dim\left(S_\bfm ^*\right)\times \dim\left(S_\bfm ^*\right)$ matrices $\bfA\left(\psib_\bfm^*\right)$ and $\Bb\left(\psib_\bfm^*\right)$, which are defined as follows:
\begin{align*}
\left[\bfA\left(\psib_\bfm^*\right)\right]_{k,l}& = \E{}{\frac{1}{n}\sum_{i=1}^n \int \frac{-\partial^2 \ln s_{\psib_\bfm^*}}{\partial \psib_{\bfm,k}^*\partial \psib_{\bfm,l}^*}\left(\yv\mid \xv_i\right) s_0^* \left(\yv\mid \xv_i\right)d \yv},\nn\\
\left[\Bb\left(\psib_\bfm^*\right)\right]_{k,l}& = \E{}{\frac{1}{n}\sum_{i=1}^n \int \frac{\partial \ln s_{\psib_\bfm^*}}{\partial \psib_{\bfm,k}^*}\left(\yv\mid \xv_i\right)\frac{\partial \ln s_{\psib_\bfm^*}}{\partial \psib_{\bfm,l}^*}\left(\yv\mid \xv_i\right) s_0^* \left(\yv\mid \xv_i\right)d \yv}.
\end{align*}
Then, the results from \cite{white1982maximum} and \cite{cohen2011conditional} imply that $\E{}{\tkl \left(s_0^*,\widehat{s}_\bfm^*\right)}$ is asymptotically equivalent to 
$$\tkl \left(s_0^*,s_{\psib^{**}_\bfm}\right)+\frac{1}{2n}\tr\left(\Bb\left(\psib^{**}_\bfm\right)\bfA\left(\psib^{**}_\bfm\right)^{-1}\right),$$
where we defined $\psib^{**}_\bfm = \argmin_{s_{\psib_\bfm^*} \in S_\bfm ^*} \tkl \left(s_0^*,s_{\psib^{*}_m}\right)$.

In particular, $\E{}{\tkl \left(s_0^*,\widehat{s}_\bfm^*\right)}$ is asymptotically equivalent to $\frac{1}{2n}\dim\left(S_\bfm ^*\right)$, whenever $s_0^*$ belongs to the model collection $S_\bfm ^*$. Furthermore, even though there is no theoretical guarantee, the slope of the mean error in the misspecified case seems also to grow at the same rate as $\frac{\dim\left(S_\bfm ^*\right)}{2n}$, for large enough number of mixture components ($K \ge 6$ in the WS case and $K \ge 9$ in the MS case). 


%

%


\begin{figure}
\centering
\begin{subfigure}{.4\textheight}
	\centering
	\includegraphics[width=0.85\textwidth]{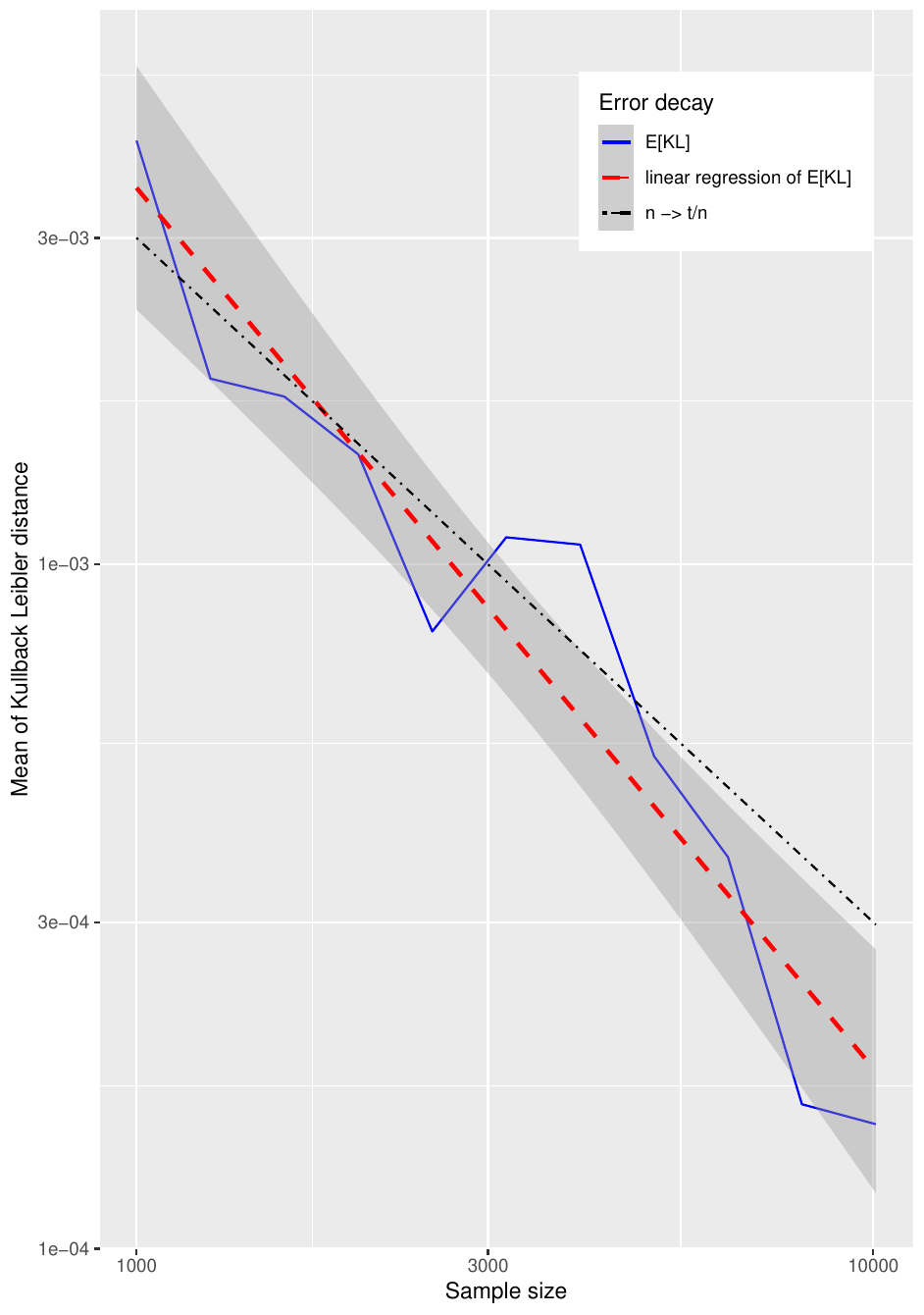}
	\caption{WS Example: The slope of the
		free regression line is $\approx -1.287 $ and $t=3$.}
	\label{KL_Sel_WS_Djump_CNLL_ErrorDecay}
\end{subfigure}
%
\begin{subfigure}{.4\textheight}
	\centering
	\includegraphics[width=0.85\textwidth]{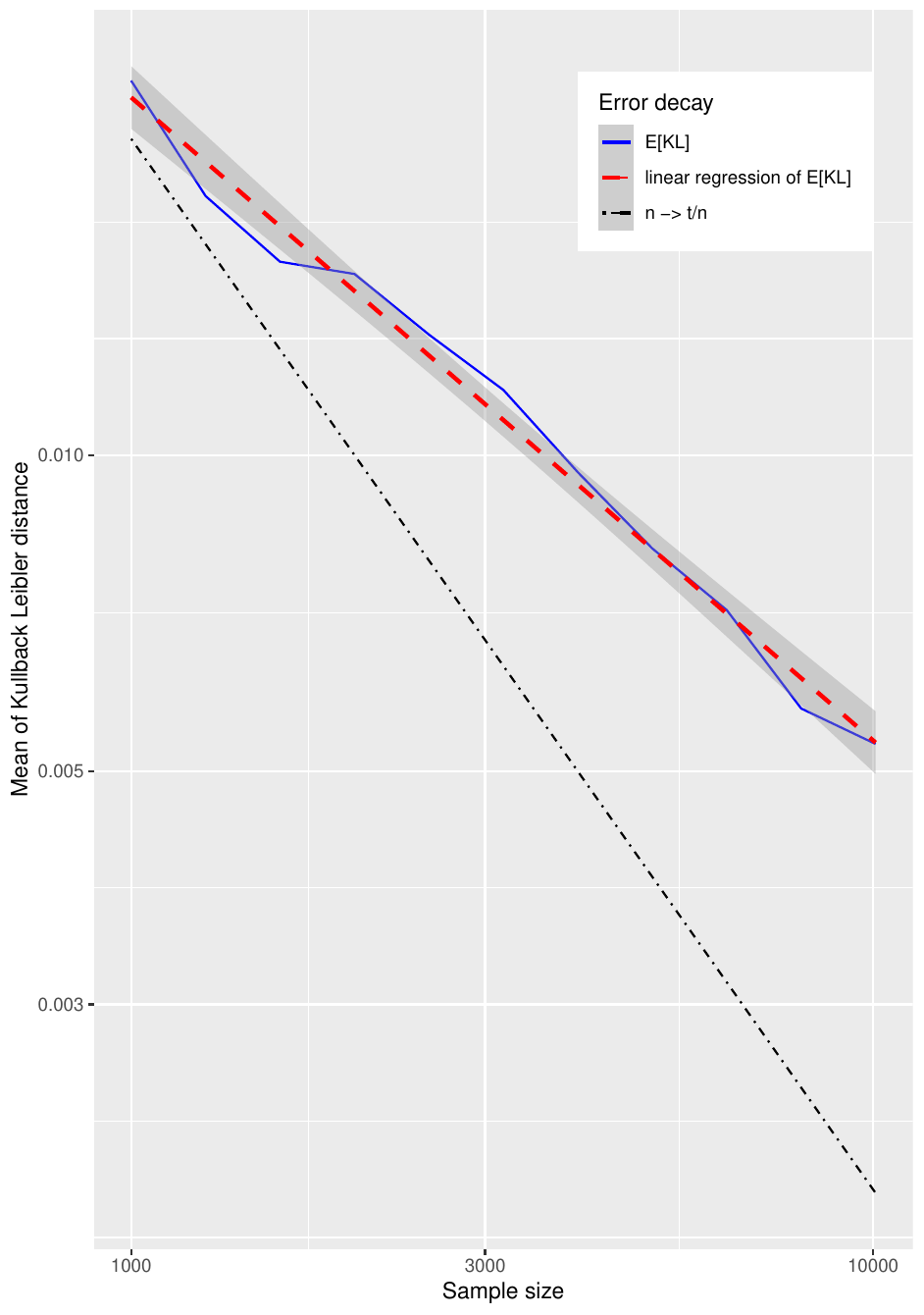}
	\caption{MS Example: The slope of the
		free regression line is $\approx -0.6120  $ and $t=20$.}
	\label{KL_Sel_MS_Djump_CNLL_ErrorDecay}
\end{subfigure}%
\caption{Tensorized Kullback--Leibler divergence between the true and selected densities based on the jump criterion, represented in a log-log scale, using 30 trials. A free least-square regression with standard error and a regression with slope $-1$ were added to stress the two different behavior for each graph.}
\label{KL_Sel_Djump_JNLL_CNLL_ErrorDecay}
\end{figure}
\cref{KL_Sel_Djump_JNLL_CNLL_ErrorDecay} shows that the error decays when the sample size $n$ grows, when using the penalty based on the jump criterion. The first remark is that we observed the error decay is of order $t/n$, as predicted in by the theory, where $t$ is some constant, as expected in the well-specified case. The rate of convergence for the misspecified case seems to be slower.

\subsection{Ethanol data set}
\begin{figure}
	\centering
	\begin{subfigure}{.48\textwidth}
		\centering
		\includegraphics[scale = .35]{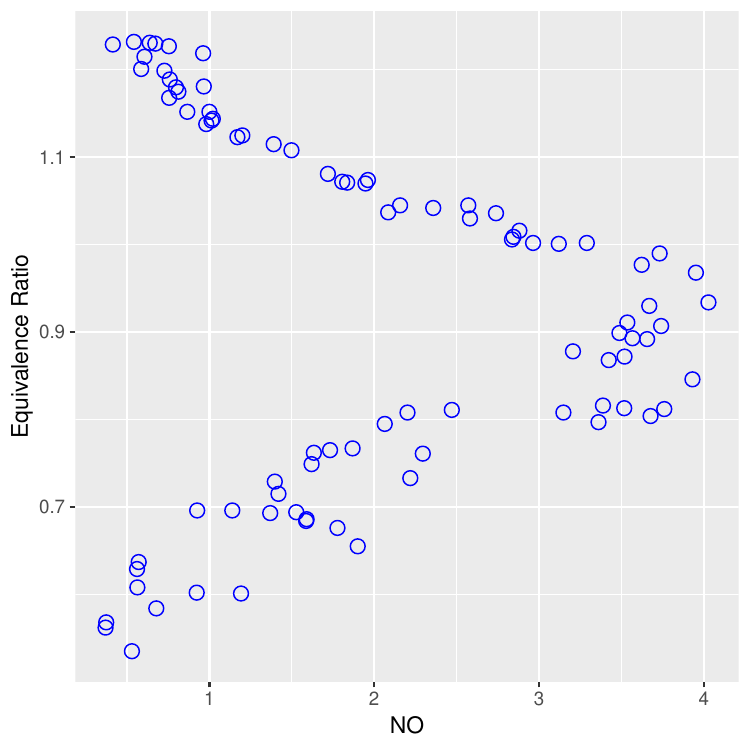}
		\caption{Raw Ethanol data set based on NO.}
		\label{subfig_Ethanol_Raw_NO_EquivRatio}
	\end{subfigure}
	\hspace{0.3cm}
	\begin{subfigure}{.48\textwidth}
		\centering
		\includegraphics[scale = .35]{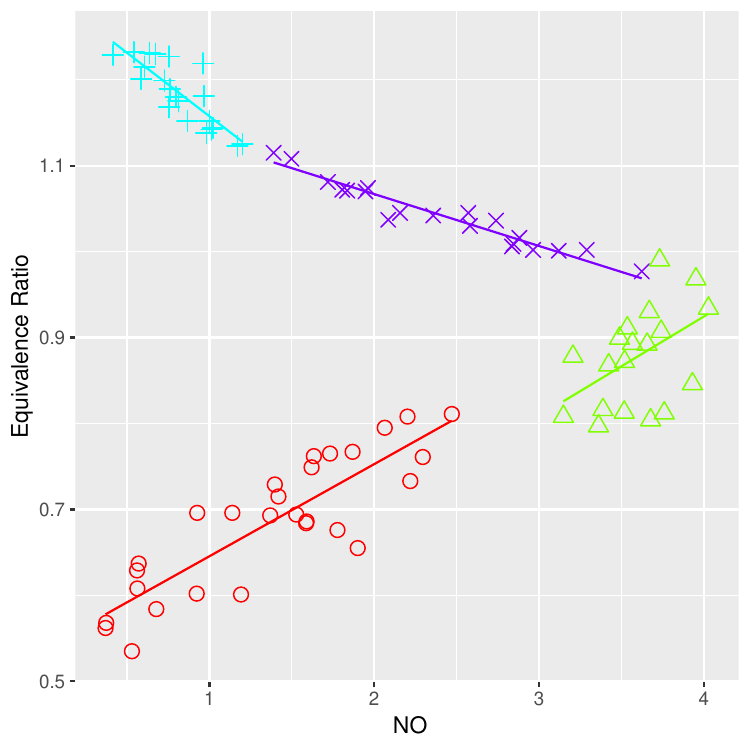}
		\caption{Clustering by GLoME based on NO.}
		\label{subfig_Ethanol_Clustering_NO_EquivRatio}
	\end{subfigure}%
	
	\begin{subfigure}{.48\textwidth}
		\centering
		\includegraphics[scale = .35]{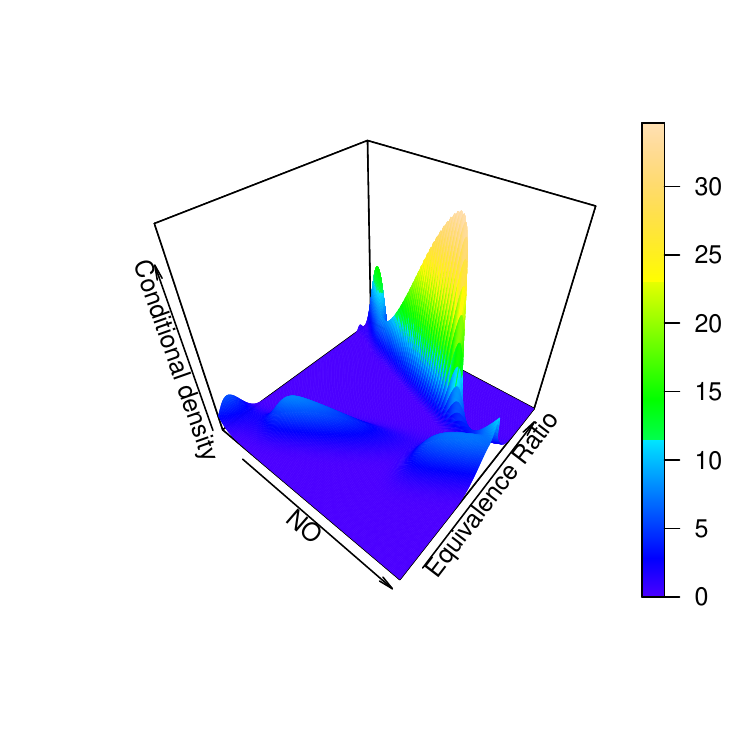}
		\caption{3D view of the resulting conditional density based on NO with the 4 clusters.}
		\label{subfig_Ethanol_Conditional_density_3D_NO_EquivRatio}
	\end{subfigure}%
	\hspace{0.3cm}
	\begin{subfigure}{.48\textwidth}
		\centering
		\includegraphics[scale = .35]{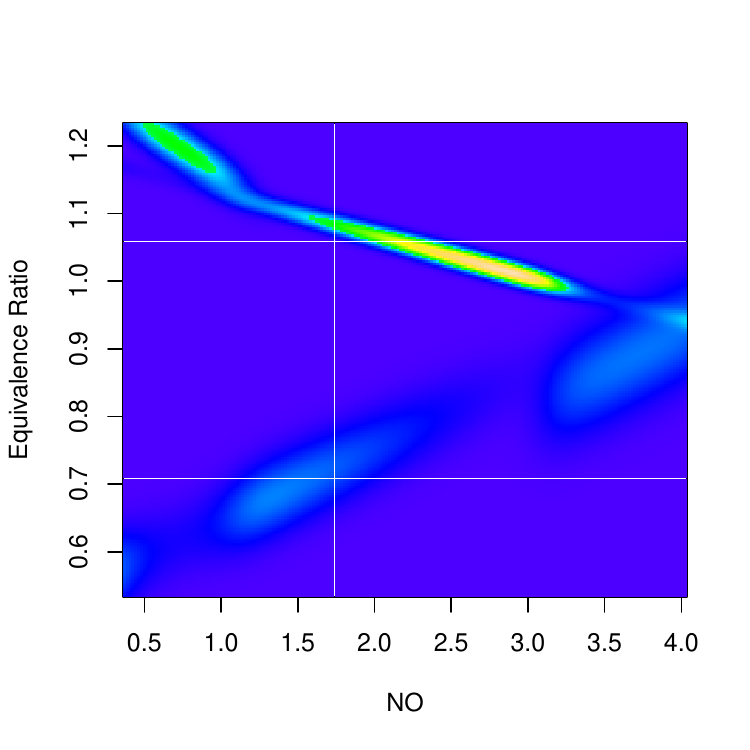}
		\caption{2D view of the same conditional density on NO.}
		\label{subfig_Ethanol_Conditional_density_2D_NO_EquivRatio}
	\end{subfigure}%
	
	%
	%
	
	\begin{subfigure}{.48\textwidth}
		\centering
		\includegraphics[scale = .35]{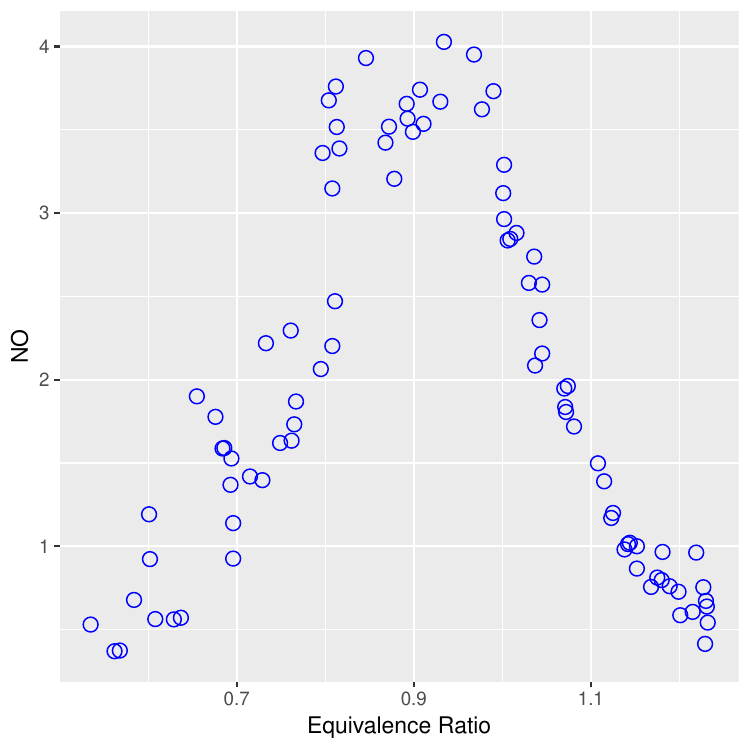}
		\caption{Raw Ethanol data set based on ER.}
		\label{subfig_Ethanol_Raw_EquivRatio_NO}
	\end{subfigure}%
	\hspace{0.3cm}
	\begin{subfigure}{.48\textwidth}
		\centering
		\includegraphics[scale = .35]{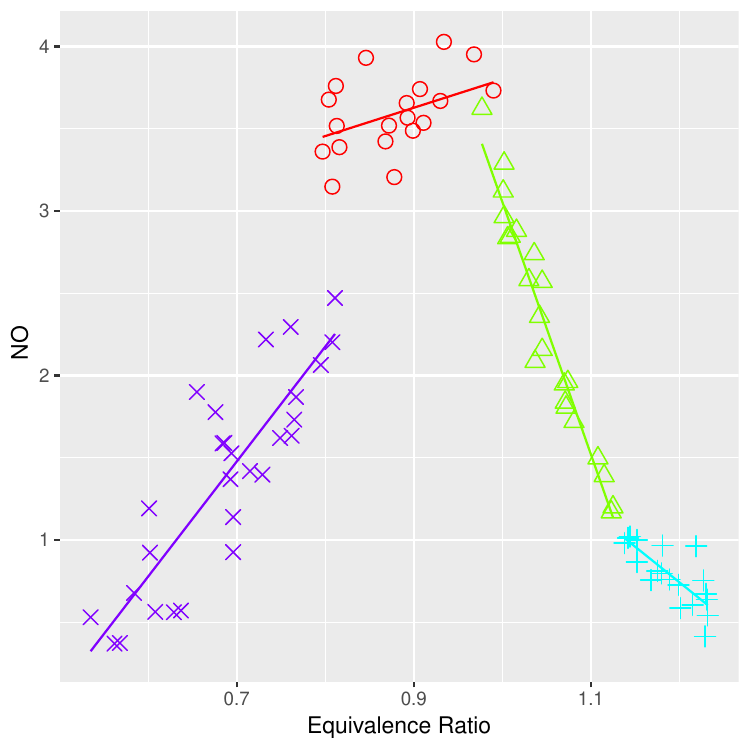}
		\caption{Clustering by GLoME based on ER.}
		\label{subfig_Ethanol_Clustering_EquivRatio_NO}
	\end{subfigure}%

	\begin{subfigure}{.48\textwidth}
		\centering
		\includegraphics[scale = .3]{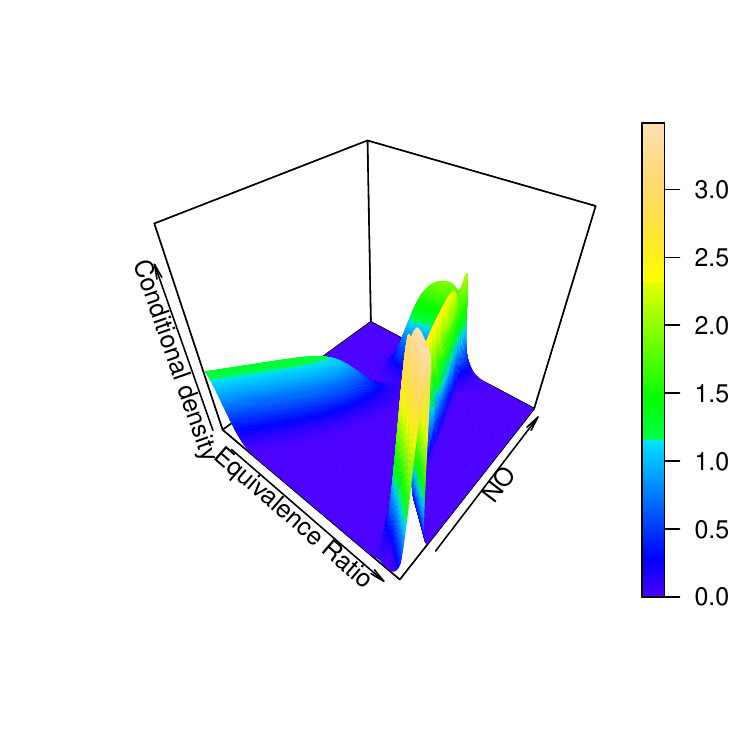}
		\caption{3D view of the estimated conditional density with the 4 clusters.}
		\label{subfig_Ethanol_Conditional_density_3D_EquivRatio_NO}
	\end{subfigure}%
	\hspace{0.3cm}
	\begin{subfigure}{.48\textwidth}
		\centering
		\includegraphics[scale = 0.3]{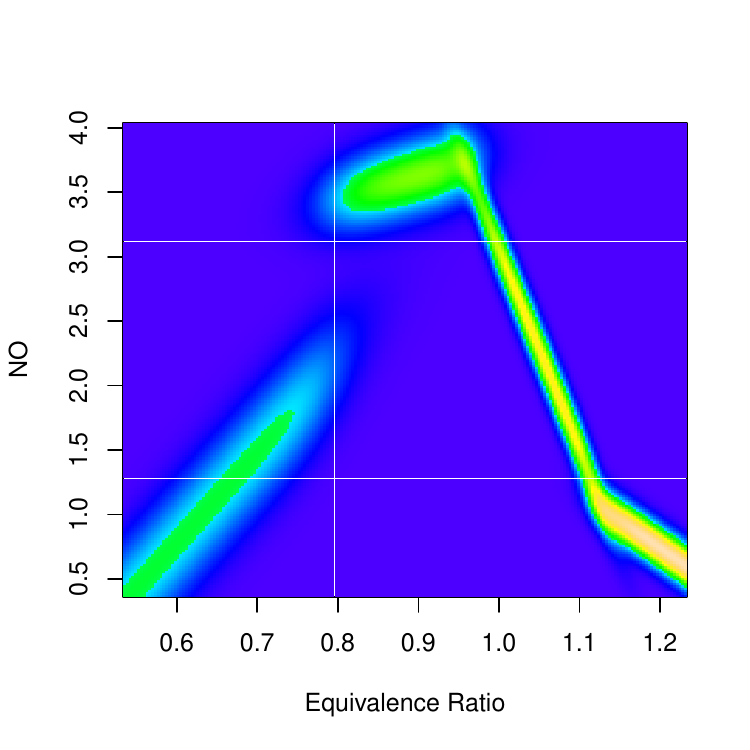}
		\caption{2D view of the same conditional density.}
		\label{subfig_Ethanol_Conditional_density_2D_EquivRatio_NO}
	\end{subfigure}%
	
	\caption{Estimated conditional density with 4 components based upon on the covariate NO or ER from the Ethanol data set.}
	\label{fig_Ethanol_Conditional_density_EquivRatio_NO}
\end{figure}

We now consider the use of GLoME models for performing clustering and regression tasks on a real data set. Following the numerical experiments from \cite{young2014mixtures} and \cite{montuelle2014mixture}, we demonstrate our model on the ethanol data set of \cite{brinkman1981ethanol}. The data comprises of 88 observations, which represent the relationship between the engine’s concentration of nitrogen oxide (NO) emissions and the equivalence ratio (ER), a measure of the air-ethanol mix, used as a spark-ignition engine fuel in a single-cylinder automobile test (\cref{subfig_Ethanol_Raw_NO_EquivRatio,subfig_Ethanol_Raw_EquivRatio_NO}).
Our goal is then to estimate the parameters of a GLoME model, as well as the number of mixture components.

\begin{figure}
	\centering
	\begin{subfigure}{.45\textheight}
		\centering
		\includegraphics[width=0.8\textwidth]{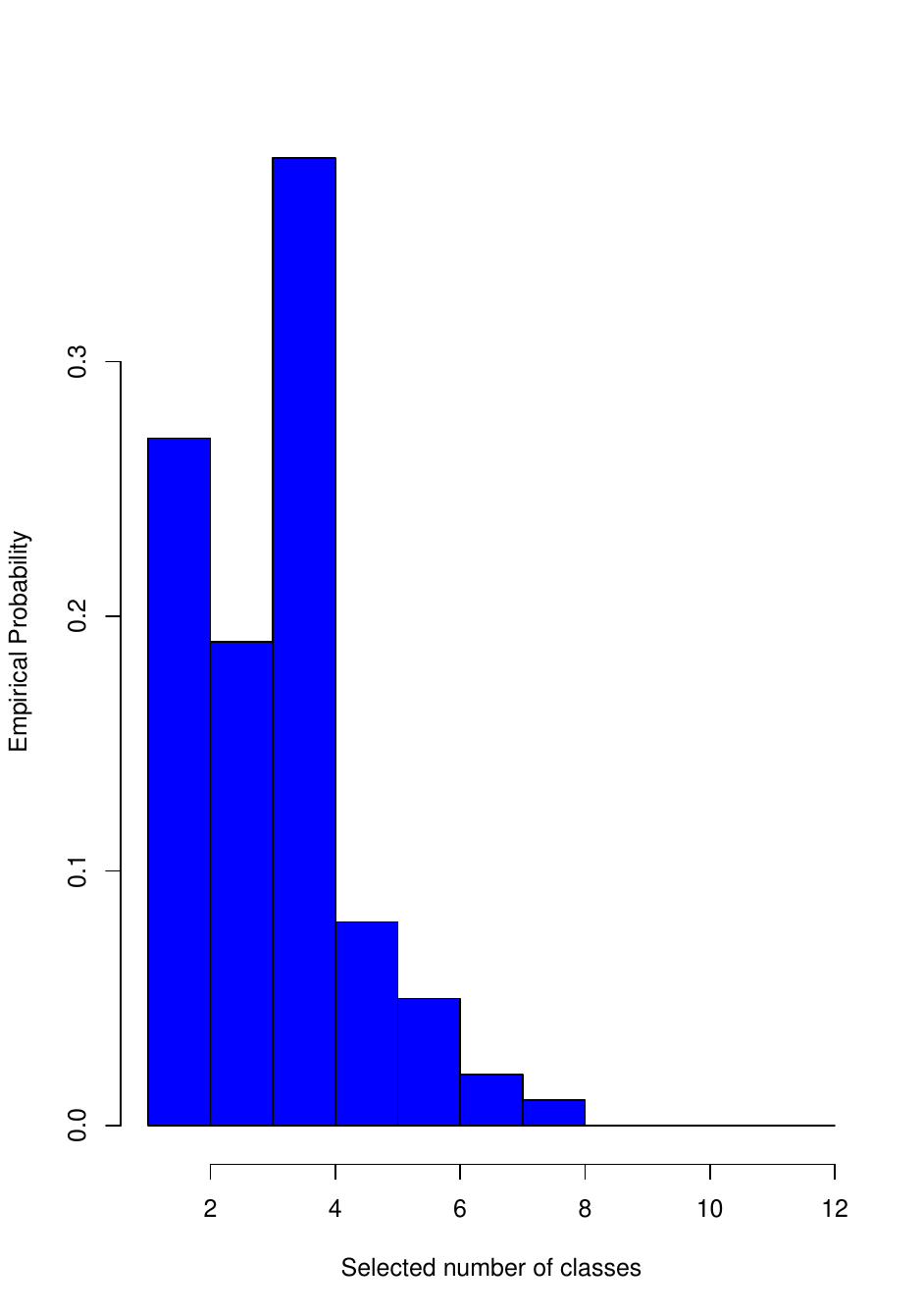}
		\caption{Jump criterion.}	
		\label{subfig_selectedK_Djump_Ethanol_CNLL}
	\end{subfigure}
	\begin{subfigure}{.45\textheight}
		\centering
		\includegraphics[width=.8\textwidth]{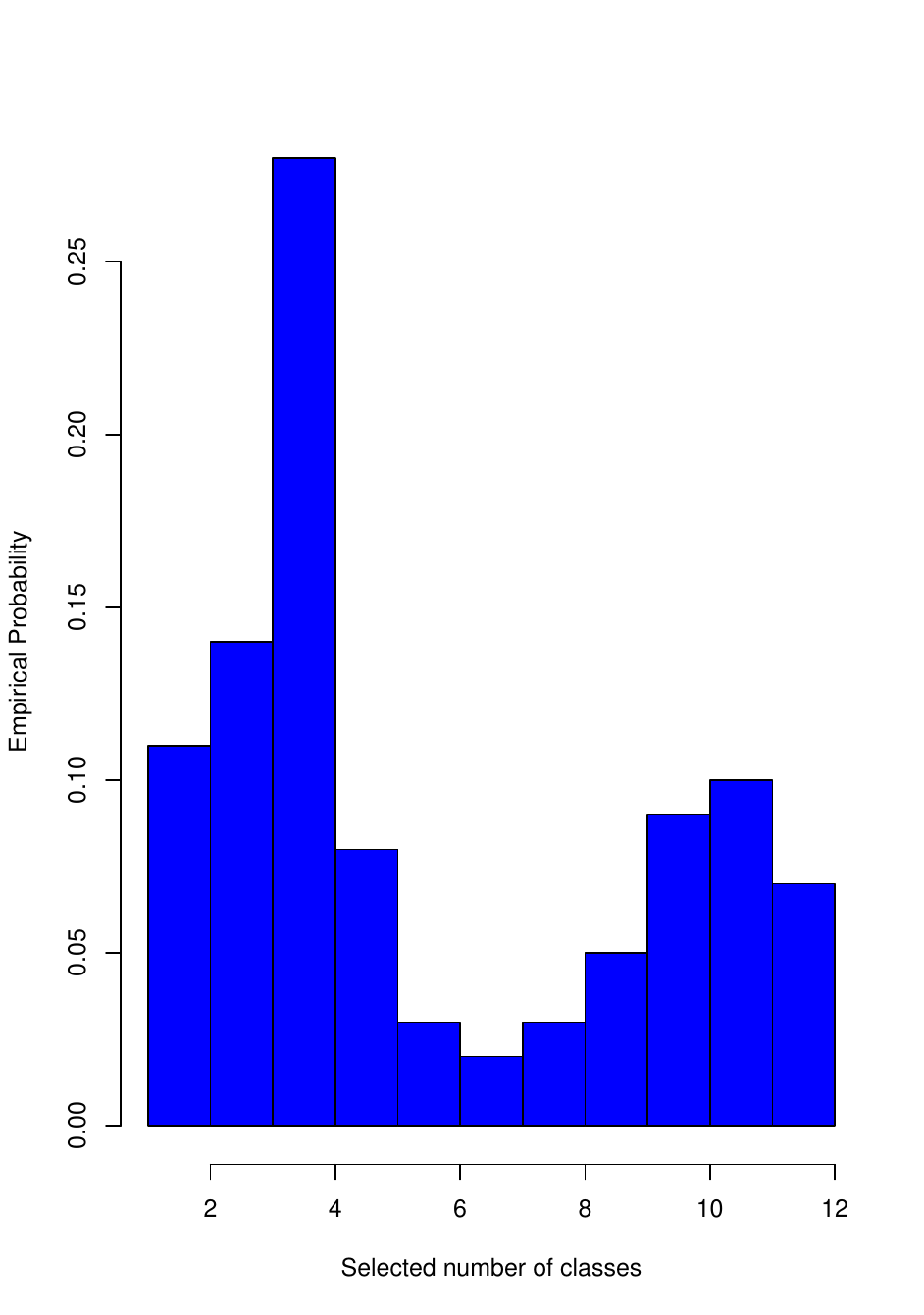}
		\caption{Slope criterion.}	
		\label{subfig_selectedK_DDSE_Ethanol_CNLL}
	\end{subfigure}
	\caption{Histogram of selected $K$ of GLoME on Ethanol data set based on NO using slope heuristic.}
	\label{fig_selectedK_Djump_Ethanol_NO}
\end{figure}

\begin{figure}
	\centering
		\begin{subfigure}{.45\textheight}
				\centering
				\includegraphics[width=0.8\textwidth]{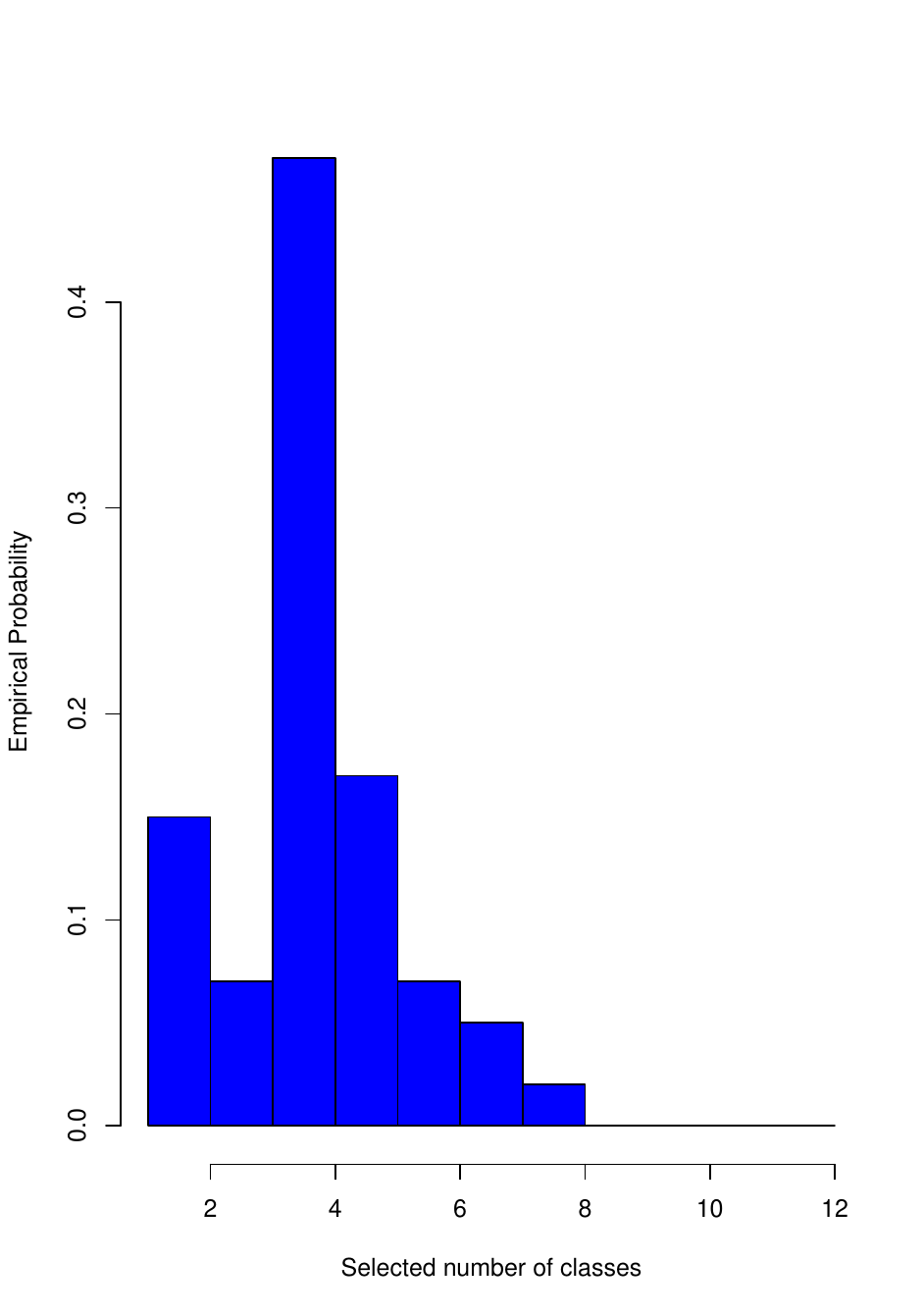}
				\caption{Jump criterion.}
				\label{subfig_selectedK_Djump_Ethanol_inverse_CNLL}
			\end{subfigure}%
		
		\begin{subfigure}{.45\textheight}
			\centering
			\includegraphics[width=0.8\textwidth]{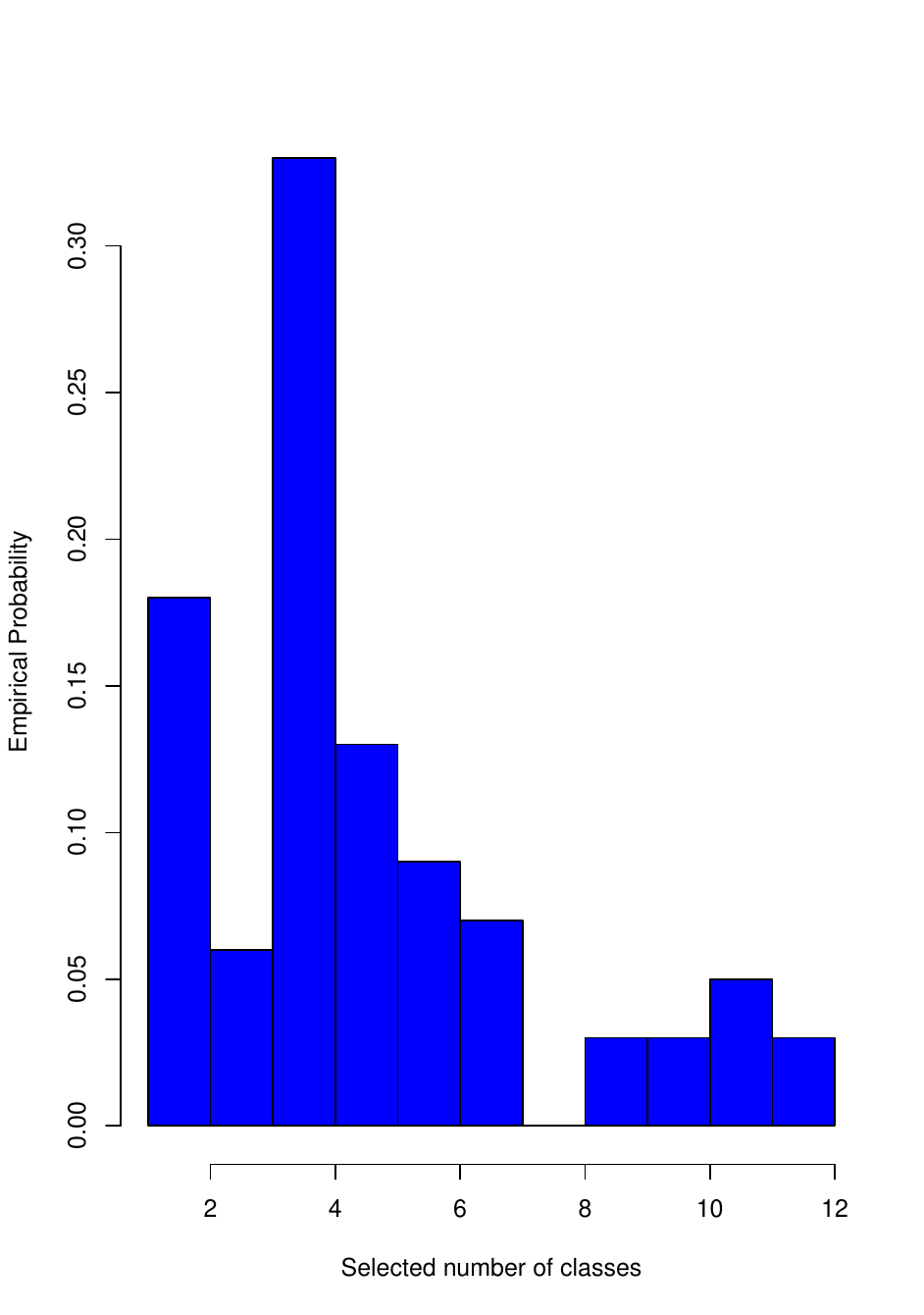}
				\caption{Slope criterion.}
				\label{subfig_selectedK_DDSE_Ethanol_inverse_CNLL}
			\end{subfigure}%
	\caption{Histogram of selected $K$ of GLoME on Ethanol data set based on ER using slope heuristic.}
	\label{fig_selectedK_Djump_Ethanol_ER}
\end{figure}

\begin{figure}
	\centering	\begin{subfigure}{.45\textheight}
		\centering
		\includegraphics[width=0.8\textwidth]{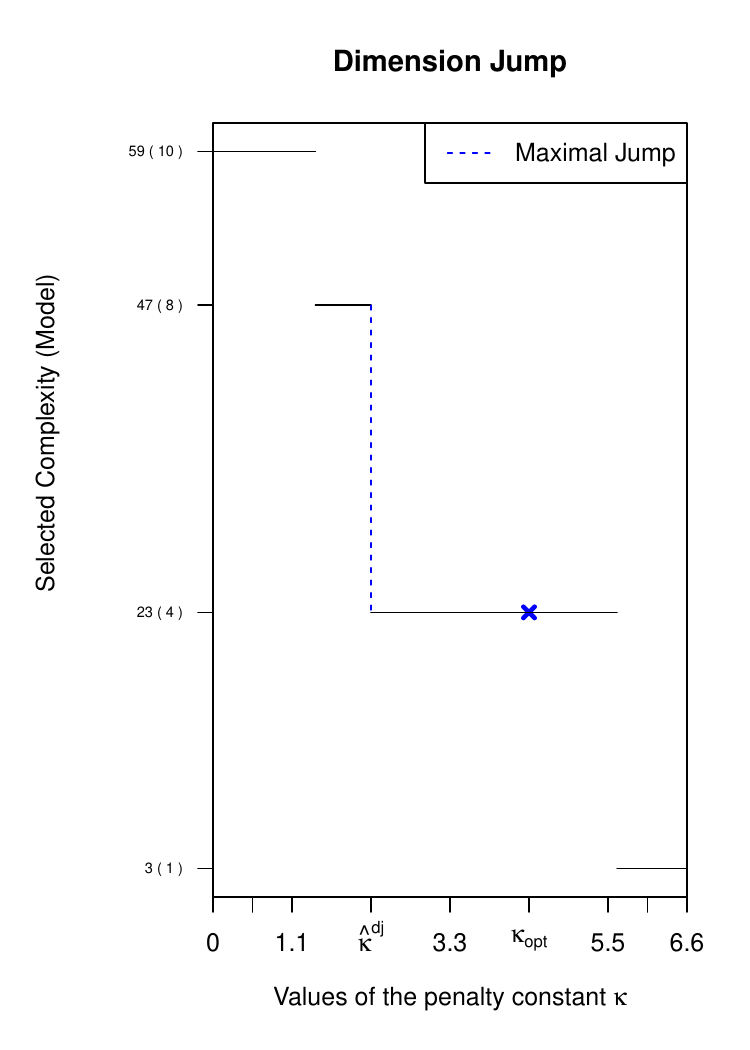}
		\caption{Jump criterion based on NO}	
		\label{subfig_selectedK_Djump_Plot_Ethanol_CNLL}
	\end{subfigure}
	\begin{subfigure}{.45\textheight}
		\centering
		\includegraphics[width=0.8\textwidth]{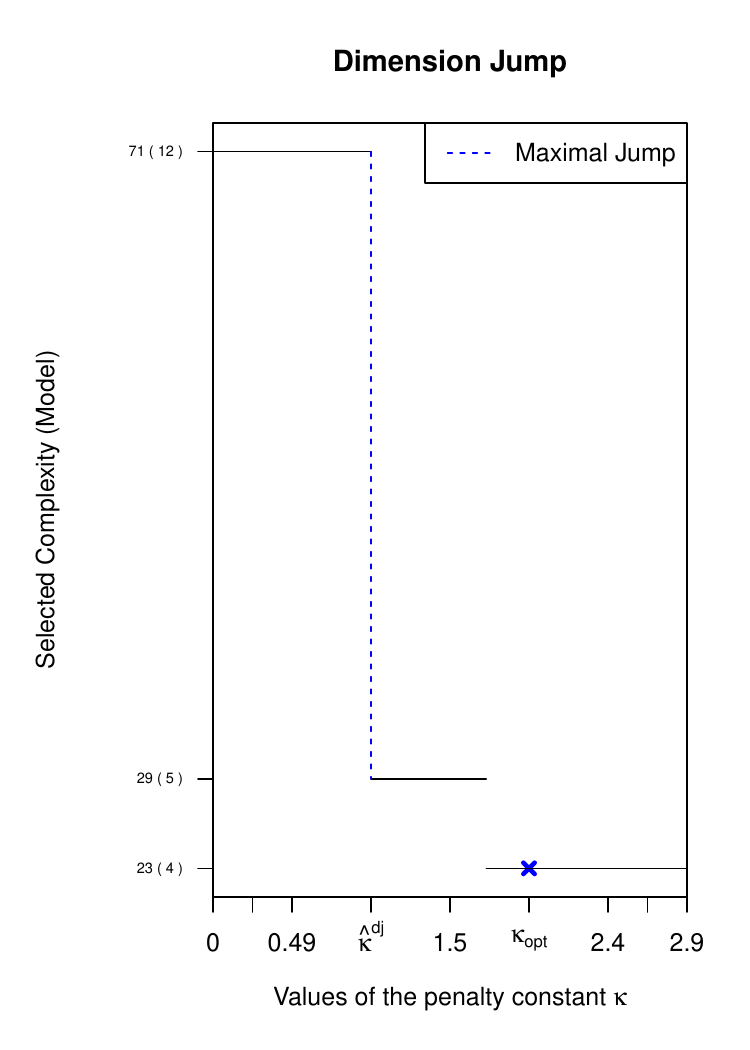}
			\caption{Jump criterion based on ER}	
			\label{subfig_selectedK_Djump_Plot_EquivRatio_NO_Ethanol_CNLL}
		\end{subfigure}
		
	\caption{The jump criterion corresponding to the models chosen with highest empirical probabilities.}
	\label{fig_selectedK_Djump_Ethanol_plot}
\end{figure}

\begin{figure}
	\centering	
		\begin{subfigure}{.4\textheight}
			\centering
			\includegraphics[trim={0 0 11cm 0},clip,width=\linewidth]{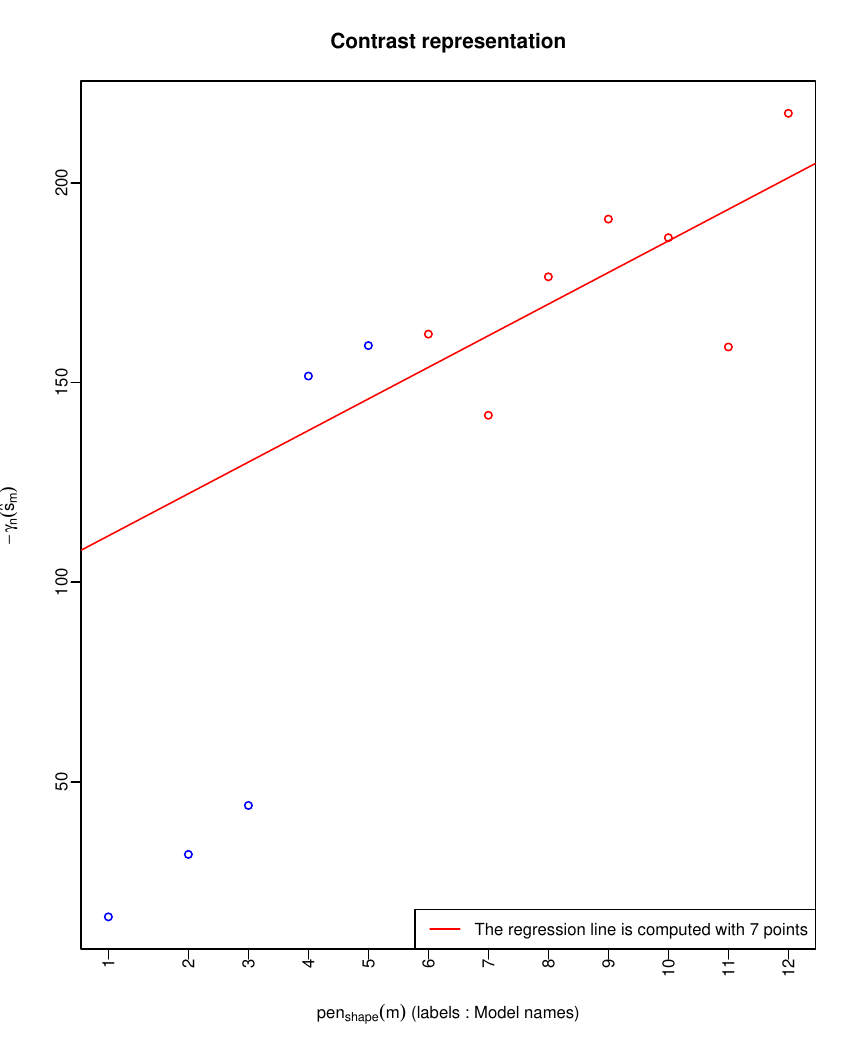}
				\caption{Slope criterion based on NO}
				\label{subfig_selectedK_DDSE_Plot_Ethanol_CNLL}
			\end{subfigure}
		\begin{subfigure}{.4\textheight}
			\centering
			\includegraphics[trim={0 0 11cm 0},clip,width=\linewidth]{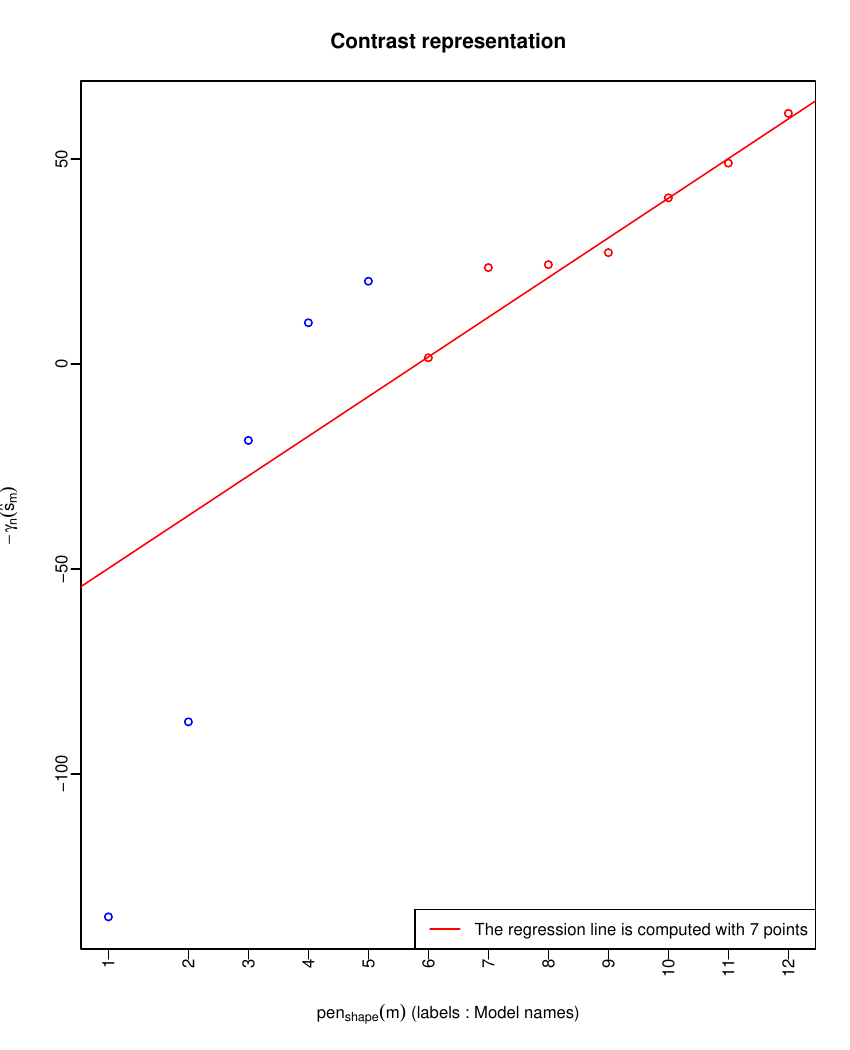}
				\caption{Slope criterion based on ER}
				\label{subfig_selectedK_DDSE_Plot_EquivRatio_NO_Ethanol_CNLL}
			\end{subfigure}%
	\caption{The slope criterion corresponding to the models chosen with highest empirical probabilities.}
	\label{fig_selectedK_DDSE_Ethanol_plot}
\end{figure}

\begin{figure}
	\centering
	\begin{subfigure}{.48\textwidth}
		\centering
		\includegraphics[scale = .3]{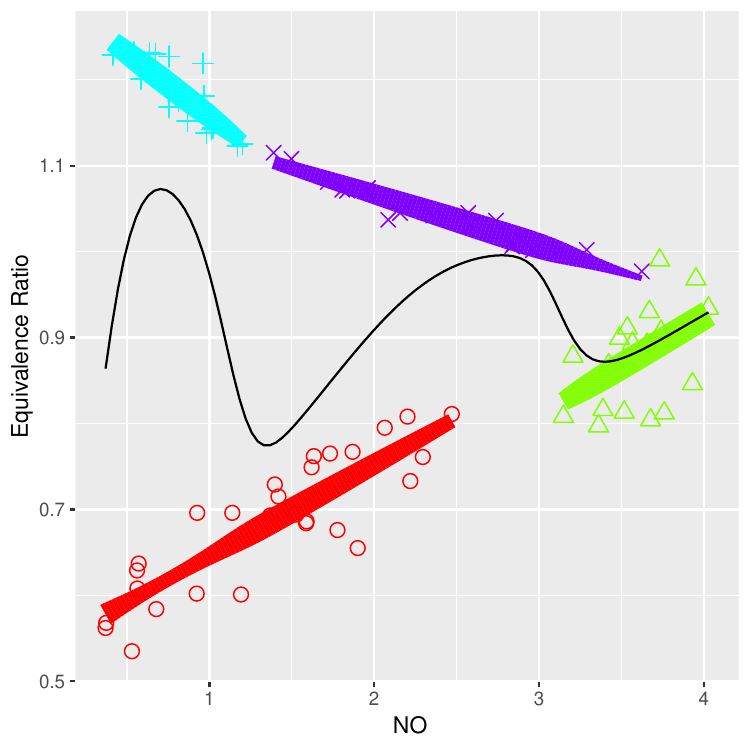}
		\caption{Clustering based on NO}	
		\label{subfig_Ethanol_Clustering_NO_EquivRatio_proportion}
	\end{subfigure}
	\hspace{0.3cm}
	\begin{subfigure}{.48\textwidth}
		\centering
		\includegraphics[scale = .3]{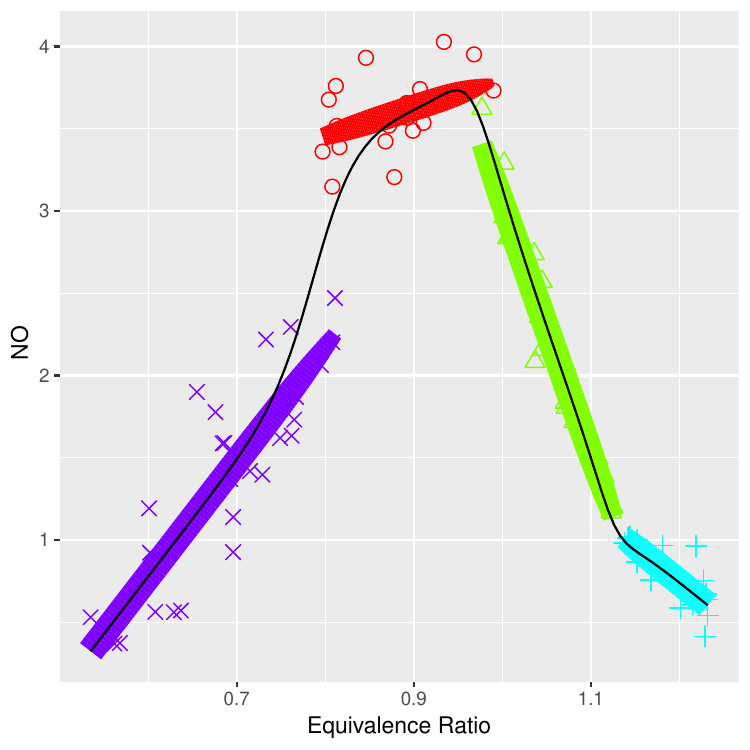}
		\caption{Clustering based on ER}
		\label{subfig_Ethanol_Clustering_EquivRatio_NO_proportion}
	\end{subfigure}%
	
	\begin{subfigure}{.48\textwidth}
		\centering
		\includegraphics[scale = .3]{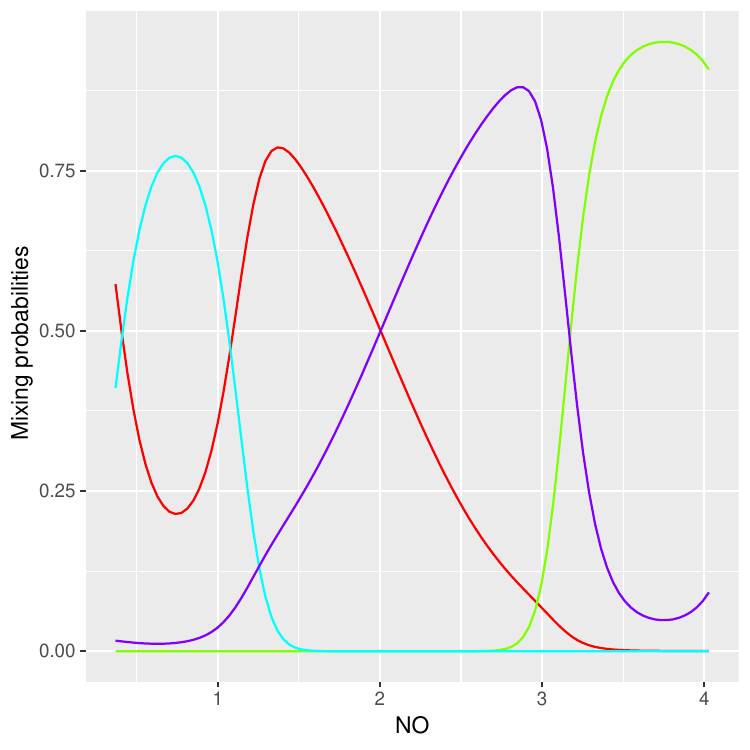}
		\caption{Gating network probabilities based on NO}
		\label{subfig_Ethanol_NO_EquivRatio_Posterior}
	\end{subfigure}%
	\hspace{0.3cm}
	\begin{subfigure}{.48\textwidth}
		\centering
		\includegraphics[scale = .3]{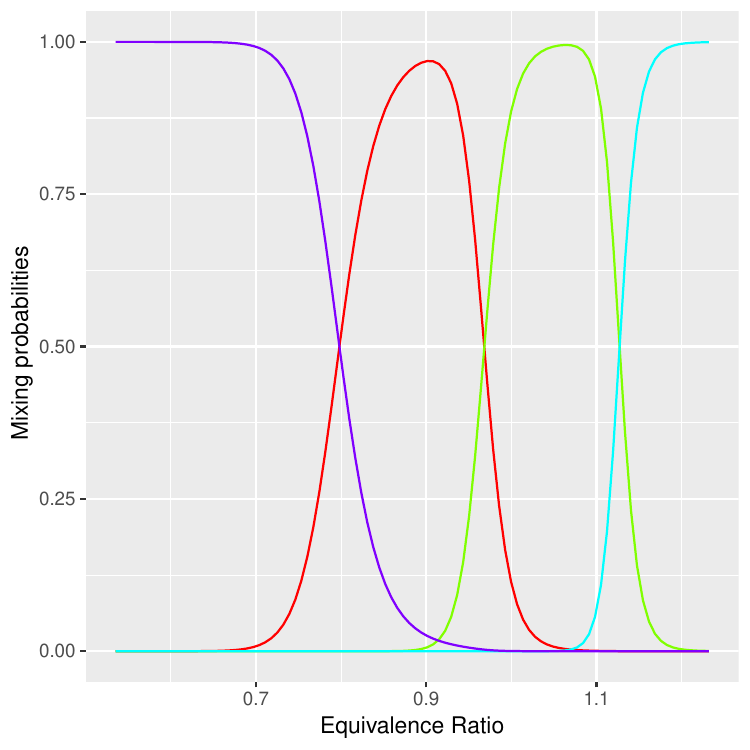}
		\caption{Gating network probabilities based on ER}
		\label{subfig_Ethanol_EquivRatio_NO_Posterior}
	\end{subfigure}%
	
	\caption{Clustering of Ethanol data set. The black curves present the estimated mean functions. The size of the component mean functions corresponds to the posterior mixture proportions.}
	\label{fig_Ethanol_EquivRatio_NO_Posterior}
\end{figure}

More precisely, we first use the EM algorithm from the \href{https://cran.r-project.org/web/packages/xLLiM/index.html}{xLLiM} package to compute the forward PMLE of \eqref{eq_define_GLoME}, for each $K \in [12]$, on the Ethanol data set. Then, based on the slope heuristic (\cref{fig_selectedK_Djump_Ethanol_NO,fig_selectedK_Djump_Ethanol_ER,fig_selectedK_Djump_Ethanol_plot,fig_selectedK_DDSE_Ethanol_plot}), we select the best model. Given the estimators of the model chosen, we obtain the estimated conditional density and clustering by applying the maximum a posteriori probability (MAP) rule (\cref{fig_Ethanol_Conditional_density_EquivRatio_NO,fig_Ethanol_EquivRatio_NO_Posterior}).

Because we only have 88 data points and $6$ parameters per class, the EM algorithm is strongly dependent on the random initialization of the parameters. One solution is that we can modify slightly that procedure in order to guarantee that at least $10$ points are assigned to each class so that the estimated parameters are more stable (cf. \cite{montuelle2014mixture}). In this work, we wish to investigate how well our proposed PMLE performs for detecting the best number of mixture components for the GLoME model. Thus,  we run our experiment over 100 trials with different initializations for the EM algorithm.
Histograms of selected values of $K$ are presented in \cref{fig_selectedK_Djump_Ethanol_NO,fig_selectedK_Djump_Ethanol_ER}. Notice that it is quite unlikely that the true conditional PDF of Ethanol data set belongs to our hypothesised collection of GLoME models. In fact, this phenomenon has been observed in the MS case, \cref{fig_histogramK_Djump_CNLL_WS,fig_histogramK_Djump_CNLL_MS,fig_histogramK_DDSE_CNLL_WS,fig_histogramK_DDSE_CNLL_MS}, on the simulated data set. We believe that this is due to the simplistic affine models used in our experiments. Furthermore, it seems that the jump criterion outperformed the slope criterion in the stability of order selection for GLoME models, as previously observed.

Based on the highest empirical probabilities in all situations, our procedure selects $K=4$ components, which is consistent with the results from \cite{montuelle2014mixture}.
It is worth noting that if we consider the regression of NO with respect to ER, our proposed PMLE of GLoME performs very well for both the clustering and regression tasks (\cref{fig_Ethanol_EquivRatio_NO_Posterior}). Here, instead of considering the variable NO as the covariate, we use it as the response variable. Then, the resulting clustering, the estimated mean function (black curve) and mixing probabilities are more easily interpretable. This is very similar to the results obtained in \cite{montuelle2014mixture}.

\section{Conclusion and perspectives} \label{conclusions}
%
We have studied the PMLEs for GLoME and BLoME regression models. Our main contributions are to establish non-asymptotic risk bounds that take the form of finite-sample weak oracle inequalities, provided that lower bounds on the penalties hold true. We believe that our contributions help to popularize GLoME  and BLoME models as well as slope heuristics, by giving some non-asymptotic theoretical foundations for model selection technique in this area and demonstrating some interesting numerical schemes and experiments. In particular, we aim to provide extensions of the current finite-sample oracle inequalities, \cref{thm_Oracle_Inequality_GLoME,thm_Oracle_Inequality_BLoME}, to more general frameworks where Gaussian experts are replaced by the elliptical distributions in the future work.


%
Some important theoretical issues regarding the tightness of the bounds of the PMLEs are still open. As per the SGaME models from \cite{montuelle2014mixture}, we also have to deal with three potential issues: the differences of divergence on the left ($\jtkl$) and the right ($\tkl$) hand side, $C_1 > 1$, and the relationship between  $\frac{\pen(\bfm)}{n}$ and the variance.  
The first issue is important as in general we have $\jtkl(s_0,s_\bfm ) \le \tkl(s_0,s_\bfm )$. However, \eqref{montuelle.KLandJKL.1} ensures that the two divergences are \emph{equivalent} under regularity conditions. Namely, when 
\begin{align*}
	\sup_{\bfm \in \cM} \sup_{s_\bfm  \in S_\bfm }
	\norm{s_0/s_\bfm }_\infty < \infty.
\end{align*}
Such a strong assumption is satisfied as long as $\cX$ is compact, $s_0$ is compactly supported, and the regression functions are uniformly bounded, and under the condition that there is a uniform lower bound on the eigenvalues of the covariance matrices. 
When $s_0 \notin \left(S_\bfm\right)_{\bfm\in\cM}$, the error bound converges to $C_1 \inf_{s_\bfm  \in S_\bfm } \tkl\left(s_0,s_\bfm \right)$, which may be large. Nevertheless, as we consider GLoME models, some recent results from \cite{nguyen2019approximation,nguyen2020approximationMoE} imply that if we take a sufficiently large number of mixture components, we can approximate a broad class of densities, and thus the term on the right hand side is small for $K$ sufficiently large. This improves the error bound even when $s_0$ does not belong to $S_\bfm$ for any $\bfm \in \cM$. 
%
For the last issue, it holds that $\frac{\pen(\bfm)}{n}$ is approximately proportional to the asymptotic variance, $\frac{\dim\left(S_\bfm \right)}{n}$, in the parametric case, see its justification in \cref{sec_practical_application}. 
Finally, it is also interesting to point out that the weak inequalities in \cref{thm_Oracle_Inequality_GLoME,thm_Oracle_Inequality_BLoME} make use of different divergences and require some regularity assumptions to be considered proper oracle inequalities. To illustrate the strictness of the compactness assumption for $s_0$, we only need to consider $s_0$ as a univariate Gaussian PDF, which obviously does not satisfy such a hypothesis. 
This motivates \cite[Theorem 3.1]{nguyen2020l1oracle} to investigate an $l_1$-ball model selection procedure and an $l_1$-oracle inequality of SGaME with the LASSO estimator.
Such an $l_1$-oracle inequality can be considered as a complementary to \cref{thm_Oracle_Inequality_GLoME}. In particular, the most intriguing property of the $l_1$-oracle inequality is that it requires only the boundedness assumption of the parameters of the model, which is also required in \cref{thm_Oracle_Inequality_GLoME}, as well as in \cite{stadler2010l1,Meynet:2013aa,Devijver:2015aa, nguyen2020l1oracle}. Note that such a mild assumption is quite common when working with MLE (cf. \cite{baudry2009selection,maugis2011non}), to tackle the problem of the unboundedness of the likelihood at the boundary of the parameter space \cite{mclachlan2000finite,redner1984mixture}, and to prevent it from diverging. Nevertheless, by using the smaller divergence: $\jtkl$ (or more strict assumptions on $s_0$ and $s_\bfm $, with the same divergence $\tkl$ as ours), \cref{thm_Oracle_Inequality_GLoME} obtain a faster convergence rate of order $1/n$ compared to $1/\sqrt{n}$ in the $l_1$-oracle inequality. Therefore, in the case where there is no guarantee of a compact support of $s_0$ or uniformly bounded regression functions, the $l_1$-oracle inequality gives a theoretical foundation for the $l_1$-ball model selection procedure, with the order of rate convergence of $1/\sqrt{n}$, with only the boundedness assumption on the parameter space. 

\section*{Acknowledgments}
TrungTin Nguyen is supported by a ``Contrat doctoral'' from the French Ministry of Higher Education
and Research. Faicel Chamroukhi is granted by the French National Research Agency (ANR) grant \href{https://anr.fr/en/funded-projects-and-impact/funded-projects/project/funded/project/b2d9d3668f92a3b9fbbf7866072501ef-f004f5ad27/?tx_anrprojects\_funded\%5Bcontroller\%5D=Funded&cHash=895d4c3a16ad6a0902e6515eb65fba37}{SMILES ANR-18-CE40-0014}. Hien Duy Nguyen is funded by Australian Research Council grant number DP180101192. This research is funded directly by the Inria \href{https://team.inria.fr/statify/projects/lander/}{LANDER} project. TrungTin Nguyen also sincerely acknowledges Inria Grenoble-Rhône-Alpes Research Centre for a valuable Visiting PhD Fellowship working with \href{https://team.inria.fr/statify/}{STATIFY} team so that this research can be completed, Erwan LE PENNEC and Lucie Montuelle for  providing the simulations for the SGaME models.  

\appendix

\section{Our contributions on lemma proofs}\label{proofLemma}

\subsection{Proof of Lemma \ref{lem_inverse_regression_mapping}: Inverse regression trick }\label{proof_1_1_GLLiM}

We  firstly aim to provide the proof of one-to-one correspondence defines the link between the inverse and forward conditional distributions not only for the special case of Gaussian distribution in \eqref{eq_inverse_regression_mapping} but also for elliptical distributions (cf.~\cite{cambanis1981theory,Fang1990}). It is worth mentioning that the multivariate normal distribution, multivariate $t$-distribution and multivariate Laplace distribution
are some instances of elliptical distributions (cf. \cite[Chapter 1]{frahm2004generalized}, \cite{hult2002multivariate}). Note that a statement similar to the following has been proved in the linear regression setting in \cite[Section 2.2]{devijver2020prediction}. We include a proof for mixture of regression models, which is an extension to the aforementioned result.
\subsubsection{Elliptically symmetric distributions} 
We will provide the proof of \cref{lem_inverse_regression_mapping} by using some general results regarding elliptical distributions.

\begin{definition}
	Let $\Xv$ be a $D$-dimensional random vector. Then $\Xv$ is said to be \emph{elliptically distributed} (or simply \emph{elliptical}) if and only if there exist a vector $\mub \in \R^D$, a positive semidefinite matrix $\Sigmab \in \R^{D\times D}$ and a function $\phi: \R^+ \rightarrow \R$ such that the characteristic function $\bft \mapsto \varphi_{\Xv -\mub}(\bft)$ of $\Xv -\mub$ corresponds to $\bft \mapsto \phi\left(\bft^\top\Sigmab \bft\right), \bft\in \R^D$. We write $\Xv \sim \cE_D\left(\mub,\Sigmab,\phi\right)$ to denote that $\mathbf{X}$ is elliptical.
\end{definition}
The function $\phi$ is referred to as the \emph{characteristic generator} of $\Xb$. When $D = 1$ the class of elliptical distributions coincides with the class of univariate symmetric distributions.
Thanks to Proposition 1 from \cite{frahm2004generalized}, it holds that every affinely transformed elliptical random
vector is elliptically distributed. Moreover, the following stochastic representation theorem, \cref{eq_elliptical_Distribution_Theorem}, shows that
the converse is true if the transformation matrix has full rank.
\begin{theorem}[Theorem 1 from \cite{cambanis1981theory}]\label{eq_elliptical_Distribution_Theorem}
	Let $\Xv$ be a $D$-dimensional random vector. Then $\Xv \sim \cE_D\left(\mub,\Sigmab,\phi\right)$ with $\rank(\Sigmab) = k$ if and only if 
	\begin{align*}
		\Xv = \mub + \cR \Lambdab \Ub^{(k)},
	\end{align*}
	where $\Ub^{(k)}$ is a $k$-dimensional random vector uniformly distributed on the unit hypersphere with $k-1$ dimensions $\cS^{k-1} = \left\{\xb\in\R^D: \left\|\xb\right\|_2 = 1\right\}$, $\cR$ is a non-negative random variable with distribution function $F$ related to $\phi$ being stochastically independent of $\Ub^{(k)}$, $\mub \in \R^D$ and $\Sigmab = \Lambdab \Lambdab^\top$ is a rank factorization of $\Sigmab$ where $\Lambdab \in \R^{D\times k}$ with $\rank(\Lambdab) = k$.
\end{theorem}

Note that via the transformation matrix $\Lambdab$, the spherical random vector $\Ub^{(k)}$ produces elliptically
contoured density surfaces, whereas the generating random variable $\cR$ determines the distribution’s shape, in particular the heaviness of the distribution's tails. Further, $\mub$ determines the
location of the random vector $\Xv$. The matrix $\Sigmab$ is called the \emph{dispersion matrix} or \emph{scatter matrix} of $\Xv$. Therefore, it holds that every elliptical distribution belongs to a location-scale-family \citep{kelker1970distribution} defined by an underlying
spherical standard distribution.

\begin{example}[Multivariate normal distribution]\label{example_mulGaussian_EllipticalDistribution}
	Let $\mub \in \R^D$ and $\Lambdab \in \R^{D \times k}$ such that $\Sigmab := \Lambdab \Lambdab^\top \in \R^{D\times D}$ is positive definite. The random vector $\Xb \sim \Phi_D \left(\mub,\Sigmab\right)$ is elliptically
	distributed since $\Xv$ is representable as
	$\Xb =\mub +\sqrt{\chi^2_k} \Lambdab \Ub^{(k)}$. The underlying spherical standard distribution is the
	standard normal distribution. Further, since $s \mapsto \exp(-s/2)$ is the characteristic
	generator for the class of normal distributions, the characteristic function of $\Xb - \mub $ corresponds to $\bft \mapsto \phi_{\Xb -\mub}(\bft) = \exp\left(-\bft^\top \Sigmab \bft\right), \bft\in \R^d$.
\end{example}

We next describe some important results on the conditional distributions of elliptical random
vectors (cf. \cite[Corollary 5]{cambanis1981theory}, \cite[Chapter 1]{frahm2004generalized}).
\begin{theorem}\label{thm_properties_elliptical_distribution}
	Let $\Xv \sim \cE_D\left(\mub,\Sigmab,\phi\right)$ with $\rank(\Sigmab) = k$. It holds that:
	\begin{enumerate}[(a)]
		\item $\Ep{X} = \mub$.
		\item $\varp{\Xv} = \frac{\Ep{\cR^2}}{k}\Sigmab = -2 \phi'(0)\Sigmab$, if $\phi$ is differentiable at $0$.
		
		\item The sum of independent elliptically distributed random vector with the same dispersion matrix $\Sigmab$ is elliptically too. Furthermore, the sum of two dependent elliptical random vectors with the same dispersion matrix, which are dependent only through their radial parts $\cR$, is also elliptical \citep[Theorem 4.1]{hult2002multivariate}. More precisely, let $\cR$ and $\widetilde{\cR}$ be nonnegative random variables and let $\Xv := \mub + \cR \Zb \sim \cE_D\left(\mub,\Sigmab,\phi\right)$ and $\widetilde{\Xv} := \widetilde{\mub} + \widetilde{\cR} \widetilde{\Zb} \sim \cE_D\left(\widetilde{\mub},\Sigmab,\widetilde{\phi}\right)$, where $\left(\cR,\widetilde{\cR}\right)$, $Z$, and $\widetilde{Z}$ are independent. Then $\Xv + \widetilde{\Xv} \sim \cE_D\left(\mub + \widetilde{\mub},\Sigmab,\phi^*\right)$.
		
		\item Affine transformation: every affinely transformed and particularly every linearly combined elliptical random
		vector is elliptical, too. More formally, for any $\bfb \in \R^L$, $\bfA \in \R^{L\times D}$, and $\Yv = \bfb+ \bfA\Xv$ where $\Xv =\mub +\cR 
		\Lambdab \Ub^{(k)}$ with $\Lambdab \in \R^{D \times k}$, it follows that $\Yv 
		= \cE_L\left(\bfb +\bfA\mub,\bfA\Sigmab \bfA^\top,\phi\right)$ since
		\begin{align} \label{eq_Elliptical_Distribution_affine}
			\Yv  =\bfb+ \bfA\left(\mub +\cR 
			\Lambdab \Ub^{(k)}\right) = \left(\bfA\mub + \bfb\right)	+\cR 
			\bfA\Lambdab \Ub^{(k)}.
		\end{align}  
		\item Marginal distribution: let $\cP_m \in \left\{0,1\right\}^{m\times D} (m \le D)$ be a \emph{permutation and deletion} matrix, \ie $\cP_m$ has only binary
		entries of $0$’s and $1$’s and $\cP_m \cP_m^\top = {\id}_{m}$. 
		So the transformation $\cP_m \Xv =: \Yv$ permutes and
		deletes certain components of $\Xv$ such that $\Yv$ is a $k$-dimensional random vector containing
		the remaining components of $\Xv$ and having a (multivariate) marginal distribution with
		respect to the joint distribution of $\Xv$. Then by \eqref{eq_Elliptical_Distribution_affine}, we obtain $\Yv \sim \cE_m\left( \cP_m\mub,\cP_m \Sigmab \cP_m^\top,\phi\right)$ since 
		\begin{align} \label{eq_Elliptical_Distribution_marginal}
			\Yv  = \cP_m\left(\mub +\cR 
			\Lambdab \Ub^{(k)}\right) = \cP_m\mub +\cR 
			\cP_m\Lambdab \Ub^{(k)}.
		\end{align} 
		\item Conditional distribution: let $\Xv = \left(\Xv_1, \Xv_2\right)$, where $\Xv_1$ is a $k$-dimensional sub-vector of $\Xv$, and let $\Sigmab =\begin{pmatrix}
			\Sigmab_{11}&\Sigmab_{12}\\
			\Sigmab_{21}&\Sigmab_{22}
		\end{pmatrix}\in\R^{D\times D}$. Provided the conditional random vector $\Xv_2\mid \Xv_1 = \xv_1$ exists, it is also elliptically distributed: $\Xv_2\mid (\Xv_1 = \xv_1)\sim \cE_{D-k}\left(\mub^*,\Sigmab^*,\phi^*\right)$. Moreover, it can be presented stochastically by 
		\begin{align*}
			\Xv_2\mid \left(\Xv_1 = \xv_1\right)=\mub^*+ \cR^* \Ub^{(D-k)} \Gamma^*
		\end{align*}
		and $\Ub^{(D-k)}$ is uniformly distributed on $\cS^{(D-k-1)}$, and 
		\begin{align*}
			\cR^* &= \cR \sqrt{1-\beta}\mid  \left(\cR \sqrt{\beta}\Ub^{(k)} = \Sigmab^{-1}_{11}\left(\xv_1-\mub_1\right)\right),\\
			\mub^* &= \mub_2 + \Sigmab_{21} \Sigmab^{-1}_{11} \left(\xv_1 - \mub_1\right), \\
			\Sigmab^* &= \Sigmab_{22} - \Sigmab_{21} \Sigmab^{-1}_{11} \Sigmab_{12}.
		\end{align*}
		where $\beta \sim \Beta \left(k/2,(D-k)/2\right)$ and $\cR,\beta,\Ub^{(k)}$ and $\Ub^{(D-k)}$ are mutually independent, and $\Sigmab^* = \left(\Gamma^*\right)^\top \Gamma^*$. 
	\end{enumerate}
	
\end{theorem}
\subsubsection{Relation between forward and inverse regression}
Proposition 1 from \cite{deleforge2015high}, a multivariate extension of \cite{ingrassia2012local}, leads to a link between GLLiM, defined in \eqref{eq_inverseGLLiM} models, and a Gaussian mixture model on the joint variable $\left[\Xv;\Yv\right]$. This result motivates us to establish the general proof for the relationship between forward and inverse mixture of elliptical regression models. More precisely, we consider the following generative model, conditionally on the cluster label:
\begin{eqnarray}
	\left[\Xv;\Yv\right] \mid  \left(Z=k\right) \sim \cE_{L+D} \left(\bfm_k, \bfV_k, \phi\right), \label{def_joint} 
\end{eqnarray}
where $\cE_{L+D}$ denotes an elliptical distribution of dimension $D+L$, and are $\bfm_k$ and $\bfV_k$ its location and scale parameters, respectively. 

When applying the inverse regression strategy in the context of mixture of elliptical locally-linear mapping, the key point is to account for \eqref{eq_locallyaffine_inverse}:
\begin{align}\label{eq_locallyaffine_inverse}
	\Xv = \sum_{k=1}^K \Indi \left(Z = k\right) \left(\bfA_k\Yv +\bfb_k + \bfE_k\right),
\end{align}
where $\bfA_k \in \R^{D\times L}$ and vector $\bfb_k \in \R^D$, and $\bfE_k$ is an error term capturing both the reconstruction error due to the local affine approximation and the observation noise in $\R^D$, into the parameterization of $\bfm_k$ and $\bfV_k$. Given $Z = k$, it follows from \eqref{def_joint} that $\Yv$ is also elliptical distribution by using \cref{thm_properties_elliptical_distribution} (e) and $\Yv$ can be assumed to have a location $\bfc_k \in \R^L$ and a scale matrix $\Gammab_k \in \R^{L \times L}$. We assume further that the error term $\bfE_k \sim \cE\left(\zero,\phi_{e_k},\Sigmab_k\right)$ is an unobserved centered elliptical random noise with residual covariance matrix $\Sigmab_k$ of type $\phi_{e_k}$, and is independent of $\Yv$. Then, using \eqref{eq_locallyaffine_inverse} and \cref{thm_properties_elliptical_distribution}, we have
\begin{equation}
	\begin{array}{rl}
		\Ep{\Xv\mid \left(Z=k\right)} =& \Ep{\bfA_k\Yv +\bfb_k +\bfE_k} = \bfA_k\bfc_k+\bfb_k,\\
		\varp{\Xv\mid \left(Z=k\right)} =& 
		\varp{\bfA_k\Yv} + \varp{\bfE_k} = \bfA_k\Gammab_k \bfA_k^\top +\Sigmab_k,\\
		\covp{\Xv,\Yv\mid \left(Z=k\right)}=& 
		\bfA_k\Gammab_k,\covp{\Yv,\Xv}= \covp{\Yv,\bfA_k\Yv} = \Gammab_k^\top\bfA_k^\top,\\
		\bfm_k=&
		\begin{bmatrix}
			\bfA_k\bfc_k+\bfb_k\\
			\bfc_k
		\end{bmatrix},  \vspace{0,1cm} \\
		\bfV_k=&
		\begin{bmatrix}
			\Sigmab_k+\bfA_k\Gammab_k \bfA_k^\top & \bfA_k\Gammab_k \\
			(\bfA_k\Gammab_k)^\top & \Gammab_k
		\end{bmatrix}. \label{eq_JGMMtoGLM}
	\end{array}
\end{equation}


Note that in the forward and inverse regression problems of elliptical locally-linear mapping (containing the Gaussian case \eqref{eq_conditional_forward}, \eqref{eq_marginal_forward}, \eqref{eq_conditional_inverse} and \eqref{eq_marginal_inverse}), by using \cref{thm_properties_elliptical_distribution}, the joint distribution defined in \eqref{def_joint} allows us to consider a mixture of linear regression problem, characterized by the following marginal and conditional distributions:
\begin{align}
	\Xv\mid \left(Z=k\right) &\sim \cE_{D} \left(\bfc_k^*,\Gammab_k^*,\phi\right),\label{eq_marginal_forward_distribution}\\ 
	\Yv\mid \left(\Xv,Z=k\right) & =  \bfA^*_k\Xv +\bfb^*_k+ \bfE^*_k,\label{eq_conditional_forward_distribution}\\
	\Yv\mid \left(Z=k\right) &\sim \cE_{D} \left(\bfc_k,\Gammab_k,\phi\right),\label{eq_marginal_inverse_distribution}\\ 
	\Xv\mid \left(\Yv,Z=k\right) & =  \bfA_k\Yv +\bfb_k+ \bfE_k, \label{eq_conditional_inverse_distribution}
\end{align} 
where $\bfE_k \sim \cE\left(\zero,\phi_{e_k},\Sigmab_k\right)$ and $\bfE^*_k \sim \cE\left(\zero,\phi_{e^*_k},\Sigmab^*_k\right)$.

Then, we claim that the joint distribution, defined in \eqref{def_joint} and \eqref{eq_JGMMtoGLM}, leads to the marginal and the conditional distributions of \crefrange{eq_marginal_forward_distribution}{eq_conditional_inverse_distribution} and to a mapping between their mean and variance parameters. Indeed, by using conditioning properties of elliptical distributions, see more in \cref{thm_properties_elliptical_distribution}, implies the following marginal for $\Xv$ and conditional distribution for $\Yv$ given $\Xv$ as follows:
\begin{align}
	\Xv\mid \left(Z=k\right) &\sim \cE_{D} \left(\bfA_k\bfc_k+\bfb_k,\Sigmab_k+\bfA_k\Gammab_k \bfA_k^\top,\phi\right),\label{eq_forward_marginal_distribution}\\
	\Yv\mid \left(\Xv,Z=k\right) &\sim \cE_{L} \left(\mb^{yx}_k,\Sigmab^{yx}_k,\widetilde{\phi}\right) \label{eq_forward_conditional_distribution},
\end{align}
where the explicit expression of the characteristic function $\widetilde{\phi}$ can be found in \cite[Corollary 5]{cambanis1981theory}, and
\begin{align*}
	\mb^{yx}_k &=\bfc_k+\Gammab_k^\top \bfA_k^\top\left(\Sigmab_k+\bfA_k\Gammab_k \bfA_k^\top\right)^{-1}\left(\Xv-\bfA_k\bfc_k-\bfb_k\right),\\
	\Sigmab^{yx}_k &= \Gammab_k - \Gammab_k^\top \bfA_k^\top \left(\Sigmab_k+\bfA_k\Gammab_k \bfA_k^\top\right)^{-1}\bfA_k \Gammab_k = 	\left(\Gammab_k^{-1}+\bfA_k^\top\Sigmab_k^{-1}\bfA_k\right)^{-1}, 
\end{align*}
with the fact that the last equality is the Woodbury matrix identity. Note that the locations and scale matrices of the conditional distribution do not depend upon the third parameter of the joint distribution, and consequently, we do not describe the explicit expression for $\widetilde{\phi}$. We then utilize again the Woodbury matrix identity and the symmetric property of $\Gammab$ to identify \eqref{eq_marginal_forward_distribution} and \eqref{eq_conditional_forward_distribution} with \eqref{eq_forward_marginal_distribution} and \eqref{eq_forward_conditional_distribution}, respectively, which implies the following important connections:
\begin{align*}
	\bfc_k^* &= \bfA_k\bfc_k+\bfb_k,\Gammab_k^* = \Sigmab_k+\bfA_k\Gammab_k \bfA_k^\top, \Sigmab^*_k = \left(\Gammab_k^{-1}+\bfA_k^\top\Sigmab_k^{-1}\bfA_k\right)^{-1},\\
	\bfA_k^* &=  \Gammab_k^\top \bfA_k^\top\left(\Sigmab_k+\bfA_k\Gammab_k \bfA_k^\top\right)^{-1}\nn\\
	&= \Gammab_k^\top \bfA_k^\top\Sigmab^{-1}_k - \Gammab_k^\top \bfA_k^\top\Sigmab^{-1}_k  \bfA_k \left(\Gammab_k^{-1}+\bfA_k^\top\Sigmab^{-1}_k \bfA_k\right)^{-1}\bfA_k^\top \Sigmab^{-1}_k\nn\\
	&= \left[\Gammab_k\left(\Gammab_k^{-1}+\bfA_k^\top\Sigmab^{-1}_k \bfA_k\right)  - \Gammab_k \bfA_k^\top\Sigmab^{-1}_k  \bfA_k \right] \left(\Gammab_k^{-1}+\bfA_k^\top\Sigmab^{-1}_k \bfA_k\right)^{-1}\bfA_k^\top \Sigmab^{-1}_k \nn\\
	&=\left(\Gammab_k^{-1}+\bfA_k^\top\Sigmab^{-1}_k \bfA_k\right)^{-1}\bfA_k^\top \Sigmab^{-1}_k = \Sigmab^*_k\bfA_k^\top \Sigmab^{-1}_k,\nn\\
	\bfb_k^* &= \bfc_k+\Gammab_k^\top \bfA_k^\top\left(\Sigmab_k+\bfA_k\Gammab_k \bfA_k^\top\right)^{-1}\left(-\bfA_k\bfc_k-\bfb_k\right)\nn\\
	& = \bfc_k+ \left(\Gammab_k^{-1}+\bfA_k^\top\Sigmab^{-1}_k \bfA_k\right)^{-1}\bfA_k^\top \Sigmab^{-1}_k \left(-\bfA_k\bfc_k-\bfb_k\right)\nn\\
	& = \left(\Gammab_k^{-1}+\bfA_k^\top\Sigmab^{-1}_k \bfA_k\right)^{-1} \left[\left(\Gammab_k^{-1}+\bfA_k^\top\Sigmab^{-1}_k \bfA_k\right) \bfc_k+ \bfA_k^\top \Sigmab^{-1}_k \left(-\bfA_k\bfc_k-\bfb_k\right)\right]\nn\\
	&= \Sigmab_k^* (\Gammab_k^{-1}\bfc_k-\bfA_k^\top\Sigmab_k^{-1}\bfb_k).
\end{align*}
Therefore, the desired results then are obtained by using the fact that the multivariate normal distribution not only has the property of having Gaussian marginal and conditional
distributions \citep[Sections 2.3.1 and 2.3.2]{bishop2006pattern} but also belongs to elliptical distributions (detailed in \cref{example_mulGaussian_EllipticalDistribution}). Furthermore, it should be stressed that several versions of multivariate $t$-distributions (\eg Section 5.5, page 94 of \cite{kotz2004multivariate}, \cite{ding2016conditional}) have the previous property. This leads to the inverse regression model based on the multivariate $t$-distributions \citep{perthame2018inverse}.
It will be interesting to find other sub-classes of elliptically contoured distributions that have the closedness property on marginal and conditional distribution so that the previous inverse regression trick can be applied. 

\subsection{Proof of Lemma \ref{lem_Bracketing_Entropy_Gates_intro_mine}: Bracketing entropy of Gaussian gating networks} \label{sec_proof_lem_Bracketing_Entropy_Gates_intro_mine}

Note that the first inequality of \cref{lem_Bracketing_Entropy_Gates_intro_mine} comes from \cite[Lemma 4]{montuelle2014mixture} and describes relationship between the bracketing entropy of $\cP_K$ and the entropy of $\cW_K$. Therefore, \cref{lem_Bracketing_Entropy_Gates_intro_mine} is obtained by proving that there exists a constant $C_{\cW_K}$ such that $\forall \delta \in (0,2]$,
\begin{align}\label{eq_assumption_DIM_intro}
	\entropy(\delta,\cW_K) &\le \dim\left(\cW_K\right) \left(C_{\cW_K} + \ln \left(\frac{1}{\delta}\right)\right),
\end{align}
where $\dim\left(\cW_K\right) = K-1 + KL + K \frac{L(L+1)}{2}$.

In order to establish the proof for \eqref{eq_assumption_DIM_intro}, we have to construct firstly the $\delta_\pib$-covering $\Pib_{K-1,\omegab}$ of $\Pib_{K-1}$ via \cref{lem_MetricEntropySimplex_maximum_norm_intro}, which is proved in \cref{sec_proof_lem_MetricEntropySimplex_maximum_norm_intro}.
\begin{lem}[Covering number of probability simplex with maximum norm]\label{lem_MetricEntropySimplex_maximum_norm_intro}
	Given any $\delta_{\pib} >0$, any $\pib \in \Pib_{K-1}$, we can choose $\widehat{\pib}\in \Pib_{K-1,\omegab}$, an $\delta_\pib$-covering of $\Pib_{K-1}$, so that $\max_{k\in[K]}\left|\pi_k - \widehat{\pib}_k\right| \le \delta_{\pib}$. Furthermore, it holds that
	\begin{align}
		\cN\left(\delta_{\pib},\Pib_{K-1},\left\|\cdot\right\|_\infty\right) \le \frac{K \left(2 \pi e\right)^{K/2}}{\delta^{K-1}_{\pib}}. \label{eq.MetricEntropySimplex.full_intro}
	\end{align}
\end{lem}

Then, by definition of the covering number, \eqref{eq_assumption_DIM_intro} is obtained immediately via \cref{lem_covering_number_WK_intro}, which controls the covering number of $\cW_K$ and is proved in \cref{sec_proof_lem_covering_number_WK_intro}.

\begin{lem}\label{lem_covering_number_WK_intro}
	Given a bounded set $\cY$ in $\R^L$ such that $\mathcal{\cY} = \left\{\yv \in \R^L: \left\|\yv\right\|_\infty \le C_\cY \right\}$, it holds that $\cW_K$ has a covering number satisfied $\cN\left(\delta,\cW_K,\dsup \right) \le C \delta^{-\dim\left(\cW_K\right)}$, for some constant $C$.
\end{lem}
Indeed, \cref{lem_covering_number_WK_intro} implies the desired result by noting that 
\begin{align*}
	\entropy\left(\delta,\cW_K \right)
	&= \ln \cN\left(\delta,\cW_K,d_{\norm{\sup}_\infty} \right)
	\le  \ln\left[\frac{C}{\delta^{\dim\left(\cW_K\right)}}\right]\nn\\
	&=  \dim\left(\cW_K\right) \left[\frac{1}{\dim\left(\cW_K\right)}\ln C+ \ln \left(\frac{1}{\delta}\right)\right]\nn\\
	&= \dim\left(\cW_K\right) \left(C_{\cW_K} + \ln \left(\frac{1}{\delta}\right)\right).
\end{align*}
\subsubsection{Proof of Lemma \ref{lem_MetricEntropySimplex_maximum_norm_intro}}\label{sec_proof_lem_MetricEntropySimplex_maximum_norm_intro}
Note that \citep[Lemma 2]{genovese2000rates} provide a result for controlling a $\delta_\pib$-Hellinger bracketing of $\Pib_{K-1}$. However, such result can not be applied for our  \cref{lem_MetricEntropySimplex_maximum_norm_intro} since they use $\delta_\pib$-Hellinger bracketing entropy while we use $\delta_\pib$-covering number for the probability complex with maximum norm.

Given any $\pib = \left(\pi_k\right)_{k\in[K]} \in \Pib_{K-1}$, let $\xib = \left(\xib_k\right)_{k\in [K]}$ where $\xib_k = \sqrt{\pi_k},\forall k \in [K]$. Then $\pib\in\Pib_{K-1}$ if and only if $\xib \in Q^+ \cap U$, where $U$ is the surface of the unit sphere and $Q^+$ is the positive quadrant of $\R^K$. Next, we divide the unit cube in $\R^K$ into disjoint cubes with sides parallel to the axes and sides of length $\delta_\pib/\sqrt{K}$. Let $ \left(\bfC_j\right)_{j\in [N]}$ is the subset of these cubes that have non-empty intersection with $Q^+ \cap U$. For any $j \in [N]$, let $\nub_j = \left(\nub_{j,k}\right)_{k\in[K]}$ be the center of the cube $\bfC_j$ and $\nub^2_j = \left(\nub^2_{j,k}\right)_{k\in[K]}$ . 

Then $\left\{\nub_j\right\}_{j \in [N]}$ is a $\delta_{\pib}/\left(2\sqrt{K}\right)$-covering of $Q^+ \cap U$, since we have for any $\xib = \left(\xib_k\right)_{k\in [K]} \in Q^+ \cap U$, there exists $j_0 \in [N]$ such that $\xib \in \bfC_{j_0}$, and 
\begin{align}
	\left\|\xib-\nub_{j_0}\right\|_\infty= \max_{k\in[K]}\left|\xib_k-\nub_{j_0,k}\right| \le \frac{\delta_{\pib}}{2\sqrt{K}}.\label{eq.defineCoveringUnitSphere.full_intro}
\end{align}
Therefore, it follows that $\Pib_{K-1,\omegab}:=\left\{\nub^2_j\right\}_{j \in [N]}$ is a $\delta_{\pib}$-covering of $\Pib_{K-1}$, since for any $\pib = \left(\pi_k\right)_{k\in[K]} \in \Pib_{K-1}$, \eqref{eq.defineCoveringUnitSphere.full_intro} leads to the existence of $j_0\in[N]$, such that
\begin{align*}
	\left\|\pib-\nub^2_{j_0}\right\|_\infty&=\max_{k\in[K]}\left|\xib^2_k-\nub^2_{j_0,k}\right|= \max_{k\in[K]}\left\{\left|\xib_k-\nub_{j_0,k}\right|\left|\xib_k+\nub_{j_0,k}\right|\right\} \nn\\
	&\quad \le \frac{\delta_{\pib}}{2\sqrt{K}} \max_{k\in[K]}\left|\xib_k+\nub_{j_0,k}\right|\le \frac{\delta_{\pib}}{\sqrt{K}} \le   \delta_{\pib},
\end{align*}
where we used the fact that $\max_{k\in[K]}\left|\xib_k+\nub_{j_0,k}\right| \le 2$. Now, it remains to count the number of cubes $N$. Let $\cT_a = \left\{\zb \in Q^+: \left\|\zb\right\|_2 \le a\right\}$ and let $\cC = \bigcup_{ j \in [N]}\bfC_j$. Note that $\cC \subset \cT_{1+\delta_{\pib}}-\cT_{1-\delta_{\pib}} \equiv \cT$, and so
\begin{align*}
	\text{Volume}(\cT) \ge \text{Volume}(\cC)= N \left(\frac{\delta_{\pib}}{\sqrt{K}}\right)^K.
\end{align*}
Note that here we use the notation $\pi$ for the Archimedes' constant, which differs from $\pib = \left(\pi_k\right)_{k\in[K]}$ for the mixing proportion of the GLoME model. Then, we define $\cV_K\left(a\right) = a^K\pi^{K/2}$ as the volume of a sphere of radius $a$. Since $z! \ge z^z e^{-z}$ and $\left(1+\delta_{\pib}\right)^K-\left(1-\delta_{\pib}\right)^K = K \int_{1-\delta_{\pib}}^{1+\delta_{\pib}}z^{K-1}dz \le 2 \delta_{\pib}K\left(1+\delta_{\pib}\right)^{K-1}$, it follows that
\begin{align*}
	&\cN\left(\delta_{\pib},\Pib_{K-1},\left\|\cdot\right\|_\infty\right) \nn\\
	&\le N \le \frac{\text{Volume}\left(\cC\right)}{\left(\delta_{\pib}/\sqrt{K}\right)^K} = \frac{1}{2^K} \frac{\cV_K\left(1+\delta_{\pib}\right)-\cV_K\left(1-\delta_{\pib}\right)}{\left(\delta_{\pib}/\sqrt{K}\right)^K}\\
	&=\frac{1}{2^K} \frac{\left[\left(1+\delta_{\pib}\right)^K-\left(1-\delta_{\pib}\right)^K\right]}{\left(\delta_{\pib}/\sqrt{K}\right)^K}\frac{\pi^{K/2}}{\left(K/2\right)!} \le \left(\frac{\pi e}{2}\right)^{K/2} 
	\frac{\left[\left(1+\delta_{\pib}\right)^K-\left(1-\delta_{\pib}\right)^K\right]}{\delta_{\pib}^K}\\
	& \le \frac{K \left(2 \pi e\right)^{K/2}}{\delta^{K-1}_{\pib}}.
\end{align*}

\subsubsection{Proof of Lemma \ref{lem_covering_number_WK_intro}} \label{sec_proof_lem_covering_number_WK_intro}
In order to find an upper bound for a covering number of $\cW_K$, we wish to construct a finite $\delta$-covering $\bfW_{K,\omegab}$ of $\cW_K$, with respect to the distance $\dsup$. That is, given any $\delta > 0,\bfw\left(\cdot;\omegab\right) \in \cW_K$, we aim to prove that there exists $\bfw\left(\cdot;\widehat{\omegab}\right) \in \bfW_{K,\omegab}$ such that
\begin{align}
	\dsup\left(\bfw(\cdot;\omegab),\bfw(\cdot;\widehat{\omegab})\right) =\max_{k\in[K]} \sup_{\yv \in \cY}\left|\cW_K(\yv;\omegab)-\cW_K\left(\yv;\widehat{\omegab}\right)\right| \le \delta. \label{eq_defineCoveringWeights_intro}
\end{align}
In order to accomplish such task, given any positive constants $\delta_{\cb},\delta_{\Gammab},\delta_{\pib}$, and any $k \in [K]$, let us define
\begin{align}
	\cF &= \left\{\cY \ni \yv \mapsto \ln\left(\phi_L\left(\yv;\cb,\Gammab\right)\right):\left\|\cb\right\|_\infty  \le A_{\cb},a_{\Gammab} \le m\left({\Gammab}\right) \le  M\left({\Gammab}\right) \le A_{\Gammab}\right\},\nn\\
	\cF_{\bfc_k} &= \Big\{\ln\left(\phi_L\left(\cdot;\bfc_k,\Gammab_k\right)\right): \ln\left(\phi_L\left(\cdot;\bfc_k,\Gammab_k\right)\right) \in \cF,\nn\\
	&\hspace{1cm} \bfc_{k,j} \in \left\{-C_\cY + l\delta_{\cb}/L: l = 0,\ldots,\lceil 2C_\cY L/\delta_{\cb} \rceil \right\},j \in [L] \Big\},\label{eq_defineNet_mu}\\
	\cF_{\bfc_k,\Gammab_k} &= \Bigg\{\ln\left(\phi_L\left(\cdot;\bfc_k,\Gammab_k\right)\right): \ln\left(\phi_L\left(\cdot;\bfc_k,\Gammab_k\right)\right) \in \cF_{\bfc_k}; \left[\vect\left(\Gammab_k\right)\right]_{i,j} =
	\gamma_{i,j}\frac{\delta_{\Gammab}}{L^2},\nn\\
	&\hspace{1cm}\gamma_{i,j} = \gamma_{j,i} \in \Z \cap \left[-\Big\lfloor \frac{L^2A_\Gammab}{\delta_{\Gammab}} \Big\rfloor,\Big\lfloor \frac{L^2A_\Gammab }{\delta_{\Gammab}} \Big\rfloor\right],i\in [L],j \in [L] \Bigg\},\label{eq_defineNet-sigma_intro}\\
	\bfW_{K,\omegab} &= \bigg\{\bfw\left(\cdot;\omegab\right): \bfw\left(\cdot;\omegab\right) \in \cW_K, \nn\\
	&\hspace{1cm} \forall k\in[K], \ln\left(\phi_L\left(\cdot;\bfc_k,\Gammab_k\right)\right)\in \cF_{\bfc_k,\Gammab_k},\pib\in \Pib_{K-1,\omegab}\bigg\}. \label{eq.eq.defineNet-alpha-Gen.fullnew_intro}
\end{align}
Here, $\lceil\cdot\rceil$ and $\lfloor\cdot\rfloor$ are ceiling and floor functions, respectively, and $\vect\left(\cdot\right)$ is an operator that stacks matrix columns into a column vector. In particular, we denote $\Pib_{K-1,\omegab}$ as a $\delta_{\pib}$-covering of $\Pib_{K-1}$, which is defined in  \cref{lem_MetricEntropySimplex_maximum_norm_intro}. By the previous definition, it holds that $\forall k \in [K]$, $\cF_{\bfc_k,\Gammab_k} \subset \cF_{\bfc_k} \subset \cF$, and $\bfW_{K,\omegab} \subset \cW_K$.

Next, we claim that $\bfW_{K,\omegab}$ is a finite $\delta$-covering of $\cW_K$ with respect to the distance $\dsup$.
To do this, for any $\bfw\left(\cdot;\omegab\right) = \left( \ln\left(\pi_k \phi_L\left(\cdot;\bfc_k,\Gammab_k\right)\right)\right)_{k\in[K]} \in \cW_K$, $\ln\left(\phi_L\left(\cdot;\bfc_k,\Gammab_k\right)\right)\in \cF,\pib\in\Pib_{K-1}$, and for any $k\in[K]$, by \eqref{eq_defineNet_mu}, we first choose a function $\ln\left(\phi_L\left(\cdot;\widehat{\cb}_k,\Gammab_k\right)\right) \in \cF_{\bfc_k}$ so that $$\left\|\widehat{\cb}_k-\bfc_k\right\|_1 = \sum_{j=1}^L \left|\widehat{\cb}_{k,j}-\bfc_{k,j}\right|  \le L\frac{\delta_{\cb}}{L}=\delta_{\cb}.$$ 
Furthermore, by \eqref{eq_defineNet-sigma_intro}, we can obtain a result
to construct the covariance matrix lattice. That is, any $\ln\left(\phi_L\left(\cdot;\widehat{\cb}_k,\Gammab_k\right)\right) \in \cF_{\bfc_k}$ can be approximated by $\ln\left(\phi_L\left(\cdot;\widehat{\cb}_k,\widehat{\Gammab}_k\right)\right)\in \cF_{\bfc_k,\Gammab_k}$ such that 
\begin{align}
	\left\|\vect \left(\widehat{\Gammab}_k\right)-\vect\left(\Gammab_k\right)\right\|_1&\equiv\left\|\vect \left(\widehat{\Gammab}_k-\Gammab_k\right)\right\|_1\nn\\
	&=\sum_{i=1}^L \sum_{j=1}^L\left|\left[\vect\left(\widehat{\Gammab}_k-\Gammab_{k}\right)\right]_{i,j}\right| \le \frac{L^2\delta_{\Gammab}}{L^2} = \delta_{\Gammab}.
\end{align}

Note that since for any $k\in[K]$, $\left(\yv,\bfc_k,\vect\left(\Gammab_k\right)\right) \mapsto\ln\left(\phi_L\left(\yv;\bfc_k,\Gammab_k\right)\right)$ is differentiable, it is also continuous \wrt~$\yv$ and its parameters $\bfc_k$ and $\Gammab_k$. Thus, for every fixed $\yv\in \cY$, for every $\widehat{\cb}_k,\bfc_k \in \cX$ with $\widehat{\cb}_k \le \bfc_k$, and for every $\widehat{\Gammab}_k$, $\Gammab_k$, where $\vect \left(\widehat{\Gammab}_k\right) \le \vect\left(\Gammab_k\right)$, we can apply the mean value theorem (see \citep[Lemma 2.5.1]{duistermaat2004multidimensional}) to $\ln\left(\phi_L\left(\yv;\cdot,\Gammab_k\right)\right)$ and $\ln\left(\phi_L\left(\yv;\widehat{\cb}_k,\cdot\right)\right)$ on the intervals $\left[\widehat{\cb}_k,\bfc_k\right]$ and $\left[\vect\left(\widehat{\Gammab}_k\right),\vect\left(\Gammab_k\right)\right]$ for some $z_{\bfc_k} \in \left(\widehat{\cb}_k,\bfc_k\right)$ and $z_{\Gammab_k} \in \left(\vect\left(\widehat{\Gammab}_k\right),\vect\left(\Gammab_k\right)\right)$, respectively, to get
\begin{align*}
	\ln\left(\phi_L\left(\yv;\widehat{\cb}_k,\Gammab_k\right)\right)-\ln\left(\phi_L\left(\yv;\bfc_k,\Gammab_k\right)\right) &= \left(\widehat{\cb}_k-\bfc_k\right)^\top \nabla_{\bfc_k} \ln\left(\phi_L\left(\yv;z_{\bfc_k},\Gammab_k\right)\right) ,
\end{align*}
\begin{align*}
	&\ln\left(\phi_L\left(\yv;\widehat{\cb}_k,\widehat{\Gammab}_k\right)\right)-\ln\left(\phi_L\left(\yv;\widehat{\cb}_k,\Gammab_k\right)\right)\nn\\
	&= \left(\vect\left(\widehat{\Gammab}_k\right)-\vect\left(\Gammab_k\right)\right)^\top \nabla_{\vect \left(\Gammab_k\right)} \ln\left(\phi_L\left(\yv;\widehat{\cb}_k,z_{\Gammab_k}\right)\right).
\end{align*}
Moreover, $\left(\yv,\bfc_k,\vect\left(\Gammab_k\right)\right) \mapsto \nabla_{\bfc_k} \ln\left(\phi_L\left(\yv;z_{\bfc_k},\Gammab_k\right)\right)$ and $\left(\yv,\bfc_k,\vect\left(\Gammab_k\right)\right) \mapsto \nabla_{\vect \left(\Gammab_k\right)} \ln\left(\phi_L\left(\yv;\widehat{\cb}_k,z_{\Gammab_k}\right)\right)$ are continuous functions on the compact set $\cU:=\cY \times \cY \times \left[a_\Gammab,A_\Gammab\right]^{L^2}$ leads to they attain minimum and maximum values (see \citep[Theorem 1.8.8]{duistermaat2004multidimensional}).
That is, we can set
\begin{align*}
	\zero &< \left(C_{\cb}\right)_{1,\ldots,L}^\top:=\max_{k\in[K]}\sup_{ \left(\yv,\bfc_k,\vect\left(\Gammab_k\right)\right)\in \cU} \left|\nabla_{\bfc_k} \ln\left|\phi_L\left(\yv;\bfc_k,\Gammab_k\right)\right|\right| < \infty,\\
	\zero &< \left(C_{\Gammab}\right)_{1,\ldots,L^2}^\top:=\max_{k\in[K]}\sup_{\left(\yv,\bfc_k,\vect\left(\Gammab_k\right)\right) \in \cU} \left|\nabla_{\vect\left(\Gammab_k\right)}\ln\left| \phi_L\left(\yv;\bfc_k,\Gammab_k\right)(x)\right|\right| < \infty.
\end{align*}
Therefore, by the Cauchy–Schwarz inequality, we have
\begin{align*}
	\sup_{\yv \in \cY} \left|\ln\left(\phi_L\left(\yv;\widehat{\cb}_k,\Gammab_k\right)\right)-\ln\left(\phi_L\left(\yv;\bfc_k,\Gammab_k\right)\right) \right| &\le \left|\widehat{\cb}_k-\bfc_k\right|^\top \left(C_{\cb}\right)_{1,\ldots,L}^\top \nn\\
	& = C_{\cb}  \left\|\widehat{\cb}_k-\bfc_k\right\|_1 \le C_{\cb} \delta_{\pib},\\
	\sup_{\yv \in \cY}\left|\ln\left(\phi_L\left(\yv;\widehat{\cb}_k,\widehat{\Gammab}_k\right)\right)-\ln\left(\phi_L\left(\yv;\widehat{\cb}_k,\Gammab_k\right)\right)\right| &
	\le C_{\Gammab} \left\|\vect \left(\widehat{\Gammab}_k-\Gammab_k\right)\right\|_1 \le C_{\Gammab} \delta_{\Gammab},
\end{align*}
and by using the triangle inequality, it follows that
\begin{align}
	&\max_{k \in [K]}\sup_{\yv \in \cY}\left|\ln\left(\phi_L\left(\yv;\widehat{\cb}_k,\widehat{\Gammab}_k\right)\right)-\ln\left(\phi_L\left(\yv;\bfc_k,\Gammab_k\right)\right) \right|
	\le  C_{\cb} \delta_{\cb}+C_{\Gammab} \delta_{\Gammab}.\label{eq.upperBoundLogGaussian_intro}
\end{align}
Moreover, for every $\pib \in \Pib_{K-1}$, \cref{lem_MetricEntropySimplex_maximum_norm_intro} implies that we can choose $\widehat{\pib}\in \Pib_{K-1,\omegab}$ so that $\max_{k\in[K]}\left|\pi_k - \widehat{\pib}_k\right| \le \delta_{\pib}$. Notice that $[a_\pib,\infty) \ni t \mapsto \ln(t), a_\pib > 0$ is a Lipschitz continuous function on $[a_\pib,\infty]$. Indeed, by the mean value theorem, it holds that there exists $c \in \left(t_1,t_2\right)$, such that 
\begin{align}
	\left|\ln\left(t_1\right)-\ln\left(t_2\right)\right|=\ln'(c)\left|t_1-t_2\right| \le\frac{1}{a_\pib}\left|t_1-t_2\right| \text{, for all } t_1, t_2 \in [a_\pib,\infty).\label{eq_lipschitzLog}
\end{align}	
Therefore, \eqref{eq_defineCoveringWeights_intro} can be obtained by the following evaluation
\begin{align*}
	&\dsup\left(\bfw(\cdot;\omegab),\bfw(\cdot;\widehat{\omegab})\right)\nn\\
	&= \max_{k\in[K]}\sup_{\yv \in \cY} \left|\ln \left(\pi_k \phi_L\left(\yv;\bfc_k,\Gammab_k\right)\right)-\ln\left(\widehat{\pib}_k \phi_L\left(\yv;\widehat{\cb}_k,\widehat{\Gammab}_k\right)\right)\right| \nn\\ 
	&=  \max_{k\in[K]} \sup_{\yv \in \cY} \left|\ln \left(\pi_k\right)-\ln \left(\widehat{\pib}_k\right)+\ln \left(\phi_L\left(\yv;\bfc_k,\Gammab_k\right)\right) - \ln\left(\phi_L\left(\yv;\widehat{\cb}_k,\widehat{\Gammab}_k\right)\right)\right| \nn\\ 
	&\le \max_{k\in[K]} \left|\ln \left(\pi_k\right)-\ln \left(\widehat{\pib}_k\right)\right|+  \max_{k\in[K]}\sup_{\yv \in \cY} \left| \ln \left(\phi_L\left(\yv;\bfc_k,\Gammab_k\right)\right) - \ln\left(\phi_L\left(\yv;\widehat{\cb}_k,\widehat{\Gammab}_k\right)\right)\right|
	\nn\\ 
	&\le \frac{1}{a_\pib} \max_{k\in[K]}\left|\pi_k-\widehat{\pib}_k\right| + C_{\cb} \delta_{\cb}+C_{\Gammab} \delta_{\Gammab} \left(\text{using \eqref{eq_lipschitzLog} and \eqref{eq.upperBoundLogGaussian_intro}}\right)
	\nn\\ 
	&\le \frac{\delta_{\pib}}{a_\pib} + C_{\cb} \delta_{\cb}+C_{\Gammab} \delta_{\Gammab} \left(\text{using \cref{lem_MetricEntropySimplex_maximum_norm_intro}}\right)
	\le  \frac{\delta}{3} + \frac{\delta}{3}+ \frac{\delta}{3}
	= \delta,
\end{align*}
where we choose $\delta_{\pib} =\frac{\delta a_\pib}{3},\delta_{\cb} =\frac{\delta}{3 C_\mu}, \delta_{\Gammab}=\frac{\delta}{3 C_\Gammab}$. Finally, we get the covering number
\begin{align*}
	&\cN\left(\delta,\cW_K,\dsup \right)\nn\\
	&\le \card \left(\bfW_{K,\omegab}\right) = \left[\frac{2C_\cY L}{\delta_{\cb}}\right]^{KL} \left[\frac{2A_\Gammab L^2}{\delta_{\Gammab}}\right]^{\frac{L\left(L+1\right)}{2}K} \cN\left(\delta_{\pib},\Pib_{K-1},\left\|\cdot\right\|_\infty\right)\nn\\ 
	&=\frac{C}{\delta^{KL+K\frac{L(L+1)}{2}+K-1}}=\frac{C}{\delta^{\dim\left(\cW_K\right)}},
\end{align*}
where $$C = \left(6C_\cb C_\cY L\right)^{KL}\left(6C_\Gammab A_\Gammab L^2\right)^{\frac{L(L+1)}{2}K} \left(\frac{3}{a_\pib}\right)^{K-1} K \left(2 \pi e\right)^{K/2}.$$

\subsection{Proof of Lemma \ref{lem_Bracketing_Entropy_Experts}: Bracketing entropy of Gaussian experts with 1-block covariance matrices} \label{sec_proof_lem_Bracketing_Entropy_Experts}
Note that the proof of \cref{lem_Bracketing_Entropy_Experts} is adapted from \cite{montuelle2014mixture} and is given here for the sake of completeness.

Indeed, Lemmas 1 and 2 from \cite{montuelle2014mixture} imply that there exists a constant $C_{\UpsilondK} = \ln \left(\sqrt{2}+ \sqrt{D} d T_\Upsilonb\right)$ (when $\UpsilondK = \Upsilonb_{b,d}^K$) or $C_{\UpsilondK} = \ln \left(\sqrt{2}+ \sqrt{D} \binom{d+L}{L} T_\Upsilonb\right)$ (when $\UpsilondK = \Upsilonb_{p,d}^K$) such that, $\forall \delta \in \left(0,\sqrt{2}\right)$,
\begin{align}\label{eq_assumption_DIM_means}
	\entropy(\delta,\UpsilondK) \le \dim\left(\UpsilondK\right) \left(C_{\UpsilondK} + \ln \frac{1}{\delta}\right).
\end{align}

Next, we rely on \cref{prop_montuelle_Prop2} for constructing of Gaussian brackets to for the Gaussian experts.
\begin{prop}[Proposition 2 from \cite{montuelle2014mixture}] \label{prop_montuelle_Prop2}
	Let $\kappa \geq \frac{17}{29}$ and $\gamma_{\kappa} = \frac{25(\kappa - \frac{1}{2})}{49\left(1+\frac{2\kappa}{5}\right)}$. For any $0 < \delta \le \sqrt{2}$, any $D \geq 1$, and any  $\delta_\Sigma \le \frac{1}{5 \sqrt{\kappa^2 \cosh\left(\frac{2\kappa}{5}\right)+\frac{1}{2}}} \frac{\delta}{D}$, let $\left(\upsilonb_d,B,\bfA,\bfP\right) \in \UpsilondK \times \left[B_-,B_+\right]\times \cA\left(\lambda_-,\lambda_+\right)\times SO(D)$ and
	$\left(\widetilde{\upsilonb}_d,\widetilde{B},\widetilde{\bfA},\widetilde{\bfP}\right) \in \UpsilondK \times \left[B_-,B_+\right]\times \cA\left(\lambda_-,+\infty\right)\times SO(D)$, and define $\Sigmab = B\bfP\bfA\bfP^\top$, $\widetilde{\Sigmab} = \widetilde{B}\widetilde{\bfP}\widetilde{\bfA}\widetilde{\bfP}^\top$,
	\begin{align*}
		t^-\left(\xv,\yv\right) &= \left(1+\kappa \delta_\Sigmab\right)^{-D} \Phi_D\left(\xv;\widetilde{\upsilonb}_d(\yv),\left(1+\delta_\Sigmab\right)^{-1}\widetilde{\Sigmab}\right),\text{ and }\nn\\
		t^+\left(\xv,\yv\right) &= \left(1+\kappa \delta_\Sigmab\right)^{D} \Phi_D\left(\xv;\widetilde{\upsilonb}_d(\yv),\left(1+\delta_\Sigmab\right)\widetilde{\Sigmab}\right).
	\end{align*}
	If
	\begin{align*}
		\begin{cases}
			\forall \yv \in \cY, \norm{\upsilonb_d(\yv) - \widetilde{\upsilonb}_d(\yv)}^2 \le D \gamma_{\kappa} \cL_-\lambda_- \frac{\lambda_-}{\lambda_+}\delta^2_\Sigmab\\
			(1+\frac{2}{25}\delta_\Sigmab)^{-1}\widetilde{B} \le B \le \widetilde{B}\\
			\forall i \in [D], \left|\bfA^{-1}_{i,i} -\widetilde{\bfA}^{-1}_{i,i} \right| \le \frac{1}{10} \frac{\delta_\Sigmab}{\lambda_+}\\
			\forall \xv \in \R^D, \norm{\bfP \xv - \widetilde{\bfP}\xv}\le \frac{1}{10} \frac{\lambda_-}{\lambda_+}\delta_\Sigmab \norm{\xv}
		\end{cases}
	\end{align*}
	then $[t^-,t^+]$ is a $\frac{\delta}{5}$ Hellinger bracket such that $t^-(\xv,\yv) \le \Phi_D\left(\xv;\upsilonb_d(\yv),\Sigmab\right)\le t^+(\xv,\yv)$.
\end{prop}
Then, following the same argument as in \cite[Appendix B.2.3]{montuelle2014mixture}, \cref{prop_montuelle_Prop2} allows us to construct nets over the spaces of the means, the volumes, the eigenvector matrices, the normalized eigenvalue matrices and then control the bracketing entropy of $\gkd$. More precisely, three different contexts are considered for the mean, volume and matrix parameters. They can be all known $(\star = 0)$, unknown but common to all classes $(\star = c)$, unknown and possibly different for every class $(\star = K)$. For example, $[\vb_K,B_0,\bfP_0,\bfA_0]$ denotes a model in which only mean vectors are assumed to be free. Then, we obtain
\begin{align} \label{bracketingEntropyGaussianFamilies}
	\cH_{[\cdot],\sup_\yv\max_k d_\xv }\left(\frac{\delta}{5},\gkd\right) \le \dim\left(\gkd\right) \left(C_{\gkd} + \ln \left(\frac{1}{\delta}\right)\right),
\end{align}
where $\dim\left(\gkd\right) = Z_{\vb,\star}+Z_{B,\star}+ \frac{D(D-1)}{2}Z_{\bfP_,\star}+ (D-1)Z_{\bfA,\star}$. Here, $Z_{\vb,K} = \dim\left(\UpsilondK\right), Z_{\vb,c} = \dim\left(\Upsilonb_1\right),Z_{\vb,0} = 0$, $Z_{B,0} = Z_{\bfP_,0} = Z_{\bfA,0} = 0$, $Z_{B,c} = Z_{\bfP_,c} = Z_{\bfA,c} = 1$, $Z_{B,K} = Z_{\bfP_,K} = Z_{\bfA,K} = K$, and given a universal constant $c_U$,
\begin{align*}
	C_{\gkd}& = \ln\left(5D \sqrt{\kappa^2 \cosh \left(\frac{2\kappa}{5}\right)+ \frac{1}{2}}\right) + C_{\UpsilondK} + \frac{1}{2} \ln \left(\frac{\lambda_+}{D\gamma_kB_-\lambda^2_-}\right)\nn\\
	&\quad + \ln\left(\frac{4+ 129 \ln\left(\frac{B_+}{B_-}\right)}{10}\right) +  \frac{D(D-1)}{2}\ln\left(c_U\right) \nn\\
	& \quad + \ln\left(\frac{10\lambda_+}{\lambda_-}\right)  + \ln \left(\frac{4}{5}+ \frac{52\lambda_+}{5\lambda_-}\ln \left(\frac{\lambda_+}{\lambda_-}\right)\right).
\end{align*}


\subsection{Proof of Lemma \ref{lem_Inequality_Hellinger_Supx_dy_dk_intro}}
\label{sec_lem_Inequality_Hellinger_Supx_dy_dk_intro}
We first aim to prove that $\thel(s,t) \le \sup_\yv d^2_\xv(s,t)$. Indeed, by definition, it follows that
\begin{align*}
	\thel\left(s,t\right) &= \E{\Yv_{[n]}}{\frac{1}{n} \sum_{i=1}^n d^2_\xv\left(s\left(\cdot\mid \Yv_i\right),t\left(\cdot\mid \Yv_i\right)\right)}
	=\frac{1}{n} \sum_{i=1}^n \E{\Yv_{[n]}}{d^2_\xv\left(s\left(\cdot\mid \Yv_i\right),t\left(\cdot\mid \Yv_i\right)\right)}\nn\\
	&=\frac{1}{n} \sum_{i=1}^n \int_{\cY}d^2_\xv\left(s\left(\cdot\mid \yv\right),t\left(\cdot\mid \yv\right)\right)s_{\xv,0}(\yv)d\yv \nn\\
	&\le  \sup_\yv d^2_\xv~(s,t) \frac{1}{n} \sum_{i=1}^n\int_{\cY}s_{\xv,0}(\yv)d\yv = \sup_\yv d^2_\xv~(s,t),
\end{align*}
where $s_{\xv,0}$ denotes that marginal PDF of $s_0$, \wrt~$\xv$.
Consequently, it holds that $\thell(s,t) = \sqrt{\thel(s,t)} \le \sqrt{\sup_\yv d^2_\xv(s,t)}= \sup_\yv d_\xv(s,t)$.
To prove that  $$\cH_{[\cdot],\thell}\left(\delta,S_\bfm\right) \le \cH_{[\cdot],\sup_\yv d_\xv }\left(\delta,S_\bfm\right),$$ 
it is sufficient to check that $$\cN_{[\cdot],\thell}\left(\delta,S_\bfm\right) \le \cN_{[\cdot],\sup_\yv d_\xv }\left(\delta,S_\bfm\right).$$
By using the definition of bracketing entropy in \eqref{eq_define_BracketingEntropy_intro} and $\thell(s,t) \le \sup_\yv d_\xv(s,t)$, given
\begin{align*}
	A &= \left\{n \in \Ns: \exists t^-_1,t^+_1,\ldots,t^-_n,t^+_n \text{ s.t }\sup_\yv d_\xv(s,t)\left(t^-_k,t^+_k\right) \le \delta,S_\bfm \subset \bigcup_{k=1}^n \left[t^-_k,t^+_k\right] \right\},\\
	B &=\left\{n\in \Ns: \exists t^-_1,t^+_1,\ldots,t^-_n,t^+_n \text{ s.t }\thell\left(t^-_k,t^+_k\right) \le \delta,S_\bfm \subset \bigcup_{k=1}^n \left[t^-_k,t^+_k\right] \right\},
\end{align*}
it leads to that $A \subset B$ and then \eqref{eq_Inequality_Hellinger_Supx_dy_intro} follows, since
\begin{align*}
	\cN_{[\cdot],\sup_\yv d_\xv(s,t)}\left(\delta,S_\bfm\right)= \min A\ge \min B =\cN_{[\cdot],\thell}\left(\delta,S_\bfm\right).
\end{align*}

With the similar argument as in the proof of \eqref{eq_Inequality_Hellinger_Supx_dy_intro}, it holds that $d_{\cP_K}\left(g^+,g^-\right) \le \sup_\yv d_k(g^+,g^-)$ and \eqref{eq_Inequality_Hellinger_Supx_dk_intro} is proved.

\subsection{Proof of Lemma \ref{lem_bracketingEntropyGaussianBlock_intro}: Bracketing entropy of Gaussian experts with multi-block-diagonal covariance matrices}\label{sec_proof_lem_bracketingEntropyGaussianBlock_intro}
It is worth mentioning that without any structures on covariance matrices of Gaussian experts from the collection $\cM$, \cref{lem_bracketingEntropyGaussianBlock_intro} can be proved using \citep[Proposition 2 and Appendix B.2.3]{montuelle2014mixture}, for constructing of Gaussian brackets to deal with the Gaussian experts. However, dealing with block-diagonal covariance matrices with random subcollection is not a trivial extension. We have to establish more constructive bracketing entropies in the spirits of \cite{maugis2011non,devijver2015finite,devijver2018block}. 

Given any $k \in [K]$, by defining
\begin{align} \label{eq_define_GaussianExpert_1D_intro}
	\gdbk &= \big\{
	\left(\xv,\yv\right)\mapsto \phi\left(\xv;\upsilonb_{k,d}(\yv),\Sigmab_k\left(\bfB_k\right)\right)=:\phi_k: \nn\\
	& \hspace{4cm} \upsilonb_{k,d} \in \Upsilondk, \Sigmab_k\left(\bfB_k\right) \in \bfV_k(\bfB_k) \big\},
\end{align}
it follows that $\gkdb = \prod_{k=1}^K \gdbk$, where $\prod$ stands for the cartesian product. By using \cref{lem_HGK_HG_intro}, which is proved in \cref{sec.prooflem_HGK_HG_intro}, it follows that
\begin{align} \label{eq_HGK_HG_intro}
	\cH_{[\cdot],d_{\gkdb} }\left(\frac{\delta}{2},\gkdb\right) \le \sum_{k=1}^K \cH_{[\cdot],d_{\gdbk} }\left(\frac{\delta}{2\sqrt{K}},\gdbk\right).
\end{align}
\begin{lem} \label{lem_HGK_HG_intro}
	Given $\gkdb = \prod_{k=1}^K \gdbk$, where $\gdbk$ is defined in \eqref{eq_define_GaussianExpert_1D_intro}, it holds that
	\begin{align*}
		\cN_{[\cdot],d_{\gkdb}}\left(\frac{\delta}{2},\gkdb\right) \le \prod_{k=1}^K\cN_{[\cdot],d_{\gdbk}}\left(\frac{\delta}{2\sqrt{K}},\gdbk\right),
	\end{align*}
	where for any $\phi^+,\phi^- \in \gkdb$ and any $\phi_k^+,\phi_k^- \in \gdbk,k\in[K]$,
	\begin{align*}
		d^2_{\gkdb}\left(\phi^+,\phi^-\right) &= \E{\Yv_{[n]}}{\frac{1}{n}\sum_{i=1}^n \sum_{k=1}^K \hel\left(\phi^+_k \left(\cdot,\Yv_i\right),\phi^-_k\left(\cdot,\Yv_i\right)\right)},\nn\\
		d^2_{\gdbk}\left(\phi_k^+,\phi_k^-\right) &= \E{\Yv_{[n]}}{\frac{1}{n}\sum_{i=1}^n  \hel\left(\phi_k^+ \left(\cdot,\Yv_i\right),\phi_k^-\left(\cdot,\Yv_i\right)\right)}
		.
	\end{align*}
\end{lem}
\cref{lem_bracketingEntropyGaussianBlock_intro} is proved via \eqref{eq_HGK_HG_intro} and \cref{lem.bracketingEntropyGaussianBlock1-D_intro}, which is proved in \cref{sec.lem.bracketingEntropyGaussianBlock1-D_intro}.
\begin{lem}\label{lem.bracketingEntropyGaussianBlock1-D_intro}
	By defining $\gdbk$ as in \eqref{eq_define_GaussianExpert_1D_intro}, for all $\delta \in (0,\sqrt{2}]$, it holds that
	\begin{align}
		\cH_{[\cdot],d_{\gdbk} }\left(\frac{\delta}{2},\gdbk\right) &\le \dim\left(\gdbk\right) \left(C_{\gdbk} + \ln \left(\frac{1}{\delta}\right)\right),\label{eq_bracketingEntropyGaussianBlock1_D_intro} \text{ where }\\
		D_{\bfB_k} &=  \sum_{g=1}^{G_k} \frac{\card\left(d_k^{[g]}\right)\left(\card\left(d_k^{[g]}\right)-1\right)}{2},\nn\\
		C_{\gdbk} &= \frac{D_{\bfB_k} \ln \left(\frac{6 \sqrt{6} \lambda_M D^2\left(D-1\right)}{ \lambda_m D_{\bfB_k}}\right) + \dim\left(\Upsilondk\right) \ln\left(\frac{6 \sqrt{2D}\exp\left(C_{\Upsilondk}\right)}{ \sqrt{\lambda_m}}\right)}{\dim\left(\gdbk\right)}\nn.
	\end{align}
\end{lem}
Indeed, \eqref{eq_HGK_HG_intro} and \eqref{eq_bracketingEntropyGaussianBlock1_D_intro} lead to
\begin{align*}
	\cH_{[\cdot],d_{\gkdb} }\left(\frac{\delta}{2},\gkdb\right) &\le \sum_{k=1}^K \cH_{[\cdot],d_{\gdbk} }\left(\frac{\delta}{2\sqrt{K}},\gdbk\right)\nn\\
	&\le  \sum_{k=1}^K  \dim\left(\gdbk\right)\left(C_{\gdbk} + \ln \left(\sqrt{K}\right) + \ln \left(\frac{1}{\delta}\right)\right)\\
	& \le \dim\left(\gkdb\right) \left(C_{\gkdb} + \ln \left(\frac{1}{\delta}\right)\right).
\end{align*}
Here, $C_{\gkdb} = \sum_{k=1}^K C_{\gdbk} + \ln \left(\sqrt{K}\right)$ and note that $\dim\left(\gkdb\right) = \sum_{k=1}^K  \dim\left(\gdbk\right)$, $\dim\left(\gdbk\right) = D_{\bfB_k}+\dim\left(\Upsilondk\right)$, $\dim\left(\Upsilondk\right) = D d_{\Upsilondk}, C_{\Upsilondk}=\sqrt{D}d_{\Upsilondk}T_{\Upsilondk}$ (in cases where linear combination of bounded functions are used for means, \ie $\Upsilondk = \Upsilonb_b$) or $\dim\left(\Upsilondk\right) = D \binom{d_{\Upsilondk} + L}{ L}$, $C_{\Upsilondk} = \sqrt{D} \binom{d_{\Upsilondk} + L}{L} T_{\Upsilondk}$ (in cases where we use polynomial means, \ie $\Upsilondk = \Upsilonb_p$).

\subsubsection{Proof of Lemma \ref{lem_HGK_HG_intro}}\label{sec.prooflem_HGK_HG_intro}
By the definition of the bracketing entropy in \eqref{eq_define_BracketingEntropy_intro}, for each $k\in[K]$, let $\left\{\left[\phi^{l,-}_k,\phi^{l,+}_k\right]\right\}_{ 1 \le  l \le  \cN_{\gdbk}}$ be a minimal covering of $\delta_k$ brackets for $d_{\gdbk}$ of $\gdbk$, with cardinality $
\cN_{\gdbk}$. This leads to
\begin{align*}
	\forall l \in \left[\cN_{\gdbk}\right], d_{\gdbk}\left(\phi^{l,-}_k,\phi^{l,+}_k\right) \le  \delta_k.
\end{align*} 
Therefore, we claim that the set $\left\{\prod_{k=1}^K\left[\phi^{l,-}_k,\phi^{l,+}_k\right]\right\}_{ 1 \le  l \le  \cN_{\gdbk}}$ is a covering of $\frac{\delta}{2}$-bracket for $d_{\gkdb}$ of $\gkdb$ with cardinality $\prod_{k=1}^K \cN_{[\cdot], d_{\gdbk} }\left(\delta_k,\gdbk\right)$.
Indeed, let any $\phi = \left(\phi_k\right)_{k\in[K]} \in \gkdb$. Consequently, for each $k\in[K], \phi_k \in \gdbk$, there exists $l(k)\in\left[ \cN_{\gdbk}\right]$, such that
\begin{align*}
	\phi^{l(k),-}_k\le  \phi_k \le \phi^{l(k),+}_k,
	d^2_{\gdbk}\left(\phi^{l(k),+}_k,\phi^{l(k),-}_k\right) \le \left(\delta_k\right)^2.
\end{align*}
Then, it follows that $\phi \in \left[\phi^{-},\phi^{+}\right] \in \left\{\prod_{k=1}^K\left[\phi^{l,-}_k,\phi^{l,+}_k\right]\right\}_{ 1 \le  l \le  \cN_{\gdbk}}$, with $\phi^{-} = \left(\phi^{l(k),-}_k\right)_{k\in[K]},\phi^{+} = \left(\phi^{l(k),+}_k\right)_{k\in[K]}$, which implies that  $\left\{\prod_{k=1}^K\left[\phi^{l,-}_k,\phi^{l,+}_k\right]\right\}_{ 1 \le  l \le  \cN_{\gdbk}}$ is a bracket covering of $\gkdb$.

Now, we want to verify that the size of this bracket is $\delta/2$ by choosing $\delta_k = \frac{\delta}{2\sqrt{K}}, \forall k \in [K]$. It follows that 
\begin{align*}
	d^2_{\gkdb}\left(\phi^{-},\phi^{+}\right) &= \E{\Yv_{[n]}}{\frac{1}{n}\sum_{i=1}^n \sum_{k=1}^K \hel\left(\phi^{l(k),-}_k \left(\cdot,\Yv_i\right),\phi^{l(k),+}_k\left(\cdot,\Yv_i\right)\right)}\nn\\
	&= \sum_{k=1}^K\E{\Yv_{[n]}}{\frac{1}{n}\sum_{i=1}^n  \hel\left(\phi^{l(k),-}_k \left(\cdot,\Yv_i\right),\phi^{l(k),+}_k\left(\cdot,\Yv_i\right)\right)} \nn\\
	&= \sum_{k=1}^K d^2_{\gdbk}\left(\phi^{l(k),-}_k ,\phi^{l(k),+}_k\right)
	\le K \left(\frac{\delta}{2\sqrt{K}}\right)^2 = \left(\frac{\delta}{2}\right)^2.
\end{align*}
To this end, by definition of a minimal $\frac{\delta}{2}$-bracket covering number for $\gkdb$, \cref{lem_HGK_HG_intro} is proved.
%

\subsubsection{Proof of Lemma \ref{lem.bracketingEntropyGaussianBlock1-D_intro}}\label{sec.lem.bracketingEntropyGaussianBlock1-D_intro}
%

To provide the upper bound of
the bracketing entropy in  \eqref{eq_bracketingEntropyGaussianBlock1_D_intro}, our technique is adapted from the work of \cite{genovese2000rates} for unidimensional Gaussian mixture families, which is recently generalized to multidimensional case by \cite{maugis2011non} for Gaussian mixture models. Furthermore, we make use of the results from \cite{devijver2018block} to deal with block-diagonal covariance matrices, $\bfV_k\left(\bfB_k\right),k\in[K]$, and from \cite{montuelle2014mixture} to handle the means of Gaussian experts $\Upsilondk,k\in[K]$.

The main idea in our approach is to define firstly a net over the parameter spaces of Gaussian experts, $\Upsilondk\times \bfV_k\left(\bfB_k\right), k\in[K]$, and to construct a bracket covering of $\gdbk$ according to the tensorized Hellinger distance. 
Note that $\dim\left(\gdbk\right)=\dim\left(\Upsilondk\right)+\dim\left(\bfV_k\left(\bfB_k\right)\right)$. 

\paragraph*{Step 1: Construction of a net for the block-diagonal covariance matrices.}
Firstly, for $k\in[K]$, we denote by $\adj\left(\Sigmab_k\left(\bfB_k\right)\right)$ the adjacency matrix associated to the covariance matrix $\Sigmab_k\left(\bfB_k\right)$. Note that this matrix of size $D^2$ can be defined by a vector of concatenated upper triangular vectors. We are going to make use of the result from \cite{devijver2018block} to handle the block-diagonal covariance matrices $\Sigmab_k\left(\bfB_k\right)$, via its corresponding adjacency matrix. To do this, we need to construct a discrete space for $\left\{0,1\right\}^{D(D-1)/2}$, which is a one-to-one correspondence (bijection) with
$$\cA_{\bfB_k} = \left\{\bfA_{\bfB_k} \in \cS_D\left(\left\{0,1\right\}\right):\exists \Sigmab_k\left(\bfB_k\right) \in \bfV_k\left(\bfB_k\right) \text{ s.t } \adj\left(\Sigmab_k\left(\bfB_k\right)\right) = \bfA_{\bfB_k} \right\},$$
where $\cS_D\left(\left\{0,1\right\}\right)$ is the set of symmetric matrices of size $D$ taking values on $\left\{0,1\right\}$. 

Then, we want to deduce a discretization of the set of covariance matrices. Let $h$ denotes Hamming distance on $\left\{0,1\right\}^{D(D-1)/2}$ defined by
$$d(z,z') = \sum_{i=1}^n \Indi\left\{z \neq z'\right\}, \text{ for all } z,z' \in \left\{0,1\right\}^{D(D-1)/2}.$$
Let $\left\{0,1\right\}_{\bfB_k}^{D(D-1)/2}$ be the subset of $\left\{0,1\right\}^{D(D-1)/2}$ of vectors for which the corresponding graph has structure $\bfB_k = \left(d^{[g]}_k\right)_{g\in\left[G_k\right]}$. Corollary 1 and Proposition 2 from Supplementary Material A of \cite{devijver2018block} imply that there exists some subset $\cR$ of $\left\{0,1\right\}^{D(D-1)/2}$, as well as its equivalent $\cA^{\disc}_{\bfB_k}$ for adjacency matrices such that, given $\epsilon > 0$, and
\begin{align*}
	{\tilde{S}}^{\disc}_{\bfB_k}(\epsilon) &= \bigg\{\Sigmab_k\left(\bfB_k\right) \in \cS_D^{++}(\mathbb{R}): \adj\left(\Sigmab_k\left(\bfB_k\right)\right) \in \cA^{\disc}_{\bfB_k},\nn\\
	& \hspace{2cm}\left[\Sigmab_k\left(\bfB_k\right)\right]_{i,j} = \sigma_{i,j}\epsilon,\sigma_{i,j} \in \left[\frac{-\lambda_M}{\epsilon},\frac{\lambda_M}{\epsilon}\right]\bigcap \Z \bigg\},
\end{align*} 
it holds that for all $\left(\Sigmab_k\left(\bfB_k\right),\widetilde{\Sigmab}_k\left(\bfB_k\right) \right)\in \left({\tilde{S}}^{\disc}_{\bfB_k}(\epsilon)\right)^2 \text{\st} \Sigmab_k\left(\bfB_k\right) \neq \widetilde{\Sigmab}_k\left(\bfB_k\right)$, 
\begin{align}
	\left\|\Sigmab_k\left(\bfB_k\right)-\widetilde{\Sigmab}_k\left(\bfB_k\right)\right\|_2^2 &\le \frac{D_{\bfB_k}}{2} \wedge \epsilon^2,\nn\\
	\card\left({\tilde{S}}^{\disc}_{\bfB_k}(\epsilon)\right) & \le  \left(\Bigg\lfloor \frac{2\lambda_M}{\epsilon}\Bigg\rfloor\frac{D\left(D-1\right)}{2 D_{\bfB_k}}\right)^{D_{\bfB_k}},\\
	D_{\bfB_k} =\dim\left(\bfV_k\left(\bfB_k\right)\right)&= \sum_{g=1}^{G_k} \frac{\card\left(d_k^{[g]}\right)\left(\card\left(d_k^{[g]}\right)-1\right)}{2}. \label{eq.upperCardCovariance_intro}
\end{align}
By choosing $\epsilon^2 \le \frac{D_{\bfB_k}}{2}$, given $\Sigmab_k\left(\bfB_k\right) \in \bfV_k\left(\bfB_k\right)$, then there exists $\widetilde{\Sigmab}_k\left(\bfB_k\right) \in {\tilde{S}}^{\disc}_{\bfB_k}(\epsilon)$, such that 
\begin{align}
	\left\|\Sigmab_k\left(\bfB_k\right)-\widetilde{\Sigmab}_k\left(\bfB_k\right)\right\|_2^2 \le \epsilon^2. \label{eq.netsCovarianceG1_intro}
\end{align}

\paragraph*{Step 2: Construction of a net for the mean functions.}
Based on $\widetilde{\Sigmab}_k\left(\bfB_k\right)$, we can construct the following bracket covering of $\gdbk$ by defining the nets for the means of Gaussian experts.
The proof of Lemma 1, page 1693, from \cite{montuelle2014mixture} implies that
\begin{align*}
	\cN_{\left[\cdot\right],\sup_\yv \left\|\cdot\right\|_2}\left(\delta_{\Upsilondk},\Upsilondk\right) \le \left(\frac{\exp\left(C_{\Upsilondk}\right)}{\delta_{\Upsilondk}}\right)^{\dim\left(\Upsilondk\right)}.
\end{align*} 
Here $\dim\left(\Upsilondk\right) = D d_{\Upsilondk}$, and $C_{\Upsilondk}=\sqrt{D}d_{\Upsilondk}T_{\Upsilondk}$ in the general case or
$\dim\left(\Upsilondk\right) = D \binom{d_{\Upsilondk} + L}{ L}$, and $C_{\Upsilondk} = \sqrt{D} \binom{d_{\Upsilondk} + L}{L} T_{\Upsilondk}$ in the special case of polynomial means. Then, by the definition of bracketing entropy in \eqref{eq_define_BracketingEntropy_intro}, for any minimal $\delta_{\Upsilondk}$-bracketing covering of the means from Gaussian experts, denoted by $G_{\Upsilondk}\left(\delta_{\Upsilondk}\right)$, it is true that
\begin{align}
	\card\left(G_{\Upsilondk}\left(\delta_{\Upsilondk}\right)\right) \le \left(\frac{\exp\left(C_{\Upsilondk}\right)}{\delta_{\Upsilondk}}\right)^{\dim\left(\Upsilondk\right)}.\label{eq.cardMeanGaussianEx_intro}
\end{align}
Therefore, given $\alpha >0$, which is specified later, we claim that the set
\begin{align*}
	\left\{\left[l,u\right]\left| 
	\begin{array}{l} 
		l(\xv,\yv) = \left(1+2\alpha\right)^{-D} \phi\left(\xv;\widetilde{\upsilonb}_{k,d}(\yv),\left(1+\alpha\right)^{-1}\widetilde{\Sigmab}_k\left(\bfB_k\right)\right), \\
		u(\xv,\yv) = \left(1+2\alpha\right)^{D} \phi\left(\xv;\widetilde{\upsilonb}_{k,d}(\yv),\left(1+\alpha\right)\widetilde{\Sigmab}_k\left(\bfB_k\right)\right), \\
		\widetilde{\upsilonb}_{k,d} \in G_{\Upsilondk}\left(\delta_{\Upsilondk}\right),\widetilde{\Sigmab}_k\left(\bfB_k\right) \in {\tilde{S}}^{\disc}_{\bfB_k}(\epsilon)
	\end{array}
	\right.\right\},
\end{align*}
is a $\delta_{\Upsilondk}$-brackets set over $\gdbk$. Indeed, let $\cX\times\cY \ni  (\xv,\yv) \mapsto f(\xv,\yv) = \phi\left(\xv;\upsilonb_{k,d}(\yv),\Sigmab_k\left(\bfB_k\right)\right)$ be a function of $\gdbk$, where $\upsilonb_{k,d} \in \Upsilondk$ and $\Sigmab_k\left(\bfB_k\right) \in \bfV_k\left(\bfB_k\right)$. According to \eqref{eq.netsCovarianceG1_intro}, there exists $\widetilde{\Sigmab}_k\left(\bfB_k\right) \in {\tilde{S}}^{\disc}_{\bfB_k}(\epsilon)$, such that 
\begin{align*}
	\left\|\Sigmab_k\left(\bfB_k\right)-\widetilde{\Sigmab}_k\left(\bfB_k\right)\right\|_2^2 \le \epsilon^2.
\end{align*}
By definition of $G_{\Upsilondk}\left(\delta_{\Upsilondk}\right)$, there exists $\widetilde{\upsilonb}_{k,d} \in G_{\Upsilondk}\left(\delta_{\Upsilondk}\right)$, such that 
\begin{align}
	\sup_{\yv \in \cY} \left\|\widetilde{\upsilonb}_{k,d}(\yv)-\upsilonb_{k,d}(\yv)\right\|_2^2 \le \delta^2_{\Upsilondk}. \label{eq.netMeanGaussian_intro}
\end{align}

\paragraph*{Step 3: Upper bound of the number of the bracketing entropy.}
Next, we wish to make use of \cref{lem.ratioGaussian_intro} to evaluate the ratio of two Gaussian densities.
\begin{lem}[Proposition C.1 from \cite{maugis2011non}]\label{lem.ratioGaussian_intro}
	Let $\phi\left(\cdot;\mub_1,\Sigmab_1\right)$ and $\phi\left(\cdot;\mub_2,\Sigmab_{2}\right)$ be two Gaussian densities. If $\Sigmab_{2}-\Sigmab_1$ is a positive definite matrix then for all $\xv \in \R^D$, 
	\begin{align*}
		\frac{\phi\left(\xv;\mub_1,\Sigmab_1\right)}{\phi\left(\xv;\mub_2,\Sigmab_{2}\right)} \le \sqrt{\frac{\left|\Sigmab_{2}\right|}{\left|\Sigmab_1\right|}} \exp\left[\frac{1}{2}\left(\mub_1-\mub_2\right)^\top \left(\Sigmab_{2}-\Sigmab_1\right)^{-1}\left(\mub_1-\mub_2\right)\right].
	\end{align*}
\end{lem}
The following \cref{lem.assumptionPDM_intro} allows us to fulfill the assumptions of \cref{lem.ratioGaussian_intro}.
\begin{lem}
	\label{lem.assumptionPDM_intro}
	Assume that $0 < \epsilon < \lambda^2_m/9$, and set $\alpha = 3 \sqrt{\epsilon} /\lambda_m$.  Then, for every $k\in[K]$, $\left(1+\alpha\right)\widetilde{\Sigmab}_k\left(\bfB_k\right) - \Sigmab_k\left(\bfB_k\right)$ and $\Sigmab_k\left(\bfB_k\right) - \left(1+\alpha\right)^{-1}\widetilde{\Sigmab}_k\left(\bfB_k\right)$ are both positive definite matrices. Moreover, for all $\xv \in \R^D$,
	\begin{align*}
		\xv^\top \left[\left(1+\alpha\right)\widetilde{\Sigmab}_k\left(\bfB_k\right) - \Sigmab_k\left(\bfB_k\right)\right]\xv &\ge \epsilon \left\|\xv\right\|_2^2,\quad
		\nn\\
		\xv^\top \left[\Sigmab_k\left(\bfB_k\right) - \left(1+\alpha\right)^{-1}\widetilde{\Sigmab}_k\left(\bfB_k\right)\right]\xv &\ge \epsilon \left\|\xv\right\|_2^2.
	\end{align*} 	
\end{lem}
\begin{proof}[Proof of \cref{lem.assumptionPDM_intro}]
	For all $\xv \neq \zero$, since $$\sup_{\lambda \in \text{vp}\left(\Sigmab_k\left(\bfB_k\right)-\widetilde{\Sigmab}_k\left(\bfB_k\right)\right)}\left|\lambda\right|=\left\|\Sigmab_k\left(\bfB_k\right)-\widetilde{\Sigmab}_k\left(\bfB_k\right)\right\|_2 \le \epsilon,$$ where $\text{vp}$ denotes the spectrum of matrix, $-\epsilon \ge -\lambda_m/3$, and $\alpha = 3 \epsilon /\lambda_m$, it follow that
	\begin{align*}
		&\xv^\top \left[\left(1+\alpha\right)\widetilde{\Sigmab}_k\left(\bfB_k\right) - \Sigmab_k\left(\bfB_k\right)\right]\xv\nn\\
		&=\left(1+\alpha\right)\xv^\top \left[\widetilde{\Sigmab}_k\left(\bfB_k\right) - \Sigmab_k\left(\bfB_k\right)\right]\xv+\alpha \xv^\top  \Sigmab_k\left(\bfB_k\right)\xv\\
		&\ge -\left(1+\alpha\right)\left\|\widetilde{\Sigmab}_k\left(\bfB_k\right) - \Sigmab_k\left(\bfB_k\right)\right\|_2\left\|\xv\right\|_2^2 + \alpha \lambda_m \left\|\xv\right\|_2^2\\ 
		&\ge \left(\alpha \lambda_m-\left(1+\alpha\right)\epsilon\right)\left\|\xv\right\|_2^2 =  \left(\alpha \lambda_m-\alpha\epsilon-\epsilon\right)\left\|\xv\right\|_2^2\\
		& \ge \left(\frac{2}{3}\alpha \lambda_m-\epsilon\right)\left\|\xv\right\|_2^2= \epsilon\left\|\xv\right\|_2^2>0,
	\end{align*}
	and 
	\begin{align*}
		&\xv^\top \left[\Sigmab_k\left(\bfB_k\right) - \left(1+\alpha\right)^{-1}\widetilde{\Sigmab}_k\left(\bfB_k\right)\right]\xv\nn\\
		&=\left(1+\alpha\right)^{-1}\xv^\top \left[ \Sigmab_k\left(\bfB_k\right)-\widetilde{\Sigmab}_k\left(\bfB_k\right)\right]\xv+\left(1-\left(1+\alpha\right)^{-1}\right) \xv^\top  \Sigmab_k\left(\bfB_k\right)\xv\\
		&\ge \left(\frac{\alpha \lambda_m-\epsilon}{1+\alpha}\right)\left\|\xv\right\|_2^2
		= \frac{2\epsilon}{1+\alpha}\left\|\xv\right\|_2^2
		\ge \epsilon\left\|\xv\right\|_2^2>0 \left(\text{ since }0 < \alpha < 1\right).
	\end{align*}
\end{proof}
By \cref{lem.ratioGaussian_intro} and the same argument as in the proof of Lemma B.9 from \cite{maugis2011non}, given  $0 < \epsilon < \lambda_m/3$, where $\epsilon$ is chosen later, and $\alpha = 3 \epsilon /\lambda_m$, we obtain
\begin{align}
	\max\left\{\frac{l(\xv,\yv)}{f(\xv,\yv)},\frac{f(\xv,\yv)}{u(\xv,\yv)}\right\} \le \left(1+2\alpha\right)^{-\frac{D}{2}} \exp\left(\frac{\left\|\upsilonb_{k,d}(\yv)-\widetilde{\upsilonb}_{k,d}(\yv)\right\|_2^2}{2\epsilon}\right).\label{eq.conditionRadiusDiscVariance_intro}
\end{align}
Because $\ln\left(\cdot\right)$ is a non-decreasing function, $\ln\left(1+2\alpha\right) \ge \alpha, \forall \alpha \in \left[0,1\right]$. Combined with \eqref{eq.netMeanGaussian_intro} where $\delta^2_{\Upsilondk} = D \alpha \epsilon$, we conclude that
\begin{align*}
	\max\left\{\ln \left(\frac{l(\xv,\yv)}{f(\xv,\yv)}\right),\ln\left(\frac{f(\xv,\yv)}{u(\xv,\yv)}\right)\right\} &\le -\frac{D}{2} \ln \left(1+2\alpha\right)+ \frac{\delta^2_{\Upsilondk}}{2\epsilon}\le -\frac{D}{2} \alpha+ \frac{\delta^2_{\Upsilondk}}{2\epsilon}=0.
\end{align*}
This means that $l(\xv,\yv) \le f(\xv,\yv) \le u(\xv,\yv), \forall (\xv,\yv) \in \cX \times \cY$. 
Hence, it remains to bound the size of bracket $\left[l,u\right]$ \wrt~$d_{\gdbk}$.
To this end, we aim to verify that $d^2_{\gdbk}\left(l,u\right) \le \frac{\delta}{2}$.
To do that, we make use of the following \cref{lem.Hellinger2Gaussian_intro}.

\begin{lem}[Proposition C.3 from \cite{maugis2011non}]\label{lem.Hellinger2Gaussian_intro}
	Let $\phi\left(\cdot;\mub_1,\Sigmab_1\right)$ and $\phi\left(\cdot;\mub_2,\Sigmab_{2}\right)$ be two Gaussian densities with full rank covariance. It holds that
	\begin{align*}
		&d^2\left(\phi\left(\cdot;\mub_1,\Sigmab_1\right),\phi\left(\cdot;\mub_2,\Sigmab_{2}\right)\right) \\
		&= 2\left\{1-2^{D/2}\left|\Sigmab_1\Sigmab_{2}\right|^{-1/4}\left|\Sigmab_1^{-1}+\Sigmab_{2}^{-1}\right|^{-1/2}
		\exp^{-\frac{1}{4}\left(\mub_1-\mub_2\right)^\top \left(\Sigmab_1+\Sigmab_{2}\right)^{-1}\left(\mub_1-\mub_2\right)}\right\}.
	\end{align*}
\end{lem}
Therefore, using the fact that $\cosh(t) = \frac{e^{-t}+e^{t}}{2}$, \cref{lem.Hellinger2Gaussian_intro} leads to, for all $\yv \in \cY$:
\begin{align*}
	\hel(l(\cdot,\yv),u(\cdot,\yv)) 
	&
	= \int_{\cX}\left[ l(\xv,\yv) + u(\xv,\yv) - 2\sqrt{l(\xv,\yv)u(\xv,\yv)} \right]d\xv \\
	%
	%
	&= \left(1+2\alpha\right)^{-D}+ \left(1+2\alpha\right)^{D}-2\\
	&+\hel\left(\phi\left(\cdot;\widetilde{\upsilonb}_{k,d}(\yv),\left(1+\alpha\right)^{-1}\widetilde{\Sigmab}_k\left(\bfB_k\right)\right),\phi\left(\cdot;\widetilde{\upsilonb}_{k,d}(\yv),\left(1+\alpha\right)\widetilde{\Sigmab}_k\left(\bfB_k\right)\right)\right)\\
	&= 2\cosh\left[D\ln\left(1+2\alpha\right)\right]-2\\
	&\quad+2\left[1-2^{D/2}\left[\left(1+\alpha\right)^{-1}+\left(1+\alpha\right)\right]^{-D/2}\left|\widetilde{\Sigmab}_k\left(\bfB_k\right)\right|^{-1/2}\left|\widetilde{\Sigmab}_k\left(\bfB_k\right)\right|^{1/2}\right] 
	\\
	&= 2\cosh\left[D\ln\left(1+2\alpha\right)\right]-2+2-2\left[\cosh\left(\ln\left(1+\alpha\right)\right)\right]^{-D/2}\\
	& = 2g\left(D\ln\left(1+2\alpha\right)\right) + 2 h\left(\ln\left(1+\alpha\right)\right),
\end{align*}
where $g(t)=\cosh(t)-1 = \frac{e^{-t}+e^{t}}{2}  - 1$, and $h(t) = 1- \cosh(t)^{-D/2}$.
The upper bounds of terms $g$ and $h$ separately imply that, for all $\yv \in \cY$,
\begin{align*}
	\hel(l(\cdot,\yv),u(\cdot,\yv)) \le 2\left(2\cosh\left(\frac{1}{\sqrt{6}}\right)\alpha^2D^2 + \frac{1}{4}\alpha^2D^2\right) \le 6 \alpha^2 D^2 = \frac{\delta^2}{4},
\end{align*} 
where we choose $\alpha = \frac{3 \epsilon}{\lambda_m}, \epsilon = \frac{\delta \lambda_m}{6\sqrt{6}D}$, $\forall \delta \in (0,1],D \in \Ns, \lambda_m > 0$, which appears in \eqref{eq.conditionRadiusDiscVariance_intro} and satisfies $\alpha=  \frac{\delta}{2\sqrt{6}D}$ and $0 < \epsilon < \frac{\lambda_m}{3}$. Indeed, studying functions $g$ and $h$ yields
\begin{align*}
	g'(t) &= \sinh(t), g''(t) = \cosh(t)\le  \cosh(c), \forall t \in [0,c], c \in\R_+,
	\\
	h'(t) &=\frac{D}{2} \cosh(t)^{-D/2-1} \sinh(t),\\ 
	h''(t) &=\frac{D}{2} \left(-\frac{D}{2}-1\right) \cosh(t)^{-D/2-2} \sinh^2(t) +\frac{D}{2} \cosh(t)^{-D/2}\\
	&=\frac{D}{2}\left(1-\left(\frac{D}{2}+1\right)\left(\frac{\sinh(t)}{\cosh(t)}\right)^2\right)  \cosh(t)^{-D/2} \le \frac{D}{2},
\end{align*}
where we used the fact that $\cosh(t) \ge 1$. Then, since $g(0) = 0,g'(0) =0,h(0) = 0,h'(0) = 0$, by applying Taylor's Theorem, it is true that 
\begin{align*}
	g(t) &= g(t) - g(0) - g'(0)t =  R_{0,1}(t) \le \cosh(c)\frac{t^2}{2},\forall t \in [0,c],\\
	h(t) &= h(t) - h(0) - h'(0)t =  R_{0,1}(t) \le \frac{D}{2}\frac{t^2}{2} \le \frac{D^2}{2}\frac{t^2}{2},\forall t \ge 0.
\end{align*}
We wish to find an upper bound for $t = D \ln\left(1+2 \alpha\right)$, $D \in \Ns$, $\alpha = \frac{\delta}{2\sqrt{6}D}$, $\delta \in (0,1]$. Since $\ln$ is an increasing function, then we have 
$$t = D \ln\left(1+\frac{\delta}{\sqrt{6}D}\right) \le D \ln\left(1+\frac{1}{\sqrt{6}D}\right) \le D \frac{1}{\sqrt{6}D}= \frac{1}{\sqrt{6}}, \forall\delta \in (0,1],$$ since $\ln\left(1+\frac{1}{\sqrt{6}D}\right) \le \frac{1}{\sqrt{6}D}$, $\forall D \in \Ns$. Then, since $\ln\left(1+ 2\alpha\right) \le 2\alpha, \forall \alpha \ge 0$, 
\begin{align*}
	g\left(D \ln\left(1+2 \alpha\right)\right) &\le \cosh\left(\frac{1}{\sqrt{6}}\right)\frac{\left(D \ln\left(1+2 \alpha\right)\right)^2}{2} \le \cosh\left(\frac{1}{\sqrt{6}}\right) \frac{D^2}{2} 4 \alpha^2,\\
	h\left(\ln\left(1+\alpha\right)\right) &\le \frac{D^2}{2}\frac{\left(\ln\left(1+\alpha\right)\right)^2}{2} \le \frac{D^2\alpha^2}{4}.
\end{align*}


%
Note that the set of $\delta/2$-brackets $[l,u]$ over $\gdbk$ is totally defined by the parameter spaces ${\tilde{S}}^{\disc}_{\bfB_k}(\epsilon)$ and $G_{\Upsilondk}\left(\delta_{\Upsilondk}\right)$. This leads to an upper bound of the $\delta/2$-bracketing entropy of $\gdbk$ evaluated from an upper bound of the two set cardinalities. Hence, given any $\delta > 0$, by choosing $\epsilon = \frac{\delta \lambda_m}{6\sqrt{6}D}$, $\alpha = \frac{3 \epsilon}{\lambda_m}=\frac{\delta}{2\sqrt{6}D}$, and $\delta^2_{\Upsilondk} = D \alpha \epsilon=D\frac{\delta}{2\sqrt{6}D}\frac{\delta \lambda_m}{6\sqrt{6}D} = \frac{\delta^2\lambda_m}{72D}$, it holds that
\begin{align*}
	&\cN_{[\cdot],d_{\gdbk} }\left(\frac{\delta}{2},\gdbk\right)
	\nn\\
	&\le \card\left({\tilde{S}}^{\disc}_{\bfB_k}(\epsilon)\right)\times \card\left(G_{\Upsilondk}\left(\delta_{\Upsilondk}\right)\right) \\
	&\le \left(\Bigg\lfloor \frac{2\lambda_M}{\epsilon}\Bigg\rfloor\frac{D\left(D-1\right)}{2 D_{\bfB_k}}\right)^{D_{\bfB_k}} \left(\frac{\exp\left(C_{\Upsilondk}\right)}{\delta_{\Upsilondk}}\right)^{\dim\left(\Upsilondk\right)} \left(\text{using \eqref{eq.upperCardCovariance_intro} and \eqref{eq.cardMeanGaussianEx_intro}}\right)\\
	&\le \left(\frac{2\lambda_M 6 \sqrt{6}D}{\delta \lambda_m}\frac{D\left(D-1\right)}{2 D_{\bfB_k}}\right)^{D_{\bfB_k}} \left(\frac{6 \sqrt{2D}\exp\left(C_{\Upsilondk}\right)}{\delta \sqrt{\lambda_m}}\right)^{\dim\left(\Upsilondk\right)}\\
	&= \left(\frac{6 \sqrt{6} \lambda_M D^2\left(D-1\right)}{ \lambda_m D_{\bfB_k}}\right)^{D_{\bfB_k}} \left(\frac{6 \sqrt{2D}\exp\left(C_{\Upsilondk}\right)}{ \sqrt{\lambda_m}}\right)^{\dim\left(\Upsilondk\right)}\left(\frac{1}{\delta}\right)^{D_{\bfB_k}+\dim\left(\Upsilondk\right)}.
\end{align*}
Finally, by definition of bracketing entropy in \eqref{eq_define_BracketingEntropy_intro}, we obtain
\begin{align*}
	&\cH_{[\cdot],d_{\gdbk} }\left(\frac{\delta}{2},\gdbk\right)
	\nn\\
	%
	&\le D_{\bfB_k} \ln \left(\frac{6 \sqrt{6} \lambda_M D^2\left(D-1\right)}{ \lambda_m D_{\bfB_k}}\right) + \dim\left(\Upsilondk\right) \ln\left(\frac{6 \sqrt{2D}\exp\left(C_{\Upsilondk}\right)}{ \sqrt{\lambda_m}}\right)\nn\\
	&\quad +\left(D_{\bfB_k}+\dim\left(\Upsilondk\right)\right)\ln\left(\frac{1}{\delta}\right)
	= \dim\left(\gdbk\right) \left(C_{\gdbk} + \ln \left(\frac{1}{\delta}\right)\right),
\end{align*}
where $\dim\left(\gdbk\right) = D_{\bfB_k}+\dim\left(\Upsilondk\right)$ and $$C_{\gdbk} = \frac{D_{\bfB_k} \ln \left(\frac{6 \sqrt{6} \lambda_M D^2\left(D-1\right)}{ \lambda_m D_{\bfB_k}}\right) + \dim\left(\Upsilondk\right) \ln\left(\frac{6 \sqrt{2D}\exp\left(C_{\Upsilondk}\right)}{ \sqrt{\lambda_m}}\right)}{\dim\left(\gdbk\right)}.$$


\bibliographystyle{imsart-number} 

\bibliography{GLoME}       


\end{document}